\documentclass[article]{article}
\usepackage{amssymb}
\usepackage{amsmath}
\usepackage{amsthm}
\usepackage[affil-it]{authblk}
\usepackage{booktabs}
\usepackage[usenames]{color}
\usepackage[hmargin={25mm,25mm},vmargin={30mm,35mm}]{geometry}
\usepackage[utf8]{inputenc}
\usepackage[colorlinks,allcolors={blue}]{hyperref}
\usepackage{tabularx}
\usepackage{xcolor}
\usepackage{multirow}
\usepackage{colortbl}
\usepackage{array}
\usepackage{algorithm2e}
\usepackage[noend]{algorithmic}
\usepackage{paralist}
\usepackage{stackrel}
\usepackage{cancel}
\usepackage{graphicx}
\usepackage{rotating}

\usepackage{newtxtext}
\usepackage{newtxmath}

\definecolor{LightCyan}{rgb}{0.88,1,1}
\definecolor{LightGoldenrod}{rgb}{0.93,0.87,0.51}
\definecolor{LightBlue}{rgb}{0.54, 0.81, 0.94}
\definecolor{BallBlue}{rgb}{0.13, 0.67, 0.8}
\definecolor{DardSeeGreen}{rgb}{0.56, 0.74, 0.56}

\newtheorem{theorem}{Theorem}

\theoremstyle{remark}
\newtheorem{remark}[theorem]{Remark}
\theoremstyle{definition}

\newcommand{\Real}{\mathbb{R}}

\newcommand{\eg}{e.g.}

\DeclareMathOperator{\DIV}{div}

\newcommand{\aT}{T}
\newcommand{\kl}{k_\ell}
\newcommand{\klpo}{k_{\ell+1}}
\newcommand{\IntOp}{\mathcal{I}}
\newcommand{\Mh}{\mathcal{M}_h}
\newcommand{\Th}{\mathcal{T}_h}
\newcommand{\Fh}{\mathcal{F}_h}

\newcommand{\Poly}[1]{\mathcal{P}^{#1}}
\newcommand{\Polyd}[2]{\mathbb{P}_{#1}^{#2}}
\DeclareMathOperator{\CARD}{card}
\newcommand{\card}[1]{\CARD(#1)}

\newcommand{\jump}[1]{\llbracket #1 \rrbracket}
\newcommand{\averg}[1]{\{ #1 \}}

\newcommand{\normal}{\boldsymbol{n}}

\newcommand{\FT}{\mathcal{F}_T}
\newcommand{\TF}{\mathcal{T}_F}
\newcommand{\FTp}{\mathcal{F}_{T'}}

\newcommand{\A}{\boldsymbol{A}}
\newcommand{\IntOpB}{\boldsymbol{\mathcal{I}}}
\newcommand{\bVec}{\boldsymbol{b}}
\newcommand{\uVec}{\boldsymbol{u}}

\newcommand{\wVec}{\boldsymbol{w}}
\newcommand{\zVec}{\boldsymbol{z}}
\newcommand{\vVec}{\boldsymbol{v}}

\newcommand{\fVec}{\boldsymbol{f}}

\everymath{\displaystyle}

\newcommand{\dir}{{\rm D}}
\newcommand{\neu}{{\rm N}}
\newcommand{\internal}{{\rm i}}

\newcommand{\momentum}{{\rm mnt}}
\newcommand{\mass}{{\rm cnt}}

\newtheorem{scheme}{Scheme}

\newcommand{\MAT}[1]{\boldsymbol{\sf #1}}
\newcommand{\VEC}[1]{{\sf #1}}
\newcommand{\UVEC}[1]{\underline{\sf #1}}
\newcommand{\trans}{^\intercal}

\newcommand{\hho}{{\tt HHO}}
\newcommand{\hhodp}{{\tt HHO-dp}}
\newcommand{\hhohp}{{\tt HHO-hp}}
\newcommand{\vcond}{{\tt v-cond}}
\newcommand{\vc}{{\tt v}}
\newcommand{\uncond}{{\tt uncond}}
\newcommand{\vpcond}{{\tt v\&p-cond}}
\newcommand{\vpc}{{\tt v\&p}}
\newcommand{\dg}{{\tt DG}}

\usepackage[normalem]{ulem}
\newcounter{corr}
\definecolor{violet}{rgb}{0.580,0.,0.827}
\newcommand{\corr}[3]{\typeout{Warning : a correction remains in page \thepage}
  \stepcounter{corr}        
	      {\color{blue}\ifmmode\text{\,\sout{\ensuremath{#1}}\,}\else\sout{#1}\fi}
              {\color{red}#2}
              {\color{violet} #3}
}

\newcommand{\email}[1]{\href{mailto:#1}{#1}}

\begin{document}

\title{$p$-Multilevel preconditioners for HHO discretizations of the Stokes equations with static condensation}
\author[1]{Lorenzo Botti}
\author[2]{Daniele A. Di Pietro}
\affil[1]{Department of Engineering and Applied Sciences, University of Bergamo, Italy, \email{lorenzo.botti@unibg.it}}
\affil[2]{IMAG, Univ Montpellier, CNRS, Montpellier, France, \email{daniele.di-pietro@umontpellier.fr}}

\maketitle

\begin{abstract}
  We propose a $p$-multilevel preconditioner for Hybrid High-Order discretizations (HHO) of the Stokes equation, numerically assess its performance on two variants of the method, and compare with a classical Discontinuous Galerkin scheme.
  We specifically investigate how the combination of $p$-coarsening and static condensation influences the performance of the $V$-cycle iteration for HHO.
  Two different static condensation procedures are considered, resulting in global linear systems with a different number of unknowns and non-zero elements.
  An efficient implementation is proposed where coarse level operators are inherited using $L^2$-orthogonal projections defined over mesh faces and the restriction of the fine grid operators is performed recursively and matrix-free.
  The various resolution strategies are thoroughly validated on two- and three-dimensional problems.
\end{abstract}



\section{Introduction}

In this work we develop and numerically validate $p$-multigrid solution strategies for nonconforming polytopal discretizations of the Stokes equations, governing the creeping flow of incompressible fluids.
\smallskip

For the sake of simplicity, we focus on a Newtonian fluid with uniform density and unit kinematic viscosity.
Given a polygonal or polyhedral domain $\Omega\subset\Real^d$, $d\in\{2,3\}$, with boundary $\partial \Omega$, the Stokes problem consist in finding the velocity field $\uVec: \Omega \rightarrow \Real^d$, and the pressure field $p: \Omega \rightarrow \Real$,
such that
\begin{subequations}
  \label{stokesProb}
  \begin{alignat}{2}\label{stokesProb:momentum}
    -{\Delta \uVec} + {\nabla p} &= \boldsymbol{f} &\qquad& \text{in $\Omega$}, \\
    \label{stokesProb:mass}
          {\nabla \cdot \uVec} &= 0 &\qquad& \text{in $\Omega$}, \\ 
    \uVec &= \boldsymbol{g}_\dir &\qquad& \text{on $\partial\Omega_\dir$}. \\
    - \normal \cdot \nabla \uVec + p \normal &= \boldsymbol{g}_\neu &\qquad& \text{on $\partial\Omega_\neu$},
  \end{alignat}
\end{subequations}
where $\normal$ denotes the unit vector normal to $\partial\Omega$ pointing out of $\Omega$,
$\boldsymbol{g}_\dir$ and $\boldsymbol{g}_\neu$ denote, respectively, the prescribed velocity on the Dirichlet boundary $\partial\Omega_\dir\subset\partial\Omega$ and the prescribed traction on the Neumann boundary $\partial\Omega_\neu\coloneq\partial\Omega\setminus\partial\Omega_\dir$,
while $\fVec:\Omega\to\Real^d$ is a given body force.
For the sake of simplicity, it is assumed in what follows that both $\partial\Omega_\dir$ and $\partial\Omega_\neu$ have non-zero $(d-1)$-dimensional Hausdorff measure (otherwise, additional closure conditions are needed).
\smallskip

Our focus is on new generation discretization methods for problem \eqref{stokesProb} that support general polytopal meshes and high-order: Hybrid High-Order (HHO) and Discontinuous Galerkin (DG) methods.

Hybrid High-Order discretizations of the Stokes equations have been originally considered in \cite{Aghili.Boyaval.ea:15} and later extended in \cite{Di-Pietro.Ern.ea:16} to incorporate robust handling of large irrotational body forces.
Other extensions include their application to the Brinkman problem, considered in \cite{Botti.Di-Pietro.ea:18}, and to the full Navier--Stokes equations \cite{Di-Pietro.Krell:17,Di-Pietro.Krell:18,Botti.Di-Pietro.ea:19*1}; see also \cite[Chapters 8 and 9]{Di-Pietro.Droniou:20} for further details.
In this work, we consider two HHO schemes that are novel variations of existing schemes with improved features.
The first scheme, based on a hybrid approximation of the velocity along with a discontinuous approximation of the pressure, is a variation of the one considered in \cite[Chapter 8]{Di-Pietro.Droniou:20} including two choices for the polynomial degree of the element velocity unknowns in the spirit of \cite{Cockburn.Di-Pietro.ea:16} (see also \cite[Section 5.1]{Di-Pietro.Droniou:20}).
The second scheme, inspired by the Hybridizable Discontinuous Galerkin (HDG) method of \cite{Rhebergen.Wells:18}, hinges on hybrid approximations of both the velocity and the pressure and includes, with respect to the above reference, a different treatment of viscous terms that results in improved orders of convergence.
In both cases, the Dirichlet condition on the velocity is enforced weakly in the spirit of \cite{Botti.Di-Pietro.ea:19*1}.

Since the pioneering works \cite{Cockburn.Shu:91,Cockburn.Shu:89,Cockburn.Lin.ea:89,Cockburn.Hou.ea:90,Cockburn.Shu:98} dating back to the late 1980s, DG methods have gained significant popularity in computational fluid mechanics, boosted by the 1997 landmark papers \cite{Bassi.Rebay:97,Bassi.Rebay.ea:97} on the treatment of viscous terms.
The extension of DG methods to general polyhedral meshes was systematically considered in \cite{Di-Pietro.Ern:10} and \cite{Di-Pietro.Ern:12}.
Crucially, this extension paved the way to adaptive mesh coarsening by agglomeration, a strategy proposed in \cite{Bassi.Botti.ea:12} and exploited in  \cite{Bassi.Botti.ea:12*1,Bassi.Botti.ea:14} in practical CFD applications to provide high-order accurate geometry representation with arbitrarily coarse meshes.
More recent developments, including $hp$-versions and the support of meshes with small faces, can be found in \cite{Antonietti.Giani.ea:13,Antonietti.Cangiani.ea:16}; see also the recent monograph \cite{Cangiani.Dong.ea:17}. 
Our focus is on an equal-order approximation with stabilized pressure-velocity coupling in the spirit of \cite{Cockburn.Kanschat.ea:02} and a treatment of the viscous term based on the Bassi--Rebay 2 (BR2) method of \cite{Bassi.Rebay.ea:97}.
Related works include \cite{Bassi.Crivellini.ea:06,Di-Pietro:07}; see also \cite[Chapter 6]{Di-Pietro.Ern:12} and references therein.
\smallskip

$p$-Multilevel solvers are well suited for both HHO and DG methods
because the process of building coarse level operators based on polynomial degree reduction is straightforward and inexpensive.
The purpose of applying iterative solvers to coarse problems is twofold: 
on one hand, a coarser operator translates into a global sparse matrix of smaller size with fewer non-zero entries, resulting in cheaper matrix-vector products; 
on the other hand, coarse level iterations are best suited to smooth out the low-frequency components of the error, that are hardly dumped by fine level iterations. 
In the context of DG discretizations, $p$-multilevel solvers have been fruitfully utilized in practical applications see, e.g., \cite{Fidkowski05,Nastase06,BassiGhidoni09,ShahbaziMavriplis-MGAlgorithms:2009,Franciolini2020}.
$h$-,$p$- and $hp$-Multigrid solvers for DG discretizations of elliptic problems have been considered in \cite{AntoniettiSartiVerani15}, where uniform convergence with respect to the number of levels for the W-cycle iteration has been proved, and in \cite{hphpBotti19}.
Multigrid solvers for HDG discretizations of scalar elliptic problems were considered in \cite{CockburnMulti13} and, more recently, in \cite{FabienHDGMultigrid19,KronbichlerMultiGhdgdg18}, where a comparison with DG is carried out.
$p$-Multivel solvers for HDG methods with application to compressible flow simulations have been recently considered in \cite{FrancioliniHDG20}.
Preconditioners for DG and HDG discretizations of the Stokes problem have been considered in \cite{DeDiosStokes12,KanschatStokes15,AdlerStokes17,BOTTIhMG17,CharrierStokesdG17} and \cite{RhebergenStokesP18}, respectively.
Finally, an $h$-multigrid method for HHO discretizations of scalar diffusion problems has been recently proposed in \cite{Matalon.Di-Pietro.ea:20}.
The main novelty consists, in this case, in the use of the local potential reconstruction in the prolongation operator.
\smallskip

In this work we propose and numerically assess $p$-multilevel solution strategies for HHO discretizations of the Stokes equations.
We specifically investigate how the combination of $p$-coarsening and static condensation influences the performance of the $V$-cycle iteration.
To this end, we compare different static condensation strategies.
In order to preserve computational efficiency, statically condensed coarse level operators are inherited 
using local $L^2$-orthogonal projections defined over mesh faces. 
Restriction of fine grid operators is performed recursively and matrix-free, relying on $L^2$-orthogonal basis functions  
to further reduce the computational burden. 
Performance assessment is based on accuracy and efficiency of $p$-multilevel solvers 
considering DG discretizations as a reference for comparison.
High-order accurate solutions approximating smooth analytical velocity and pressure fields 
are computed over standard and severely graded $h$-refined mesh sequences in both two and three space dimensions.
Interestingly, the static condensation strategy plays a crucial role in case of graded meshes. 
\smallskip

The rest of this work is organized as follows.
In Section \ref{sec:threeNCD} we state the HHO and DG schemes considered in the numerical tests.
The $p$-multilevel strategy is discussed in Section \ref{sec:p-multilevel.strategies} and computational aspects are discussed in Section \ref{sec:computational.aspects}.
Section \ref{sec:numResults} contains an extensive panel of numerical results that enable one to assess and compare several solution strategies.
Finally, some conclusions are drawn in Section \ref{sec:conclusion}.


\section{Three nonconforming methods for the Stokes problem}
\label{sec:threeNCD}

In this section we describe two HHO and one DG methods for the approximation of problem \eqref{stokesProb} that will be used to assess the performance of the $p$-multilevel preconditioner.
In order to lay the ground for future works on the full nonlinear Navier--Stokes equations, the corresponding discrete problems are formulated in terms of the annihilation of residuals.

\subsection{Discrete setting}

We consider meshes of the domain $\Omega$ corresponding to couples $\Mh\coloneq(\Th,\Fh)$, where $\Th$ is a finite collection of polygonal (if $d=2$) or polyhedral (if $d=3$) elements such that $h\coloneq\max_{T\in\Th}h_T>0$ with $h_T$ denoting the diameter of $T$, while $\Fh$ is a finite collection of line segments (if $d=2$) or polygonal faces (if $d=3$).
For the sake of brevity, in what follows the term ``face'' will be used in both two and three space dimensions.
It is assumed henceforth that the mesh $\Mh$ matches the geometrical requirements detailed in \cite[Definition 1.4]{Di-Pietro.Droniou:20}.
This covers, essentially, any reasonable partition of $\Omega$ into polyhedral sets, not necessarily convex.
For each mesh element $T \in \Th$, the faces contained in the element boundary $\partial T$ are collected in the set $\FT$,
and, for each mesh face $F \in \Fh$, $\TF$ is the set containing the one or two mesh elements sharing $F$.
We define three disjoint subsets of the set $\FT$:
the set of Dirichlet boundary faces $\FT^\dir \coloneq \{F \in \FT : F \subset \partial \Omega_\dir\}$;
the set of Neumann boundary faces $\FT^\neu \coloneq \{F \in \FT : F \subset \partial \Omega_\neu\}$;
the set of internal faces $\FT^\internal \coloneq \FT\setminus\big(\FT^\dir\cup\FT^\neu\big)$.
For future use, we also let $\FT^{\internal,\dir}\coloneq\FT^\internal\cup\FT^\dir$.
For all $T\in\Th$ and all $F\in\FT$, $\normal_{TF}$ denotes the unit vector normal to $F$ pointing out of $T$.

Hybrid High-Order methods hinge on local polynomial spaces on mesh elements and faces.
For given integers $\ell\ge 0$ and $n\ge 1$, we denote by $\Polyd{n}{\ell}$ the space of $n$-variate polynomials of total degree $\le\ell$ (in short, of degree $\ell$).
For $X$ mesh element or face, we denote by $\Poly{\ell}(X)$ the space spanned by the restriction to $X$ of functions in $\mathbb{P}_d^\ell$.
When $X$ is a mesh face, the resulting space is isomorphic to $\mathbb{P}_{d-1}^\ell$ (see \cite[Proposition 1.23]{Di-Pietro.Droniou:20}).
At the global level, we will need the broken polynomial space
\[
\Poly{\ell}(\Th)\coloneq\left\{
q\in L^2(\Omega) : \text{$q_{|T}\in\Poly{\ell}(T)$ for all $T\in\Th$}
\right\}.
\]

Let again $X$ denote a mesh element or face.
The local $L^2$-orthogonal projector $\pi_X^\ell:L^2(X)\to \Poly{\ell}(X)$ is such that, for all $q\in L^2(X)$,
\[
\int_X ( q - \pi_X^\ell q ) r = 0\qquad\forall r\in\Poly{\ell}(X).
\]
Notice that, above and in what follows, we omit the measure from integrals as it can always be inferred from the context.
The $L^2$-orthogonal projector on $\Poly{\ell}(X)^d$, obtained applying $\pi_X^\ell$ component-wise, is denoted by $\boldsymbol{\pi}_X^\ell$.

\subsection{Local reconstructions and face residuals}

The HHO discretizations of the Stokes problem considered in this work hinge on velocity reconstructions devised at the element level and obtained assembling diffusive potential reconstructions component-wise.
In what follows, we let a mesh element $T\in\Th$ be fixed, denote by $k\ge 0$ the degree of polynomials attached to mesh faces, and by $k'\in\{k,k+1\}$ the degree of polynomials attached to mesh elements.

\subsubsection{Scalar potential reconstruction}

The velocity reconstruction is obtained leveraging, for each component, the \emph{scalar potential reconstruction} originally introduced in \cite{Di-Pietro.Ern.ea:14} in the context of scalar diffusion problems (see also \cite{Cockburn.Di-Pietro.ea:16} and \cite[Section 5.1]{Di-Pietro.Droniou:20} for its generalization to the case of different polynomial degrees on elements and faces).
Define the local scalar HHO space
\begin{equation}\label{eq:VTk}
  \underline{V}_T^{k',k} \coloneq \left\{ \underline{v}_T = \left(v_T, (v_F)_{F \in \mathcal{F}_T} \right) : 
  \text{$v_T \in \Poly{k'}(T)$ and $v_F \in \Poly{k}(F)$ for all $F \in\mathcal{F}_T$}
  \right\}.  
\end{equation}
The scalar potential reconstruction operator $\mathfrak{p}_T^{k+1}$: $\underline{V}_T^{k',k} \rightarrow \Poly{k+1}(T)$ maps a vector of polynomials of $\underline{V}_T^{k',k}$ onto a polynomial of degree $(k+1)$ over $T$ as follows:
Given $\underline v_{T} \in \underline{V}_T^{k',k}$, $\mathfrak{p}_T^{k+1}\underline{v}_T$ is the unique polynomial in $\Poly{k+1}(T)$ satisfying
\begin{equation*} 
  \begin{alignedat}{2}
    \int_T \nabla \mathfrak{p}_T^{k+1} \underline{v}_T \cdot \nabla {w_T}
    &=  \int_T \nabla {v_T} \cdot \nabla {w_T} + \sum_{F \in \FT} \int_F \left({v_F} - {v_T}\right) \, \nabla {w_T} \cdot \normal_{TF}
    &\qquad&  \forall w_T \in \Poly{k +1}(T),
    \\
    \int_T \mathfrak{p}_T^{k+1} \underline{v}_T &= \int_T {v_T}.
  \end{alignedat}
\end{equation*}
Computing $\mathfrak{p}_T^{k+1}$ for each $T\in\Th$ requires to solve a small linear system.
This is an embarrassingly parallel task that can fully benefit from parallel architectures.

\subsubsection{Velocity reconstruction}
\label{sec::HHOVelRec}

Define, in analogy with \eqref{eq:VTk}, the following vector-valued HHO space for the velocity:
\[
  \underline{\boldsymbol{V}}_T^{k',k} \coloneq \left\{ \underline{\vVec}_T = \big(\vVec_T, (\vVec_F)_{F \in \mathcal{F}_T} \big) : 
  \text{%
    ${\vVec}_T \in \Poly{k'}(T)^d$ and $\vVec_F  \in \Poly{k}(F)^d$ for all $F \in\mathcal{F}_T$
  }
  \right\}.
\]
The \emph{velocity reconstruction} $\mathfrak{P}_T^{k+1}$: $\underline{\boldsymbol{V}}_T^{k',k} \rightarrow \Poly{k+1}(T)^d$ is obtained setting
\begin{equation*} 
  \mathfrak{P}_T^{k+1} \underline{\vVec}_T \coloneq \big(\mathfrak{p}_T^{k+1} \underline{v}_{T,i}\big)_{i=1,\ldots,d},
\end{equation*}
where, for all $i=1,\ldots,d$, $\underline{v}_{T,i}\in\underline{V}_T^{k',k}$ is obtained gathering the $i$th components of the polynomials in $\underline{\vVec}_T$, i.e., $\underline{v}_{T,i}\coloneq \big(v_{T,i}, (v_{F,i})_{F\in\FT}\big)$ if $\vVec_T=(v_{T,i})_{i=1,\ldots,d}$ and $\vVec_F=(v_{F,i})_{i=1,\ldots,d}$ for all $F\in\FT$.

\subsubsection{Face residuals}

Let $T\in\Th$ and $F\in\FT$.
The stabilization bilinear form for the HHO discretization of the viscous term in the momentum equation \eqref{stokesProb:momentum} hinges on the \emph{face residual} $\mathfrak{R}_{TF}^k : \underline{\boldsymbol{V}}_T^{k',k} \rightarrow \Poly{\max(k',k)}(F)^d$ such that, for all $\underline{\vVec}_T\in\underline{\boldsymbol{V}}_T^{k',k}$,
\[
\mathfrak{R}_T^{k',k} \underline{\vVec}_T \coloneq \big( \mathfrak{r}_{TF}^{k',k} \underline{v}_{T,i} \big)_{i=1,\ldots,d},
\]
where the scalar face residual $\mathfrak{r}_{TF}^{k',k}:\underline{V}_T^{k',k}\to\Poly{\max(k',k)}(F)$ is such that, for all $\underline{v}_T\in\underline{V}_T^{k',k}$,
\[
\mathfrak{r}_{TF}^{k',k} \underline{v}_T
\coloneq \pi^k_F \big(v_F - \mathfrak{p}_T^{k+1} \underline{v}_T \big) - \pi^{k'}_T \big(v_T - \mathfrak{p}_T^{k+1} \underline{v}_T \big).
\]

\subsection{HHO schemes}\label{sec:HHO.schemes}

We consider two HHO schemes based, respectively, on discontinuous and hybrid approximations of the pressure.
In both cases, the Dirichlet boundary condition is enforced weakly, considering a symmetric variation of the method discussed in \cite{Botti.Di-Pietro.ea:18}.

\subsubsection{An HHO scheme with discontinuous pressure}

Let again $k\ge 0$ and $k'\in\{k,k+1\}$ denote the polynomial degrees of the face and element unknowns, respectively, and let a mesh element $T\in\Th$ be fixed.
Given $(\underline \uVec_T, p_T) \in \underline{\boldsymbol{V}}_T^{k',k}\times \Poly{k}(T)$, the local residuals
$r^\momentum_{I,T}((\underline{\boldsymbol{u}}_T,p_T);\cdot):\underline{\boldsymbol{V}}_T^{k',k}\to\Real$ of the discrete momentum conservation equation and
$r^\mass_{I,T}(\underline{\boldsymbol{u}}_T;\cdot):\Poly{k}(T)\to\Real$ of the discrete mass conservation equation are such that, respectively:
For all $\underline{\vVec}_T\in\underline{\boldsymbol{V}}_T^{k',k}$ and all $q_T\in\Poly{k}(T)$,
\begin{subequations}\label{HHOdiscr:residuals}
  \begin{alignat}{2} \label{HHOdiscrMNT}
    r^\momentum_{I,T}((\underline{\boldsymbol{u}}_T,p_T);\underline{\boldsymbol{v}}_T)
    &\coloneq
    \begin{aligned}[t]
      &\int_T \nabla \mathfrak{P}_T^{k+1} \underline{\uVec}_{T} : \nabla \mathfrak{P}_T^{k+1} \underline{\vVec}_{T}
      + \sum_{F \in \FT} \frac{1}{h_F} \int_F \mathfrak{R}_{TF}^k \underline{\uVec}_T \, \cdot \, \mathfrak{R}_{TF}^k \underline{\vVec}_T
      \\
      & - \sum_{F \in \FT^\dir} \int_F \left[
        \bigl( \normal_{TF} \cdot \nabla \mathfrak{P}_T^{k+1} \underline{\uVec}_{T} \bigr) \cdot \vVec_F   
        + \uVec_F \cdot \bigl( \normal_{TF} \cdot \nabla \mathfrak{P}_T^{k+1} \underline{\vVec}_{T}\bigr)
        \right]
      +\sum_{F \in \FT^\dir} \frac{\eta}{h_F} \int_F \uVec_F \cdot \vVec_{F}   
      \\
      & - \int_T p_T \, (\nabla {\cdot} \vVec_T) - \sum_{F \in \FT} \int_F p_T \, (\vVec_F - \vVec_T) \cdot \normal_{TF}       
      + \sum_{F \in \FT^\dir} \int_F p_T \, (\vVec_F \cdot \normal_{TF})
      \\
      &
      - \sum_{F \in \FT^\dir} \int_F
      \boldsymbol{g}_\dir\cdot \left(
      \normal_{TF} \cdot \nabla \mathfrak{P}_T^{k+1} \underline{\vVec}_{T}
      + \frac{\eta}{h_F} \vVec_{F}
      \right)
      - \sum_{F \in \FT^\neu} \int_F \boldsymbol{g}_\neu \cdot \vVec_F
      - \int_{\aT} \boldsymbol{f} \cdot \vVec_T,
    \end{aligned}
    \\
    \label{HHOdiscrCNT} 
    r^\mass_{I,T}(\underline{\boldsymbol{u}}_T;q_T)
    &\coloneq
    \begin{aligned}[t]
      &- \int_T (\nabla {\cdot} \uVec_T) \, q_T
      - \sum_{F \in \FT} \int_F (\uVec_F - \uVec_T) \cdot \normal_{TF} \, q_T
      + \sum_{F \in \FT^\dir} \int_F (\uVec_F \cdot \normal_{TF}) \, q_T
      \\
      &- \sum_{F \in \FT^\dir} \int_F \boldsymbol{g}_\dir \cdot \normal_{TF} \, q_T.
    \end{aligned}
  \end{alignat}
\end{subequations}
In the expression of $r^\momentum_{I,T}((\underline{\boldsymbol{u}}_T,p_T);\cdot)$, $\eta>0$ is a user-dependent parameter that has to be taken large enough to ensure coercivity.
The penalty term where the parameter $\eta$ appears, along with the consistency terms in the second line and the term involving the boundary datum $\boldsymbol{g}_\dir$ in the fourth line, are responsible for the weak enforcement of the Dirichlet boundary condition for the velocity.

Define the global vector HHO space
\[
\underline{\boldsymbol{V}}_h^{k',k}\coloneq \left\{ \underline{\vVec}_h = \big( (\vVec_T)_{T\in\Th} , (\vVec_F)_{F \in \Fh} \big) : 
\text{%
  ${\vVec}_T \in \Poly{k'}(T)^d$ for all $T\in\Th$ and $\vVec_F  \in \Poly{k}(F)^d$ for all $F \in \Fh$
}
\right\}.
\]
For all $\underline{\vVec}_h\in\underline{\boldsymbol{V}}_h^{k',k}$ and all $T\in\Th$, we denote by $\underline{\vVec}_T\in\underline{\boldsymbol{V}}_T^{k',k}$ the restriction of $\underline{\vVec}_h$ to $T$.
The global residuals
$r_{I,h}^\momentum\left((\underline{\boldsymbol{u}}_h,p_h);\cdot\right):\underline{\boldsymbol{V}}_h^{k',k}\to\Real$
and $r^\mass_{I,h}(\underline{\boldsymbol{u}}_h;\cdot):\Poly{k}(\Th)\to\Real$
are obtained by element-by-element assembly, i.e.:
For all $\underline{\boldsymbol{v}}_h\in\underline{\boldsymbol{V}}_h^{k',k}$ and all $q_h\in\Poly{k}(\Th)$,
\begin{equation}
\label{HHOdisrcGlobalRes}
r_{I,h}^\momentum\left((\underline{\boldsymbol{u}}_h,p_h);\underline{\boldsymbol{v}}_h\right)
\coloneq \sum_{T\in\Th} r^\momentum_{I,T}\left((\underline{\boldsymbol{u}}_T,p_T);\underline{\boldsymbol{v}}_T\right),\qquad
r_{I,h}^\mass(\underline{\boldsymbol{u}}_h;q_h) \coloneq \sum_{T\in\Th} r^\mass_{I,T}(\underline{\boldsymbol{u}}_T;q_{h|T}).
\end{equation}
\begin{scheme}[\hhodp: HHO scheme with discontinuous pressure]\label{scheme:hho}
  Find $(\underline{\boldsymbol{u}}_h,p_h)\in\underline{\boldsymbol{V}}_h^{k',k}\times\Poly{k}(\Th)$ such that
  \begin{equation}\label{eq:hho}
    \begin{alignedat}{2}
      r_{I,h}^\momentum\left((\underline{\boldsymbol{u}}_h,p_h);\underline{\boldsymbol{v}}_h\right) &= 0 &\qquad& \forall \underline{\boldsymbol{v}}_h\in\underline{\boldsymbol{V}}_h^{k',k},
      \\
      r_{I,h}^\mass(\underline{\boldsymbol{u}}_h;q_h) &= 0 &\qquad& \forall q_h\in\Poly{k}(\Th).
    \end{alignedat}
  \end{equation}
\end{scheme}

\subsubsection{An HHO scheme with hybrid pressure}

An interesting variation of Scheme \ref{scheme:hho} is obtained combining the HHO discretization of the viscous term with $k'=k+1$ with a hybrid approximation of the pressure inspired by \cite{Rhebergen.Wells:18}.
Let $T\in\Th$.
Given $(\underline \uVec_T, \underline{p}_T)\in \underline{\boldsymbol{V}}_T^{k+1,k}\times \underline{V}_T^{k,k}$, the local residuals
$r^\momentum_{II,T} ((\underline{\boldsymbol{u}}_T,\underline{p}_T);\cdot):\underline{\boldsymbol{V}}_T^{k+1,k}\to\Real$ of the discrete momentum
and $ r^\mass_{II,T}(\underline{\boldsymbol{u}}_T;\cdot):\underline{V}^{k,k}_T\to\Real$ of the discrete mass conservation equations for the HHO scheme with hybrid pressure are such that, for all $\underline \vVec_T\in\underline{\boldsymbol{V}}_T^{k+1,k}$ and all $\underline{q}_T\in\underline{V}_T^{k,k}$,
\begin{subequations}
\begin{align*} 
  r^\momentum_{II,T} ((\underline{\boldsymbol{u}}_T,\underline{p}_T);\underline{\boldsymbol{v}}_T)
  &\coloneq
  \begin{aligned}[t]
    &\int_T \nabla \mathfrak{P}_T^{k+1} \underline{\uVec}_{T} : \nabla \mathfrak{P}_T^{k+1} \underline{\vVec}_{T}
    + \sum_{F \in \FT} \frac{1}{h_F} \int_F \mathfrak{R}_{TF}^k \underline{\vVec}_T \, \cdot \, \mathfrak{R}_{TF}^k \underline{\vVec}_T
    \\
        & - \sum_{F \in \FT^\dir} \int_F \left[
      \bigl( \normal_{TF} \cdot \nabla \mathfrak{P}_T^{k+1} \underline{\uVec}_{T} \bigr) \cdot \vVec_F   
      + \uVec_F \cdot \bigl( \normal_{TF} \cdot \nabla \mathfrak{P}_T^{k+1} \underline{\vVec}_{T}\bigr)
      \right]
    +\sum_{F \in \FT^\dir} \frac{\eta}{h_F} \int_F \uVec_F \cdot \vVec_{F}   
    \\
    &
    - \int_T p_{T} \, (\nabla \cdot \vVec_T)
    + \boxed{%
      \sum_{F \in \FT} \int_F p_F \, (\vVec_T - \vVec_F) \cdot \normal_{TF}
    }
    + \sum_{F \in \FT^\dir} \int_F p_F \, (\vVec_F \cdot \normal_{TF})
\\
&
    - \sum_{F \in \FT^\dir} \int_F
    \boldsymbol{g}_\dir\cdot \left(
      \normal_{TF} \cdot \nabla \mathfrak{P}_T^{k+1} \underline{\vVec}_{T}
      + \frac{\eta}{h_F} \vVec_{F}
      \right)
- \sum_{F \in \FT^\neu} \int_F \boldsymbol{g}_\neu \cdot \vVec_F
    - \int_{\aT} \boldsymbol{f} \cdot \vVec_T,
  \end{aligned}
  \\ 
  r^\mass_{II,T}(\underline{\boldsymbol{u}}_T;\underline{q}_T)
  &\coloneq
  - \int_T  (\nabla\cdot{\uVec}_T) \, q_T
  + \boxed{  
    \sum_{F \in \FT} \int_F ( \uVec_T - \uVec_F ) \cdot \normal_{TF} \, q_F .
  }
\end{align*}
\end{subequations}
As before, $\eta>0$ is a penalty parameter that has to be taken large enough to ensure coercivity.
The boxed terms are the ones that distinguish the local residuals on the momentum and mass conservation equations for the HHO scheme with hybrid pressure from Scheme \ref{scheme:hho} with $k'=k+1$.

Define the global scalar HHO space
\[
\underline{V}_h^{k,k}\coloneq\left\{
\underline{q}_h=\big((q_T)_{T\in\Th}, (q_F)_{F\in\Fh}\big)\,:\,
\text{$q_T\in\Poly{k}(T)$ for all $T\in\Th$ and $q_F\in\Poly{k}(F)$ for all $F\in\Fh$}
\right\}.
\]
The global residuals $r^\momentum_{II,h}((\underline{\boldsymbol{u}}_h,\underline{p}_h);\cdot):\underline{\boldsymbol{V}}^{k+1,k}_h\to\Real$ and
$r^\mass_{II,h}(\underline{\boldsymbol{u}}_h;\cdot):\underline{V}^{k,k}_h\to\Real$ are obtained by element-by-element assembly of the local residuals.
\begin{scheme}[\hhohp: HHO scheme with hybrid pressure]\label{scheme:hho.II}
  Find $(\underline{\boldsymbol{u}}_h,\underline{p}_h)\in\underline{\boldsymbol{V}}_h^{k+1,k}\times\underline{V}^{k,k}_h$ such that
  \begin{equation}\label{eq:hho.II}
    \begin{alignedat}{2}
      r_{II,h}^\momentum((\underline{\boldsymbol{u}}_h,\underline{p}_h);\underline{\boldsymbol{v}}_h) &= 0 &\qquad& \forall \underline{\boldsymbol{v}}_h\in\underline{\boldsymbol{V}}_h^{k+1,k},
      \\
      r_{II,h}^\mass(\underline{\boldsymbol{u}}_h;\underline{q}_h) &= 0 &\qquad& \forall \underline{q}_h\in\underline{V}_h^{k,k}.
    \end{alignedat}
  \end{equation}
\end{scheme}

The HHO method \eqref{eq:hho.II} yields a velocity approximation that is pointwise divergence free (as can be checked adapting the argument of \cite[Proposition 1]{Rhebergen.Wells:18}) and improves by one order the $h$-convergence rates of the HDG method proposed in \cite{Rhebergen.Wells:18}, since it relies on an HHO discretization of the viscous term (cf. the discussion in \cite{Cockburn.Di-Pietro.ea:16} and also \cite[Section 5.1.6]{Di-Pietro.Droniou:20}).
A key point consists in using element unknowns for the velocity one degree higher than face unknowns.
Notice that seeking the velocity in the space $\underline{\boldsymbol{V}}_T^{k+1,k}$ as opposed to $\underline{\boldsymbol{V}}_T^{k,k}$ does not alter the number of globally coupled unknowns, as all velocity degrees of freedom attached to the mesh elements can be removed from the global linear system by static condensation procedures similar to the ones discussed in Section \ref{sec:statCondStrat}.

\subsection{DG scheme}

The third approximation of the Stokes problem is based on discontinuous approximations of both the velocity and the pressure.
Specifically, we use the BR2 formulation for the vector Laplace operator (see \cite{Bassi.Rebay.ea:97} and also \cite[Section 5.3.2]{Di-Pietro.Ern:12}) together with a stabilized equal order pressure-velocity coupling.
Fix a polynomial degree $k\ge 1$ and let $T\in\Th$.
We define the \emph{local discrete gradient} $\boldsymbol{\mathfrak{G}}_T^k : H^1(\Th)^d \rightarrow \Poly{k}(T)^{d\times d}$ such that, for all $\vVec\in H^1(\Th)^d$,
\[
\int_T \boldsymbol{\mathfrak{G}}_T^k(\vVec) : \boldsymbol{\tau}
= \int_T \nabla \vVec_{|T} : \boldsymbol{\tau}
- \sum_{F \in \FT^{\internal,\dir}} \frac12 \int_F \left( \normal_{TF} \otimes \jump{\vVec}_{TF} \right) : \boldsymbol{\tau}
\qquad
\forall\boldsymbol{\tau} \in \Poly{k}(T)^{d \times d},
\]
where, for any $F\in\FT^{\internal,\dir}$, the jump of $\vVec$ across $F$ is defined as
\begin{equation*} 
  \jump{\vVec}_{TF} \coloneq\begin{cases}
  \vVec_{|T} - \vVec_{|T'} & \text{if $F\in\FT^\internal\cap\mathcal{F}_{T'}^\internal$ with $T,T'\in\Th$, $T\neq T'$},
  \\
  2(\vVec_T - \boldsymbol{g}_\dir) & \text{if $F\in\FT^\dir$}.
  \end{cases}
\end{equation*}
Introducing, for all $F\in\FT^{\internal,\dir}$, the jump lifting operator $\boldsymbol{\mathfrak{L}}_{FT}^k: L^2(F)^d \rightarrow \Poly{k}(T)^{d \times d}$ such that, for all $\boldsymbol{\varphi}\in L^2(F)^d$ and all $\boldsymbol{\tau}\in\Poly{k}(T)^{d\times d}$,
\[
  \int_T \boldsymbol{\mathfrak{L}}_{FT}^k (\boldsymbol{\varphi}) : \boldsymbol{\tau}
  =
  \frac{1}{2} \int_F \left(\normal_{TF} \otimes \boldsymbol{\varphi} \right) : \boldsymbol{\tau},
\]
it holds, for all $\vVec\in H^1(\Th)^d$,
\[
\boldsymbol{\mathfrak{G}}_T^k(\vVec)
= \nabla \vVec_{|T}
- \sum_{F\in\FT^{\internal,\dir}} \boldsymbol{\mathfrak{L}}_{FT}^k(\jump{\vVec}_{TF}).
\]
Given $(\uVec_h, p_h) \in \Poly{k}(\Th)^d\times \Poly{k}(\Th)$,
the local residual
$r^\momentum_{III,T}((\uVec_h,p_h);\cdot):\Poly{k}(T)^d\to\Real$ of the discrete momentum equation and
$r^\mass_{III,T}((\uVec_h,p_h);\cdot):\Poly{k}(T)\to\Real$ of the discrete mass equation are such that,
for all $\vVec_T\in \Poly{k}(T)^d$ and all $q_T\in\Poly{k}(T)$,
\begin{subequations}
\begin{align*} 
  r^\momentum_{III,T}((\uVec_h,p_h);\vVec_T)
  &\coloneq
  \begin{aligned}[t]
     &\int_T \boldsymbol{\mathfrak{G}}_T^k(\uVec_h) :  \nabla \vVec_T
     -  \sum_{F \in \FT^{\internal,\dir}} \int_F \left[
      \normal_{TF} \cdot \averg{ \nabla \uVec_T - \eta_F \boldsymbol{\mathfrak{L}}_{FT}^k(\jump{\uVec_h}_{TF})}_F
      \right] \cdot \vVec_T
    \\
    &{\;} - \int_T p_T \, (\nabla \cdot \vVec_T)
    + \sum_{F \in \FT^{\internal,\dir}} \int_F \averg{p_h}_F \, (\vVec_T \cdot \normal_{TF})
    \\
    &{\;} 
    - \int_{\aT} \boldsymbol{f} \cdot \vVec_T
    - \sum_{F \in \FT^\neu} \int_F \boldsymbol{g}_\neu \cdot \vVec_T,
  \end{aligned}
  \\ 
  r^\mass_{III,T}((\uVec_h,p_h);q_T)
  &\coloneq
  \begin{aligned}[t]
    & \int_T \uVec_T \cdot \nabla q_T
    - \sum_{F \in \FT^{\internal,\dir}} \int_F \averg{\uVec_h}_F \cdot \normal_{TF} \, q_T
    + \sum_{F \in \FT^\internal} h_F \int_F \jump{p_h}_{TF}  \, q_T, 
  \end{aligned}
\end{align*}
\end{subequations}
where, for all $\varphi\in H^1(\Th)$ and all $F\in\Fh$,
\[
\averg{\varphi}_F \coloneq \begin{cases}
  \frac12\left(\varphi_{|T} + \varphi_{|T'}\right) & \text{if $F\in\FT^\internal\cap\mathcal{F}_{T'}^\internal$ with $T,T'\in\Th$, $T\neq T'$},
  \\
  \varphi_{|F} & \text{otherwise},
\end{cases}
\]
with the understanding that the average operator acts componentwise when applied to vector and tensor functions,
and  
\begin{equation*} 
\eta_F >
\begin{cases}
  \max \big(\card{\FT}, \card{\FTp}\big) & \text{if $F\in\FT^\internal\cap\mathcal{F}_{T'}^\internal$ with $T,T'\in\Th$, $T\neq T'$},
    \\  
    \card{\FT} & \text{if $F \in \FT^\dir$}.
\end{cases}
\end{equation*}
The global residuals
$r^\momentum_{III,h}((\uVec_h,p_h);\cdot):\Poly{k}(\Th)^d\to\Real$ and
$r^\mass_{III,h}((\uVec_h,p_h);\cdot):\Poly{k}(\Th)\to\Real$
are obtained by element-by-element assembly of local residuals.
\begin{scheme}[\dg: DG scheme]\label{scheme:dg}
  Find $(\boldsymbol{u}_h,p_h)\in\Poly{k}(\Th)^d\times\Poly{k}(\Th)$ such that
  \begin{equation*}
    \begin{alignedat}{2}
      r_{III,h}^\momentum\big( (\boldsymbol{u}_h,{p}_h);\boldsymbol{v}_h\big)
      &= \sum_{T \in \Th} r_{III,T}^\momentum\big( (\boldsymbol{u}_h,{p}_h);\boldsymbol{v}_h{}_{|T}\big) = 0 &\qquad& \forall \boldsymbol{v}_h\in\Poly{k}(\Th)^d,
      \\
      r_{III,h}^\mass\big((\boldsymbol{u}_h,p_h);q_h)
      &= \sum_{T \in \Th} r_{III,T}^\mass\big((\boldsymbol{u}_h,p_h);q_h{}_{|T}) = 0 &\qquad& \forall q_h\in\Poly{k}(\Th).
    \end{alignedat}
  \end{equation*}
\end{scheme}


\section{$p$-Multilevel solution strategy}\label{sec:p-multilevel.strategies}

We consider $L$ coarse problems, indexed as $\ell =1,...,L$.
Given a polynomial degree $k\ge 0$ (for Schemes \ref{scheme:hho} and \ref{scheme:hho.II}) or $k\ge 1$ (for Scheme \ref{scheme:dg}), we set
\begin{equation*}
  k_0 \coloneq k,
\end{equation*}
the reference polynomial degree on the fine level, and denote by $\kl$ the polynomial degree at level $\ell$.
Coarsening is achieved taking $\klpo < \kl$.
The notation for the three schemes discussed in Section \ref{sec:threeNCD} is summarized in Table \ref{tab:notation}.
Notice that, for the sake of simplicity, we only consider the equal-order version of Scheme \ref{scheme:hho}, where both element and face velocity unknowns have the same polynomial degree.

\begin{table}\centering
  \renewcommand{\arraystretch}{1.2}
  \begin{tabular}{ccccc}
    \toprule
    Scheme index & Scheme label & Fine discrete space & Coarse discrete spaces & Coarsest level \\
    \midrule
    \ref{scheme:hho}   &\hhodp     & $\underline{\boldsymbol{V}}_h^{k_0,k_0}$,    $\Poly{k_0}(\Th)$  & $\underline{\boldsymbol{V}}_h^{\kl,\kl}$, $\Poly{\kl}(\Th)$    & $k_L \geq 0$    \\
    \ref{scheme:hho.II}& \hhohp   & $\underline{\boldsymbol{V}}_h^{k_0+1,k_0}$,  $\Poly{k_0}(\Th)$  & $\underline{\boldsymbol{V}}_h^{\kl+1,\kl}$, $\Poly{\kl}(\Th)$  & $k_L \geq 0$    \\
    \ref{scheme:dg}    & \dg   & $\Poly{k_0}(\Th)^d$, $\Poly{k_0}(\Th)$                            & $\Poly{\kl}(\Th)^d$, $\Poly{\kl}(\Th)$                       & $k_L \geq 1$  \\
    \bottomrule
  \end{tabular}
  \caption{Notation for the $p$-multilevel solver\label{tab:notation}. We only consider the equal-order version of Scheme \ref{scheme:hho}, where both element and face velocity unknowns have the same polynomial degree.}  
\end{table}

\subsection{Intergrid transfer operators}

Denoting by $X\in\Th\cup\Fh$ a mesh element or face, the \textit{prolongation} operator $\IntOp_{\ell+1}^{\ell,X}: \Poly{\klpo}(X) \rightarrow \Poly{\kl}(X)$ from level $\ell+1$ to level $\ell$ is the injection $\Poly{\klpo}(X)\hookrightarrow\Poly{k}(X)$.
The prolongation operator $\IntOp^0_\ell$ from level $\ell$ to level $0$ can be recursively defined by the composition of one level prolongation operators:
\begin{equation*}
  \IntOp_\ell^0 = \IntOp^0_1 \, \IntOp^1_2 \,... \, \IntOp_\ell^{\ell-1}.
\end{equation*}
The \textit{restriction} operator $\IntOp_{\ell,X}^{\ell+1}: \Poly{\kl}(X) {\rightarrow} \Poly{\klpo}(X)$ from level $\ell$ to level $\ell+1$ is simply taken equal to the $L^2$-orthogonal projector on $\Poly{\klpo}(X)$, that is, for all $w_{X,\ell} \in \Poly{\kl}(X)$, we set
\begin{equation}\label{def:restOp}
  \IntOp_{\ell,X}^{\ell+1} w_{X,\ell} \coloneq \pi_{X}^{\ell+1} w_{X,\ell}.
\end{equation}
The restriction operator $\IntOp^\ell_0$ from level $0$ to level $\ell$ is again obtained by composition:
\begin{equation*}
  \IntOp^\ell_0 = \IntOp^\ell_{\ell-1} \,... \,\IntOp^2_1 \,\IntOp^1_0.
\end{equation*}
It can be checked that $\IntOp_{\ell,X}^{\ell+1}$ is the transpose of $\IntOp_{\ell+1}^{\ell,X}$ with respect to the $L^2(X)$-inner product.
When applied to vector-valued functions, intergrid transfer operators act component-wise and are denoted using boldface font by $\IntOpB_{\ell+1}^{\ell,X}$, $\IntOpB_{\ell,X}^{\ell+1}$.
The global restriction operator $\underline{\IntOpB}^{\ell+1}_{\ell} : \underline{\boldsymbol{V}}_h^{\kl',\kl}\rightarrow \underline{\boldsymbol{V}}_h^{\klpo',\klpo}$ for HHO spaces is defined setting: For all $\underline{\vVec}_{h,\ell}\in\underline{\boldsymbol{V}}_h^{\kl',\kl}$,
\[
  \underline{\IntOpB}^{\ell+1}_{\ell} \underline{\vVec}_{h,\ell} \coloneq \big(
  (\IntOpB_{\ell,T}^{\ell+1} \vVec_{T,\ell})_{T\in\Th},
  (\IntOpB_{\ell,F}^{\ell+1} \vVec_{F,\ell})_{F\in\Fh}
  \big),
\]
while the global restriction operator for DG spaces ${\IntOpB}^{\ell+1}_{\ell} : \Poly{\kl}(\Th)^d\rightarrow \Poly{\klpo}(\Th)^d$ is obtained patching the element restriction operators: For all ${\vVec}_{h,\ell}\in\Poly{\kl}(\Th)^d$,
\[
\big({\IntOpB}^{\ell+1}_{\ell} {\vVec}_{h,\ell}\big)_{|T} \coloneq \IntOpB_{\ell,T}^{\ell+1} (\vVec_{h,\ell})_{|T} \qquad\forall T\in\Th.
\]

\subsection{Inherited multilevel operators}

For any $\ell=1,\ldots,L$ set, for the sake of brevity,
\[
\boldsymbol{W}_{I,h}^\ell\coloneq\underline{\boldsymbol{V}}^{\kl,\kl}_h \times \Poly{\kl}(\Th),\qquad
\boldsymbol{W}_{II,h}^\ell\coloneq\underline{\boldsymbol{V}}^{\kl+1,\kl}_h \times \underline{V}^{\kl,\kl}_h,\qquad
\boldsymbol{W}_{III,h}^\ell\coloneq\Poly{\kl}(\Th)^d \times \Poly{\kl}(\Th).
\]
The coarse residuals for the momentum and mass continuity equations for the schemes of Section \ref{sec:threeNCD} corresponding to a velocity-pressure couple at level $\ell$ are obtained evaluating the corresponding fine residuals defined in Section \ref{sec:HHO.schemes} at the prolongation of the given function, i.e.:
For $\ell=1,\ldots,L$,
\begin{itemize}
\item {\bf Scheme \ref{scheme:hho} (\hhodp).} \newline
  Given $(\underline{\boldsymbol{u}}_{h,\ell}, p_{h,\ell} ) \in \boldsymbol{W}_{I,h}^\ell$, $r_{I,\ell}\big((\underline{\boldsymbol{u}}_{h,\ell}, p_{h,\ell} ); \cdot\big):\boldsymbol{W}_{I,h}^\ell\to\Real$ is such that, for all $(\underline{\boldsymbol{v}}_{h,\ell}, q_{h,\ell} ) \in \boldsymbol{W}_{I,h}^\ell$,
  \begin{subequations}\label{def:coarseResHHO}
    \begin{alignat}{2}
      & \; \; r_{I,\ell}&\big((\underline{\boldsymbol{u}}_{h,\ell}, p_{h,\ell} ); (\underline{\boldsymbol{v}}_{h,\ell}, q_{h,\ell}) \big)
      &\coloneq
        {r}^\momentum_{I,\ell}\big((\underline{\boldsymbol{u}}_{h,\ell},p_{h,\ell});\underline{\boldsymbol{v}}_{h,\ell}\big)
        + r^\mass_{I,\ell} (\underline{\boldsymbol{u}}_{h,\ell};q_{h,\ell}) 
      \nonumber \\
      &\text{with}&
        {r}^\momentum_{I,\ell}\big((\underline{\boldsymbol{u}}_{h,\ell},p_{h,\ell});\underline{\boldsymbol{v}}_{h,\ell}\big)
        &\coloneq {r}^\momentum_{I,h}\big( (\underline{\IntOpB}_\ell^0\underline{\boldsymbol{u}}_{h,\ell}, \IntOp_\ell^0 p_{h,\ell});\underline{\IntOpB}_\ell^0\underline{\boldsymbol{v}}_{h,\ell} \big) \;
      \nonumber \\
        &\text{and}&
      r^\mass_{I,\ell} (\underline{\boldsymbol{u}}_{h,\ell};q_{h,\ell}) &\coloneq r^\mass_{I,h} (\underline{\IntOpB}_\ell^0\underline{\boldsymbol{u}}_{h,\ell};\IntOp_\ell^0 q_{h,\ell}); \nonumber
    \end{alignat}
  \end{subequations}

\item {\bf Scheme \ref{scheme:hho.II} (\hhohp).} \newline 
  Given $(\underline{\boldsymbol{u}}_{h,\ell},\underline{p}_{h,\ell}) \in \boldsymbol{W}_{II,h}^\ell$,
  $r_{II,\ell}\big((\underline{\boldsymbol{u}}_{h,\ell},\underline{p}_{h,\ell});\cdot\big):\boldsymbol{W}_{II,h}^\ell\to\Real$ is such that,
  for all $(\underline{\boldsymbol{v}}_{h,\ell},\underline{q}_{h,\ell}) \in \boldsymbol{W}_{II,h}^\ell$,
  \begin{subequations}\label{def:coarseResHHO.II}
    \begin{alignat}{2}
        & \;\;  r_{II,\ell} &\big((\underline{\boldsymbol{u}}_{h,\ell},\underline{p}_{h,\ell}); (\underline{\boldsymbol{v}}_{h,\ell},\underline{q}_{h,\ell})\big)
        &\coloneq {r}^\momentum_{II,\ell}\big((\underline{\boldsymbol{u}}_{h,\ell},\underline{p}_{h,\ell});\underline{\boldsymbol{v}}_{h,\ell}\big)
        +  r^\mass_{II,\ell} (\underline{\boldsymbol{u}}_{h,\ell};\underline{q}_{h,\ell}) \; 
      \nonumber \\
      &\text{
        with
      } &
        {r}^\momentum_{II,\ell}\big((\underline{\boldsymbol{u}}_{h,\ell},\underline{p}_{h,\ell});\underline{\boldsymbol{v}}_{h,\ell}\big)
        &\coloneq {r}^\momentum_{II,h}\big((\underline{\IntOpB}_\ell^0\underline{\boldsymbol{u}}_{h,\ell}, \underline{\IntOp}_\ell^0 \underline{p}_{h,\ell});\underline{\IntOpB}_\ell^0\underline{\boldsymbol{v}}_{h,\ell}\big) \; 
      \nonumber \\
      &\text{
        and
      } &
      r^\mass_{II,\ell} (\underline{\boldsymbol{u}}_{h,\ell};\underline{q}_{h,\ell})
        &\coloneq r^\mass_{II,h} (\underline{\IntOpB}_\ell^0\underline{\boldsymbol{u}}_{h,\ell};\underline{\IntOp}_\ell^0 \underline{q}_{h,\ell}); \nonumber
    \end{alignat}
  \end{subequations}
\item {\bf Scheme \ref{scheme:dg} (\dg).} \newline Given $({\boldsymbol{u}}_{h,\ell} , p_{h,\ell}) \in \boldsymbol{W}_{III,h}^\ell$,
  $r_{III,\ell}\big(({\boldsymbol{u}}_{h,\ell} , p_{h,\ell});\cdot\big):\boldsymbol{W}_{III,h}^\ell\to\Real$ is such that,
  for all $({\boldsymbol{v}}_{h,\ell} , q_{h,\ell}) \in \boldsymbol{W}_{III,h}^\ell$,
  \begin{subequations}\label{def:coarseResDG}
    \begin{alignat}{2}
        & \; r_{III,\ell}&\big(({\boldsymbol{u}}_{h,\ell} , p_{h,\ell});({\boldsymbol{v}}_{h,\ell} , q_{h,\ell})\big)
        &\coloneq {r}^\momentum_{III,\ell}\big(({\boldsymbol{u}}_{h,\ell},p_{h,\ell});{\boldsymbol{v}}_{h,\ell} \big)
        + r^\mass_{III,\ell} \big(({\boldsymbol{u}}_{h,\ell}, p_{h,\ell});q_{h,\ell} \big) 
      \nonumber \\
      & \text{
        with
      } &
        {r}^\momentum_{III,\ell}\big(({\boldsymbol{u}}_{h,\ell},p_{h,\ell});{\boldsymbol{v}}_{h,\ell} \big)
        &\coloneq {r}^\momentum_{III,h}\big(({\IntOpB}_\ell^0{\boldsymbol{u}}_{h,\ell}, \IntOp_\ell^0 p_{h,\ell});{\IntOpB}_\ell^0{\boldsymbol{v}}_{h,\ell}\big) \;
      \nonumber \\
      & \text{
        and
      } &
      r^\mass_{III,\ell} \big(({\boldsymbol{u}}_{h,\ell}, p_{h,\ell});q_{h,\ell} \big)
      & \coloneq r^\mass_{III,h} \big(({\IntOpB}_\ell^0{\boldsymbol{u}}_{h,\ell}, \IntOp_\ell^0 p_{h,\ell});\IntOp_\ell^0q_{h,\ell} \big).  \nonumber
    \end{alignat}
  \end{subequations}
\end{itemize}

  Fix $\bullet\in\{{\rm I},{\rm II},{\rm III}\}$, $\ell=0,\ldots,L$, and denote by $(\cdot,\cdot)$ an inner product in $\boldsymbol{W}_{\bullet,h}^\ell$.
  Let $\A_{h,\ell}:\boldsymbol{W}_{\bullet,h}^\ell\to\boldsymbol{W}_{\bullet,h}^\ell$ be the operator corresponding to the linear part of the residual $r_{\bullet,\ell}$, i.e., for all $\wVec_{h,\ell}\in\boldsymbol{W}_{\bullet,h}^\ell$, $(\A_{h,\ell}\wVec_{h,\ell},\zVec_{h,\ell}) = r_{\bullet,\ell}(\wVec_{h,\ell};\zVec_{h,\ell}) - r_{\bullet,\ell}(\boldsymbol{0};\zVec_{h,\ell})$ for all $\zVec_{h,\ell}\in\boldsymbol{W}_{\bullet,h}^\ell$.
  Letting $\bVec_{h,\ell}\in\boldsymbol{W}_{\bullet,h}^\ell$ denote the Riesz representation of the affine part of the residual such that $(\bVec_{h,\ell},\zVec_{h,\ell}) = r_{\bullet,\ell}(\boldsymbol{0};\zVec_{h,\ell})$ for all $\zVec_{h,\ell}\in\boldsymbol{W}_{\bullet,h}^\ell$,
 the global problem at level $\ell$ reads:
Find $\wVec_{h,\ell}\in\boldsymbol{W}_{\bullet,h}^\ell$ such that
\begin{equation*}
  \A_{h,\ell} \wVec_{h,\ell} = \bVec_{h,\ell}.
\end{equation*}
Besides the formal definition given above, coarse level operators can be efficiently inherited from the fine operators relying on the restriction and prolongation operators.
This computationally efficient strategy, also known as Galerkin projection, is detailed in Section \ref{sec:inhOpA} focusing on Scheme \ref{scheme:hho}.

\subsection{Multilevel $V$-cycle iteration}
\label{sec:MGViter}

The approximate solution $\overline{\wVec}_{h,\ell}$ to the global problem at level $\ell < L$ can be improved by means of one $V$-cycle iteration, as described in the following algorithm:

\begin{algorithm}[H]
$\mathrm{\textbf{Multilevel $V$-cycle: MG}}_{V}(\ell,\boldsymbol{b}_{h,\ell},\overline{\wVec}_{h,\ell})$
 \begin{algorithmic}[0]
   \STATE{\underline{\textbf{Pre-smoothing:}}}
   \STATE{$\overline{\wVec}_{h,\ell} = \mathrm{GMRES}(\A_{h,\ell}, \overline{\wVec}_{h,\ell}, \boldsymbol{b}_{h,\ell})$} 
   \STATE{\underline{\textbf{Compute the coarse grid correction} (recursion up to level $L$):}}
   \STATE $\boldsymbol{d}_{h,\ell+1} = \boldsymbol{\mathcal{I}}_{\ell}^{\ell+1} (\boldsymbol{b}_{h,\ell} - \A_{h,\ell} \overline{\wVec}_{h,\ell})$
    \IF {$(\ell+1 = L)$}
    \STATE $\boldsymbol{c}_{h,\ell+1} = \A_{h,\ell+1}^{-1} \boldsymbol{d}_{h,\ell+1}$
    \ENDIF
    \IF {$(\ell+1 < L)$}
   \STATE ${\boldsymbol{c}}_{h,\ell+1} = \textbf{MG}_{V}(\ell+1,\boldsymbol{d}_{h,\ell+1},\boldsymbol{0})$ 
    \ENDIF
   \STATE{\underline{\textbf{Apply the coarse grid correction:}}}
   \STATE $\overline{\wVec}_{h,\ell} = \overline{\wVec}_{h,\ell} + \boldsymbol{\mathcal{I}}_{\ell+1}^{\ell} \boldsymbol{c}_{h,\ell+1}$
   \STATE{\underline{\textbf{Post-smoothing:}}}
   \STATE{$\overline{\wVec}_{h,\ell} = \mathrm{GMRES}(\A_{h,\ell}, \overline{\wVec}_{h,\ell}, \boldsymbol{b}_{h,\ell})$} 
\end{algorithmic}
\end{algorithm}
\noindent where $\boldsymbol{d}_{h,\ell+1}$ is the restriction of the defect and $\boldsymbol{c}_{h,\ell+1}$ is the coarse grid correction.
All applications of prolongation and restriction operators involved in the multilevel $V$-cycle iteration are 
performed matrix-free, that is, without assembling the global sparse matrices associated to 
the operators $\boldsymbol{\mathcal{I}}_{\ell}^{\ell+1},\boldsymbol{\mathcal{I}}^{\ell}_{\ell+1}$.

In the pre- and post- smoothing steps, a few iterations of the Generalised Minimal Residual (GMRES) method preconditioned with 
an Incomplete Lower-Upper (ILU) factorization are performed in order to reduce the error 
$\boldsymbol{e}_{h,\ell} = \wVec_{h,\ell} - \overline{\wVec}_{h,\ell}$. 
Indeed, the components of the error associated to the highest-order basis functions at level $\ell$ 
are expected to be damped very fast, while the components of the error associated to lower-order basis functions 
are smoothed at a later stage when the recursion reaches coarser levels.

In the numerical tests of Section \ref{sec:numResults} we consider one $V$-cycle iteration as a preconditioner for the 
FGMRES (Flexible GMRES) iteration applied to solve the global problem  $\A_{h,0} \wVec_{h,0} = \bVec_{h,0}$.
We employ the solver and preconditioner framework provided by the PETSc library \cite{Balay.Abhyankar.ea:19}.


\section{Computational aspects}\label{sec:computational.aspects}

In what follows, we discuss some computational aspects for the Scheme \ref{scheme:hho} (HHO with discontinuous pressure).
Algebraic objects are denoted using sans serif font, with boldface distinguishing matrices from vectors.

\subsection{Static condensation}
\label{sec:staticCondHHO}

\subsubsection{Algebraic expression for the local residuals}

We assume that local bases for each polynomial space attached to mesh elements and faces have been fixed, so that bases for the global approximations spaces for the velocity and the pressure can be obtained by taking the Cartesian product of the latter.
Possible choices of local bases are discussed in \cite[Appendix B.1]{Di-Pietro.Droniou:20}.

The unknowns for a mesh element $T\in\Th$ correspond to the coefficients of the expansions of the velocity and pressure in the selected local bases.
Assuming that the velocity unknowns are ordered so that element velocities come first and boundary velocities next, these coefficients are collected in the vectors
\begin{equation*}
  \text{
    $\UVEC{U}_T=\begin{bmatrix} \VEC{U}_T \\ \VEC{U}_{\partial T} \end{bmatrix}$
    and $\VEC{P}_T$,
  }
\end{equation*}
where the block partition of the vector $\UVEC{U}_T$ is the one naturally induced by the selected ordering of velocity unknowns.

The local matrices corresponding to the HHO discretization of the viscous term (first two lines of the right-hand side of \eqref{HHOdiscrMNT}) and of the pressure-velocity coupling (first line of the right-hand side of \eqref{HHOdiscrCNT}) are
\begin{equation*}
  \MAT{A}_T = 
  \begin{bmatrix}
    \MAT{A}_{T  T} & \MAT{A}_{T  \partial T}  \\
    \MAT{A}_{\partial T T} & \MAT{A}_{\partial T  \partial T}
  \end{bmatrix},\qquad
  \MAT{B}_T = \begin{bmatrix}
    \MAT{B}_{T  T} & \MAT{B}_{T  \partial T}
  \end{bmatrix},
\end{equation*}
where again the block partition is the one induced by the ordering of velocity unknowns.
Details on the construction of the matrix $\MAT{A}_T$ can be found in \cite[Appendix B.2]{Di-Pietro.Droniou:20}.

\begin{remark}[Block structure]\label{rem:structure.ApTpT}
  Denoting by $N$ the number of faces of $T$, the block structure of the matrix $\MAT{A}_T$ can be further detailed as follows:
  \begin{equation}\label{eq:AT.block.refined}
      \MAT{A}_T = 
      \left[\begin{array}{c|ccc}
          \MAT{A}_{T  T} & \MAT{A}_{T F_1} & \cdots & \MAT{A}_{T F_N} \\
          \hline
          \MAT{A}_{F_1 T} & \MAT{A}_{F_1 F_1} & \cdots & \MAT{A}_{F_1 F_N} \\
          \vdots & \vdots & \ddots & \vdots \\
          \MAT{A}_{F_N T} & \MAT{A}_{F_N F_1} & \cdots & \MAT{A}_{F_N F_N}
      \end{array}\right].
  \end{equation}
  Assume that the velocity unknowns attached to $T$ and its faces are ordered by component.
  Since the viscous term is modelled in \eqref{stokesProb:momentum} applying the Laplace operator to each velocity component, each block in the decomposition \eqref{eq:AT.block.refined} is itself block-diagonal, and be efficiently constructed starting from the corresponding matrix for the scalar Laplace operator.
\end{remark}

Introducing the vector representations
$\UVEC{R}_{I,T}^\momentum=\begin{bmatrix}\VEC{R}_{I,T}^\momentum\\ \VEC{R}_{I,\partial T}^\momentum\end{bmatrix}$ and $\VEC{R}_{I,T}^\mass$ of the residual linear forms defined by \eqref{HHOdiscr:residuals},
$\VEC{G}_{\partial T}$ of the terms involving the boundary data in the last line of \eqref{HHOdiscrMNT},
$\VEC{F}_T$ of the term involving the volumetric body force in the last line of \eqref{HHOdiscrMNT},
and $\widehat{\VEC{G}}_{\partial T}$ of the last term in the right-hand side of \eqref{HHOdiscrCNT}, it holds
\begin{equation}
  \label{def:HHOAB}
  \begin{bmatrix}
    \VEC{R}^\momentum_{I,T} \\
    \VEC{R}^\momentum_{I, \partial T} \\
    \VEC{R}^{\mass}_{I,T}
  \end{bmatrix}
  =
  \begin{bmatrix}
    \MAT{A}_{T T} & \MAT{A}_{T\partial T} & \MAT{B}_{T  T}\trans    \\
    \MAT{A}_{\partial T  T} &  \MAT{A}_{\partial T  \partial T} & \MAT{B}_{T\partial T}\trans  \\
    \MAT{B}_{T  T} & \MAT{B}_{T  \partial T} & \MAT{0}
  \end{bmatrix}
  \begin{bmatrix}
    \VEC{U}_T \\
    \VEC{U}_{\partial T} \\
    \VEC{P}_T
  \end{bmatrix}
  -
  \begin{bmatrix}
    \VEC{F}_T \\
    \VEC{G}_{\partial T} \\
    \widehat{\VEC{G}}_{\partial T}
  \end{bmatrix}.
\end{equation}

\subsubsection{Static condensation strategies}
\label{sec:statCondStrat}

The discrete problem \eqref{eq:hho} is obtained enforcing that the global residuals be zero, which requires the solution of a global linear system.
The size of this linear system can be reduced by statically condensing the element velocity unknowns and, possibly, the pressure unknowns corresponding to high-order modes inside each element.
In what follows, we discuss two possible static condensations procedures leading to global systems with different features.

\paragraph{\hhodp{} \vcond{}: Static condensation of velocity element unknowns.} 

The first static condensation procedure hinges on the observation that, given a mesh element $T\in\Th$, the velocity unknowns collected in $\VEC{U}_T$ are not directly coupled with unknowns attached to mesh elements other than $T$.
As a result, enforcing that the residuals in the left-hand side of \eqref{def:HHOAB} be zero, $\VEC{U}_T$ can be locally eliminated by expressing it in terms of $\VEC{U}_{\partial T}$ and $\VEC{P}_T$ by computing the Schur complement
\begin{equation}\label{eq:ST.1}
  \renewcommand{\arraystretch}{1.3}
  \MAT{S}_T^{\vc}
  = \begin{bmatrix}
    \MAT{S}_{\partial T\partial T}^{\vc} & \MAT{S}_{\partial T T}^{\vc} \\
    \MAT{S}_{T\partial T}^{\vc} & \MAT{S}_{T T}^{\vc} \\
  \end{bmatrix}
  \coloneq
  \begin{bmatrix}
    \MAT{A}_{\partial T  \partial T} & \MAT{B}_{T\partial T}\trans  \\
    \MAT{B}_{T  \partial T} & \MAT{0}
  \end{bmatrix}
  - \begin{bmatrix}
    \MAT{A}_{\partial T  T} \\
    \MAT{B}_{T  T}
  \end{bmatrix}\MAT{A}_{TT}^{-1}\begin{bmatrix}
    \MAT{A}_{T\partial T} & \MAT{B}_{T  T}\trans
  \end{bmatrix}
\end{equation}
of the block $\MAT{A}_{TT}$ in the matrix in the right-hand side of \eqref{def:HHOAB}.
With this static condensation strategy, the zero residual condition translates into
\begin{equation*}
  \MAT{S}_T^{\vc}\begin{bmatrix}
    \VEC{U}_{\partial T}
    \\
    \VEC{P}_T
  \end{bmatrix}  
  = \begin{bmatrix}
    \VEC{G}_{\partial T}
    \\
    \widehat{\VEC{G}}_{\partial T}
  \end{bmatrix}
  - \begin{bmatrix}
    \MAT{A}_{\partial T T} \\
    \MAT{B}_{TT}
  \end{bmatrix}
  \MAT{A}_{TT}^{-1}\VEC{F}_T.
\end{equation*}

\paragraph{\hhodp{} \vpcond{}: Static condensation of velocity element unknowns and pressure modes.} 

The second static condensation strategy was originally suggested in \cite{Aghili.Boyaval.ea:15} in the framework of HHO methods and later detailed in  \cite[Section 6]{Di-Pietro.Ern.ea:16}.
Assume that the basis for the pressure inside each mesh element $T\in\Th$ is selected so that the first degree of freedom corresponds to the mean value of the pressure inside $T$ and the remaining basis functions are $L^2$-orthogonal to the first (this condition typically requires the use of modal bases).
Let now a mesh element $T\in\Th$ be fixed.
The above choice for the pressure basis induces the following partitions of the pressure unknowns and of the pressure-velocity coupling matrix:
\begin{equation*}
  \VEC{P}_T = \begin{bmatrix}
    \overline{\VEC{P}}_T \\ \widetilde{\VEC{P}}_T
  \end{bmatrix},\qquad
  \MAT{B}_T = \begin{bmatrix}
    \overline{\MAT{B}}_{T  T} & \overline{\MAT{B}}_{T  \partial T} \\
    \widetilde{\MAT{B}}_{T  T} & \widetilde{\MAT{B}}_{T  \partial T} \\
  \end{bmatrix},
\end{equation*}
where $\overline{\VEC{P}}_T\in\Real$ is the mean value of the pressure inside $T$,
$\widetilde{\VEC{P}}_T$ is the vector corresponding to high-order pressure modes,
and the matrix $\MAT{B}_T$ has been partitioned row-wise according to this decomposition.
Enforcing that the residuals be zero in \eqref{def:HHOAB} and rearranging the unknowns and equations, we infer that the discrete solution satisfies
\begin{equation}\label{eq:zero-residual.T}\renewcommand{\arraystretch}{1.3}
  \left[\begin{array}{cc|cc}
      \MAT{A}_{T T} & \widetilde{\MAT{B}}_{TT}\trans & \MAT{A}_{T\partial T} & \overline{\MAT{B}}_{T  T}\trans
      \\     
      \widetilde{\MAT{B}}_{T  T} & \MAT{0} & \widetilde{\MAT{B}}_{T  \partial T} & \MAT{0}
      \\ \hline
      \MAT{A}_{\partial T  T} & \widetilde{\MAT{B}}_{T\partial T}\trans & \MAT{A}_{\partial T  \partial T} & \overline{\MAT{B}}_{T\partial T}\trans
      \\
      \overline{\MAT{B}}_{TT} & \MAT{0} & \overline{\MAT{B}}_{T\partial T} & \MAT{0}
    \end{array}
  \right]
  \left[\begin{array}{c}
    \VEC{U}_T \\
    \widetilde{\VEC{P}}_T \\ \hline
    \VEC{U}_{\partial T} \\
    \overline{\VEC{P}}_T
  \end{array}\right]
  =
  \left[\begin{array}{c}
    \VEC{F}_T \\
    \VEC{0} \\ \hline
    \VEC{G}_{\partial T} \\
    \widehat{\VEC{G}}_{\partial T}
  \end{array}\right].
\end{equation}
The only unknowns that are globally coupled are those collected in the subvector
$\begin{bmatrix}
  \VEC{U}_{\partial T} \\
  \overline{\VEC{P}}_T
\end{bmatrix}$,
while the remaining unknowns collected in
$\begin{bmatrix}
  \VEC{U}_T \\
  \widetilde{\VEC{P}}_T
\end{bmatrix}$
can be eliminated by expressing them in terms of the former.
After performing this local elimination, the condition \eqref{eq:zero-residual.T} that the residuals associated with $T$ be zero becomes:
\begin{equation}\label{eq:ST.2}
  \MAT{S}_T^{\vpc}\begin{bmatrix}
  \VEC{U}_{\partial T} \\
  \overline{\VEC{P}}_T
  \end{bmatrix}
  = \begin{bmatrix}
    \VEC{G}_{\partial T} \\
    \widehat{\VEC{G}}_{\partial T}
  \end{bmatrix}
  - \begin{bmatrix}
    \MAT{A}_{\partial T T} & \widetilde{\MAT{B}}_{\partial TT}\trans \\
    \overline{\MAT{B}}_{TT} & \MAT{0}
  \end{bmatrix}
  \begin{bmatrix}
    \MAT{A}_{T T} & \widetilde{\MAT{B}}_{TT}\trans  \\
    \widetilde{\MAT{B}}_{T  T} & \MAT{0}
  \end{bmatrix}^{-1}
  \begin{bmatrix}
    \VEC{F}_T \\
    \VEC{0}
  \end{bmatrix},
\end{equation}
where $\MAT{S}_T^{\vpc}$ denotes the Schur complement of the top left block of the matrix in \eqref{eq:zero-residual.T}, that is,
\begin{equation*}
  \renewcommand{\arraystretch}{1.3}
  \MAT{S}_T^{\vpc}
  = \begin{bmatrix}
    \MAT{S}_{\partial T\partial T}^{\vpc} & \MAT{S}_{\partial T T}^{\vpc} \\
    \MAT{S}_{T\partial T}^{\vpc} & \MAT{S}_{T T}^{\vpc} \\
  \end{bmatrix}
  \coloneq \begin{bmatrix}
    \MAT{A}_{\partial T  \partial T} & \overline{\MAT{B}}_{T\partial T}\trans
    \\
    \overline{\MAT{B}}_{T\partial T} & \MAT{0}
  \end{bmatrix} - \begin{bmatrix}
    \MAT{A}_{\partial T  T} & \widetilde{\MAT{B}}_{T\partial T}\trans
    \\
    \overline{\MAT{B}}_{TT} & \MAT{0}
  \end{bmatrix}\begin{bmatrix}
    \MAT{A}_{T T} & \widetilde{\MAT{B}}_{TT}\trans
      \\     
      \widetilde{\MAT{B}}_{T  T} & \MAT{0}
  \end{bmatrix}^{-1}\begin{bmatrix}
    \MAT{A}_{T\partial T} & \overline{\MAT{B}}_{T  T}\trans
    \\     
    \widetilde{\MAT{B}}_{T  \partial T} & \MAT{0}
  \end{bmatrix}.
\end{equation*}

\begin{remark}[Differences between the static condensation strategies]
\label{rem:StatCond}
  The two static condensation strategies outlined above coincide for $k{=}0$.
  For $k\ge 1$, the first, obvious difference is that the second results in a smaller global system, since high-order pressure unknowns are eliminated in addition to element-based velocity unknowns.
  There is, however, a second, more subtle difference.
  As a matter of fact, while the block $\MAT{S}_{\partial T \partial T}^{\vpc}$ in \eqref{eq:ST.2} is full, the block $\MAT{S}_{\partial T \partial T}^{\vc}$ in \eqref{eq:ST.1} preserves the pattern of $\A_{\partial T \partial T}$ (which is composed of block-diagonal blocks, see Remark \ref{rem:structure.ApTpT}).
  As a result, the first static condensation strategy results in a sparser, albeit larger, matrix.
  The numerical tests in the next section show that sparsity prevails over size, so that the first static condensation strategy is in fact the more efficient.
  
  Notice that this difference would disappear if we replaced the Laplace operator in the momentum equation \eqref{stokesProb:momentum} by $\DIV(\nu \boldsymbol{\nabla}_{\rm s}\cdot)$, with $\boldsymbol{\nabla}_{\rm s}$ denoting the symmetric part of the gradient operator applied to vector-valued fields, as would be required for a viscosity coefficient $\nu:\Omega\to\Real^+$ variable in space.
\end{remark}

\subsection{Inheritance by means of Galerkin projections}
\label{sec:inhOpA}

We show in this section how the operators can be inherited from level $\ell$ to $\ell+1$.
For $X$ mesh element or face, we let $\{\psi_1^{X,\ell},\psi_2^{X,\ell},...,\psi_P^{X,\ell}\}$ be a basis of $\Poly{\kl}(X)$ (with $P$ denoting the dimension of this vector space) and $\{\psi_1^{X,\ell+1},\psi_2^{X,\ell+1},...,\psi_Q^{X,\ell+1}\}$ a basis of $\Poly{\klpo}(X)$ (with $Q$ denoting the dimension of this vector space).
The algebraic counterpart $\MAT{I}_{\ell,X}^{\ell+1}$ of the local restriction operator $\IntOp_{\ell,X}^{\ell+1}$ defined by \eqref{def:restOp} reads
$$
\MAT{I}_{\ell,X}^{\ell+1} =\begin{pmatrix}
\int_X \psi_i^{X,\ell+1} \psi_j^{X,\ell}
\end{pmatrix}_{i=1,\ldots,Q,\,j=1,\ldots,P},
$$
and the algebraic counterpart $\MAT{I}_{\ell+1}^{\ell,X}$ of the local prolongation operator $\IntOp^{\ell,X}_{\ell+1}$ is
\begin{equation*}
  \MAT{I}^{\ell,X}_{\ell+1} = \big(\MAT{I}_{\ell,X}^{\ell+1}\big)\trans.
\end{equation*}
Interestingly, when using hierarchical orthonormal bases and the basis for $\Poly{\klpo}(X)$ is obtained by restriction of the basis for $\Poly{\kl}(X)$, both the prolongation and restriction operators are represented by unit diagonal rectangular matrices.
In particular, for the local restriction operator it holds
\begin{equation*}
  \big(\MAT{I}_{\ell,X}^{\ell+1}\big)_{i,j}=\delta_{ij}, \qquad\text{for all $i=1,...,Q$ and all $j=1,...,P$.}
\end{equation*}
As a result, intergrid transfer operators need not need be computed nor stored in memory. 

With a little abuse of notation, we also denote by $\MAT{I}_{\ell,X}^{\ell+1}$ and $\MAT{I}_{\ell+1}^{\ell,X}$ the local restriction and prolongation operators applied to vector-valued variables, which are obtained assembling component-wise the corresponding operators acting on scalar-valued variables.
The matrix $\MAT{A}_T^{\ell+1}$ discretizing the viscous term at level $\ell+1$ can be inherited from the corresponding matrix $\MAT{A}_T^\ell$ at the level $\ell$ applying the restriction operators block-wise (compare with \eqref{eq:AT.block.refined}):
\begin{equation*}\renewcommand{\arraystretch}{1.3}
  \MAT{A}_T^{\ell+1} = 
  \begin{bmatrix}
    \MAT{I}_{\ell,T}^{\ell+1} \,\MAT{A}^\ell_{TT}    \, \MAT{I}^{\ell,T}_{\ell+1}  & \MAT{I}_{\ell,T}^{\ell+1}  \, \MAT{A}^\ell_{TF_1}    \, \MAT{I}^{\ell,F_1}_{\ell+1} & \cdots & \MAT{I}_{\ell,T}^{\ell+1} \, \MAT{A}^\ell_{TF_N}    \, \MAT{I}^{\ell,F_N}_{\ell+1}\\
    \MAT{I}_{\ell,F_1}^{\ell+1} \,\MAT{A}^\ell_{F_1 T} \, \MAT{I}^{\ell_,T}_{\ell+1} & \MAT{I}_{\ell,F_1}^{\ell+1} \, \MAT{A}^\ell_{F_1 F_1} \, \MAT{I}^{\ell,F_1}_{\ell+1} & \cdots & \MAT{I}_{\ell,F_1}^{\ell+1} \, \MAT{A}^\ell_{F_1 F_N} \, \MAT{I}^{\ell,F_N}_{\ell+1}\\
    \vdots  &  \vdots & \ddots & \vdots \\
    \MAT{I}_{\ell,F_N}^{\ell+1} \,\MAT{A}^\ell_{F_N T} \, \MAT{I}^{\ell,T}_{\ell+1}  & \MAT{I}_{\ell,F_N}^{\ell+1}   \, \MAT{A}^\ell_{F_N F_1} \, \MAT{I}^{\ell,F_1}_{\ell+1} & \cdots & \MAT{I}_{\ell,F_N}^{\ell+1} \, \MAT{A}^\ell_{F_N F_N} \, \MAT{I}^{\ell,F_N}_{\ell+1}\\
  \end{bmatrix}
\end{equation*}
Applying this procedure recursively shows that, for any level $\ell\ge 1$, the matrix $\MAT{A}_T^{\ell}$ can be obtained from the fine matrix $\MAT{A}_T^0$.
Note that pre- and post-multiplication of the matrix blocks by the restriction and the prolongation operators, respectively, results in a block shrink.  
When using orthonormal basis functions, these matrix multiplications can be avoided altogether and replaced with inexpensive sub-block extractions.

In order to further reduce the computational costs, Galerkin projections can be performed on the statically condensed fine grid operator, so that static condensation of coarse grid operators is avoided altogether.
For example, having computed the fine-level block of the Schur complement $\MAT{S}_{\partial T\partial T}^0$ (given by either formula \eqref{eq:ST.1} or \eqref{eq:ST.2}), the corresponding block $\MAT{S}_{\partial T\partial T}^{\ell+1}$ at level $\ell+1$ is computed applying recursively the relation:
\begin{equation*}
  \MAT{S}^{\ell+1}_{\partial T \, \partial T} =
  \begin{bmatrix}
    \MAT{I}_{\ell,F_1}^{\ell+1} \, \MAT{S}^{\ell}_{F_1 F_1} \, \MAT{I}^{\ell,F_1}_{\ell+1} & \cdots & \MAT{I}_{\ell,F_1}^{\ell+1} \, \MAT{S}^{\ell}_{F_1 F_N} \, \MAT{I}^{\ell,F_N}_{\ell+1}\\
    \vdots & \ddots & \vdots \\
    \MAT{I}_{\ell,F_N}^{\ell+1} \, \MAT{S}^{\ell}_{F_N F_1} \, \MAT{I}^{\ell,F_1}_{\ell+1} & ... & \MAT{I}_{\ell,F_N}^{\ell+1} \, \MAT{S}^{\ell}_{F_N F_N} \, \MAT{I}^{\ell,F_N}_{\ell+1}
  \end{bmatrix}.
\end{equation*}
To conclude the resulting sub-blocks are assembled into the global matrix. 


\section{Numerical results}
\label{sec:numResults}
\subsection{Mesh sequences}
\label{sec:meshSeq}
\begin{figure}
\centering
\begin{tabular}{cccc}
\includegraphics[width=0.20\textwidth]{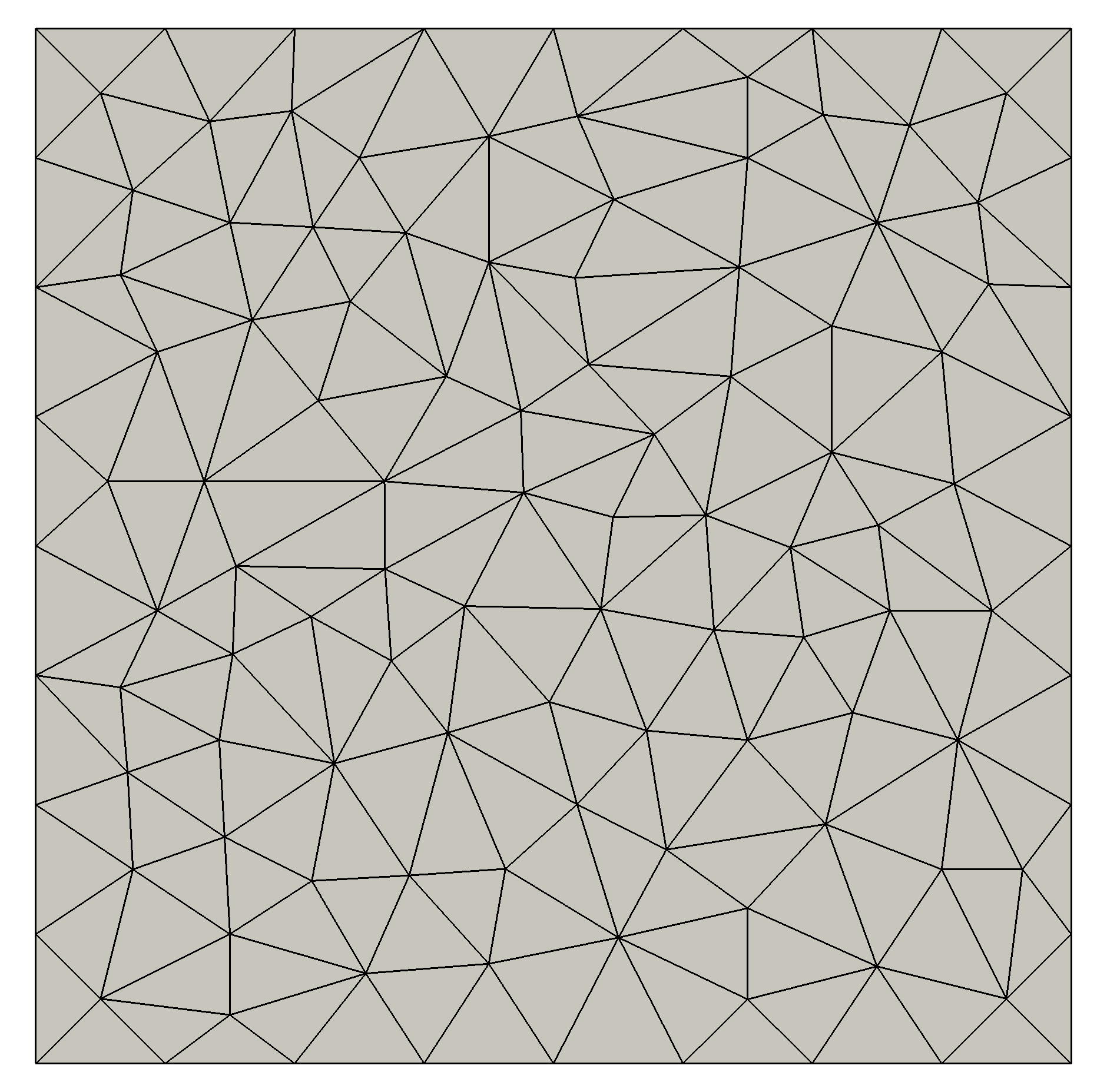} & 
\includegraphics[width=0.20\textwidth]{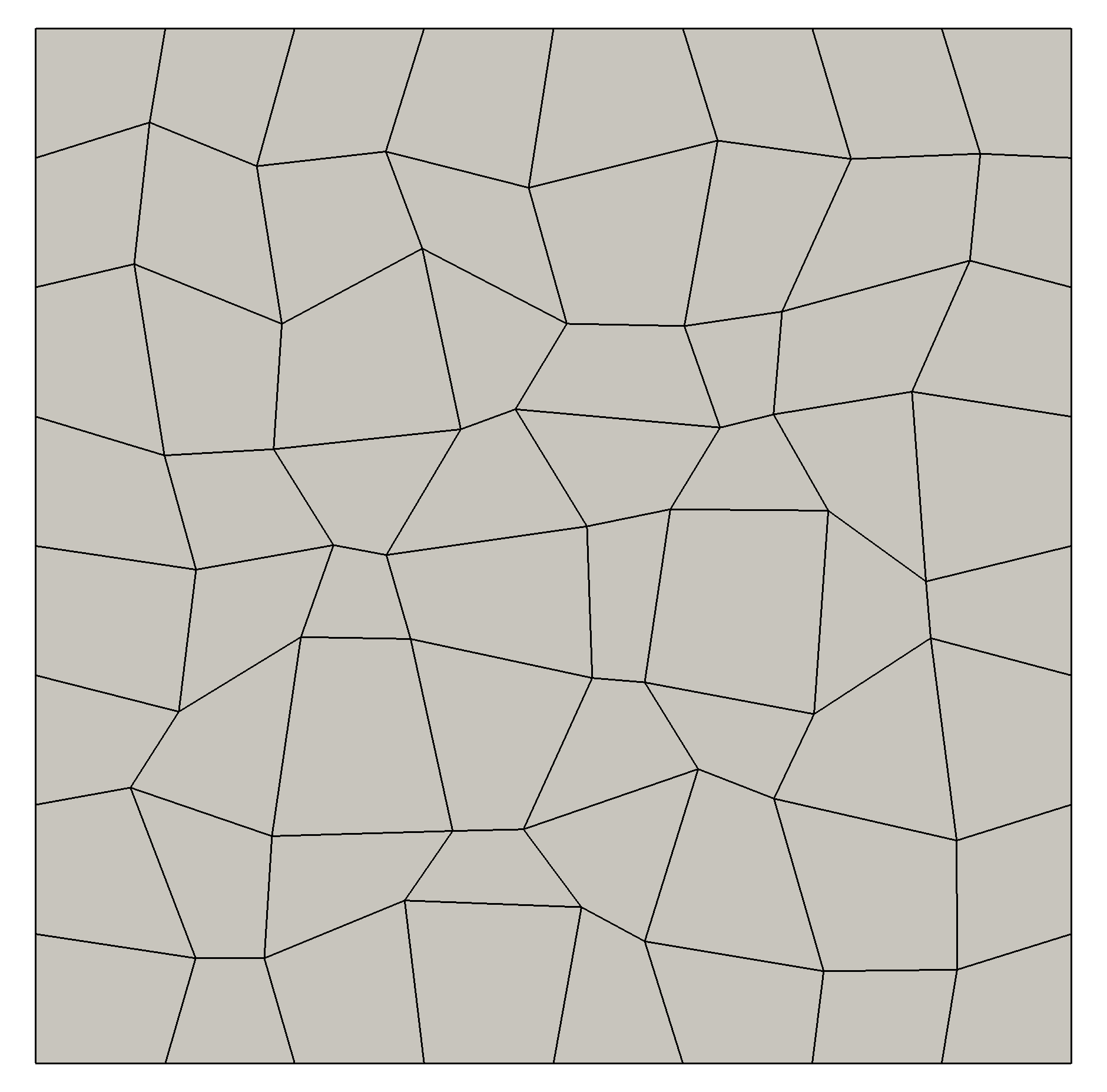}  & 
\includegraphics[width=0.20\textwidth]{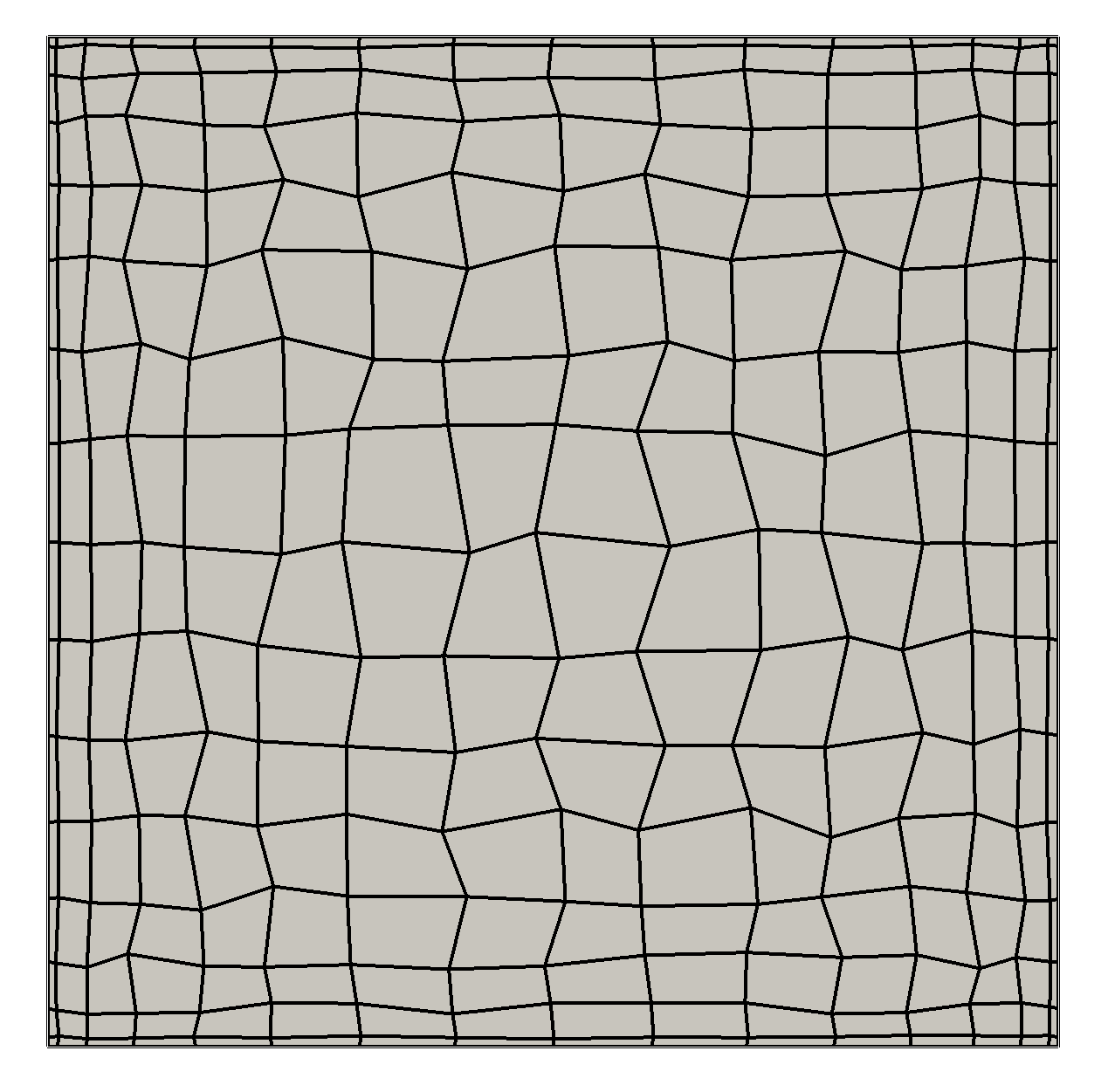} & 
\includegraphics[width=0.20\textwidth]{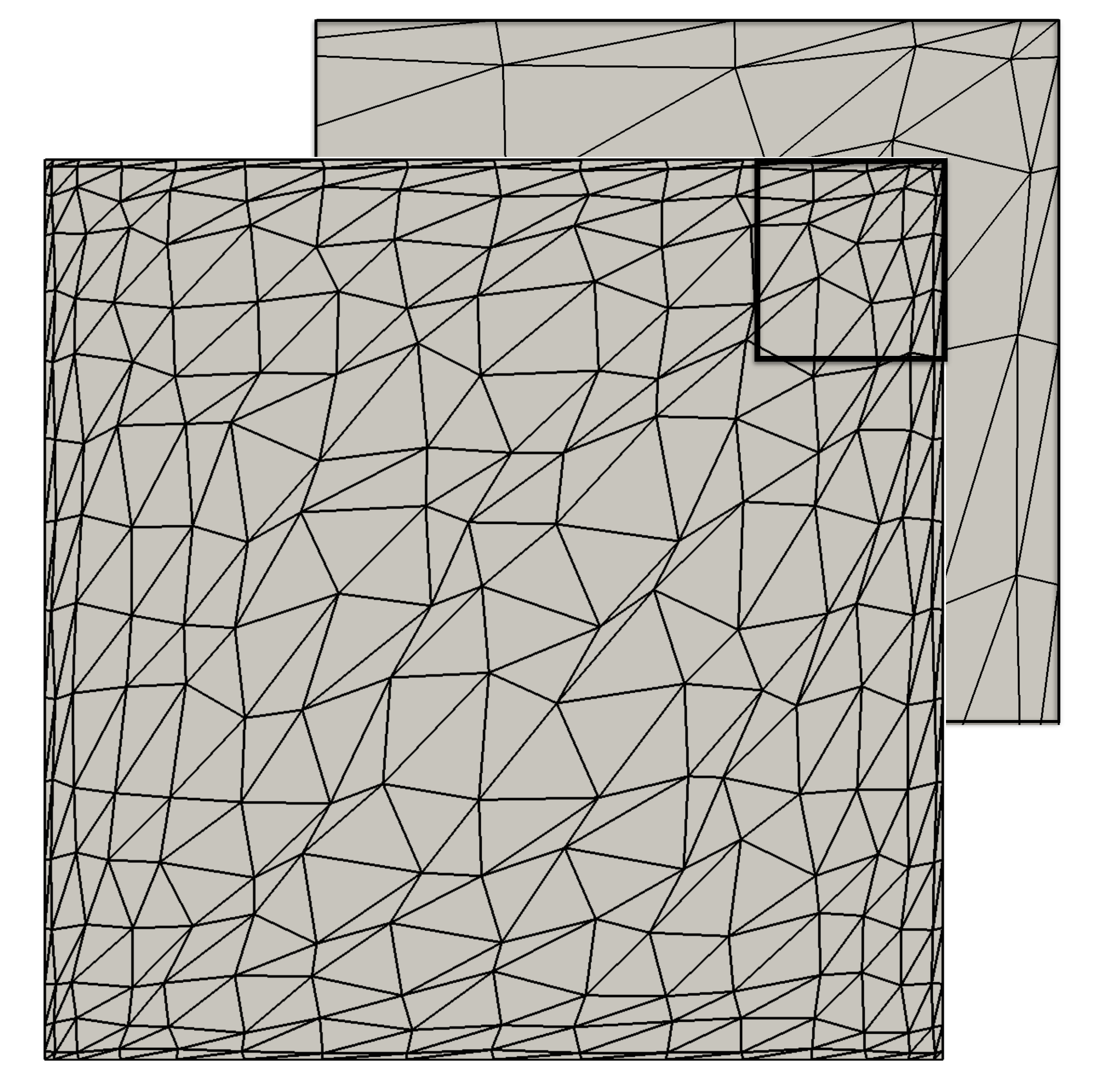} \\ 
\end{tabular}
\caption{Two-dimensional meshes (one mesh for each $h$-refined mesh sequence here considered). 
         From left to right: Delaunay triangular mesh, trapezoidal mesh, graded trapezoidal mesh, graded triangular mesh. \label{meshes2D}}
\end{figure}

\begin{figure}
\centering
\begin{tabular}{ccc}
\includegraphics[width=0.26\textwidth]{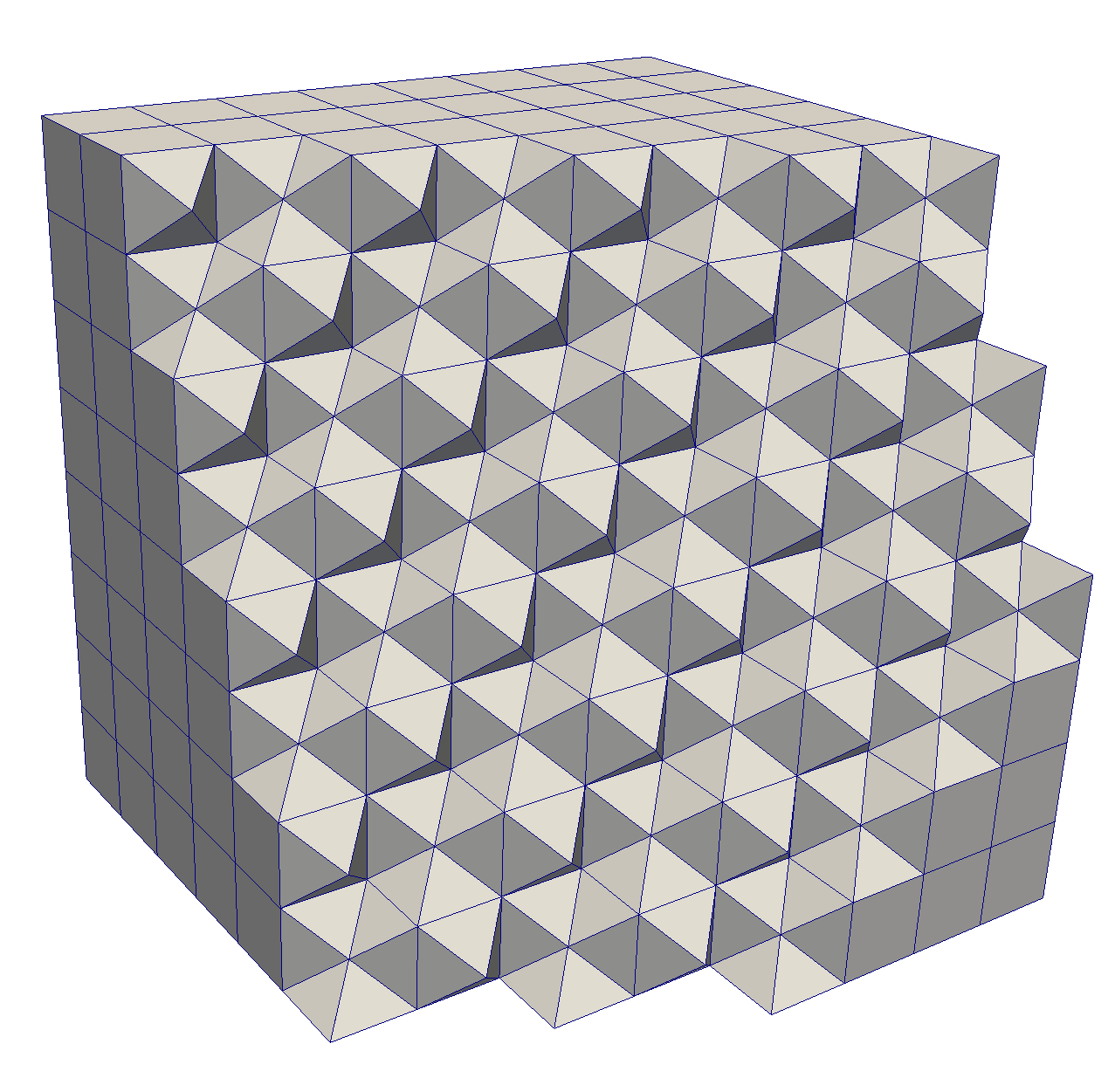} & 
\includegraphics[width=0.26\textwidth]{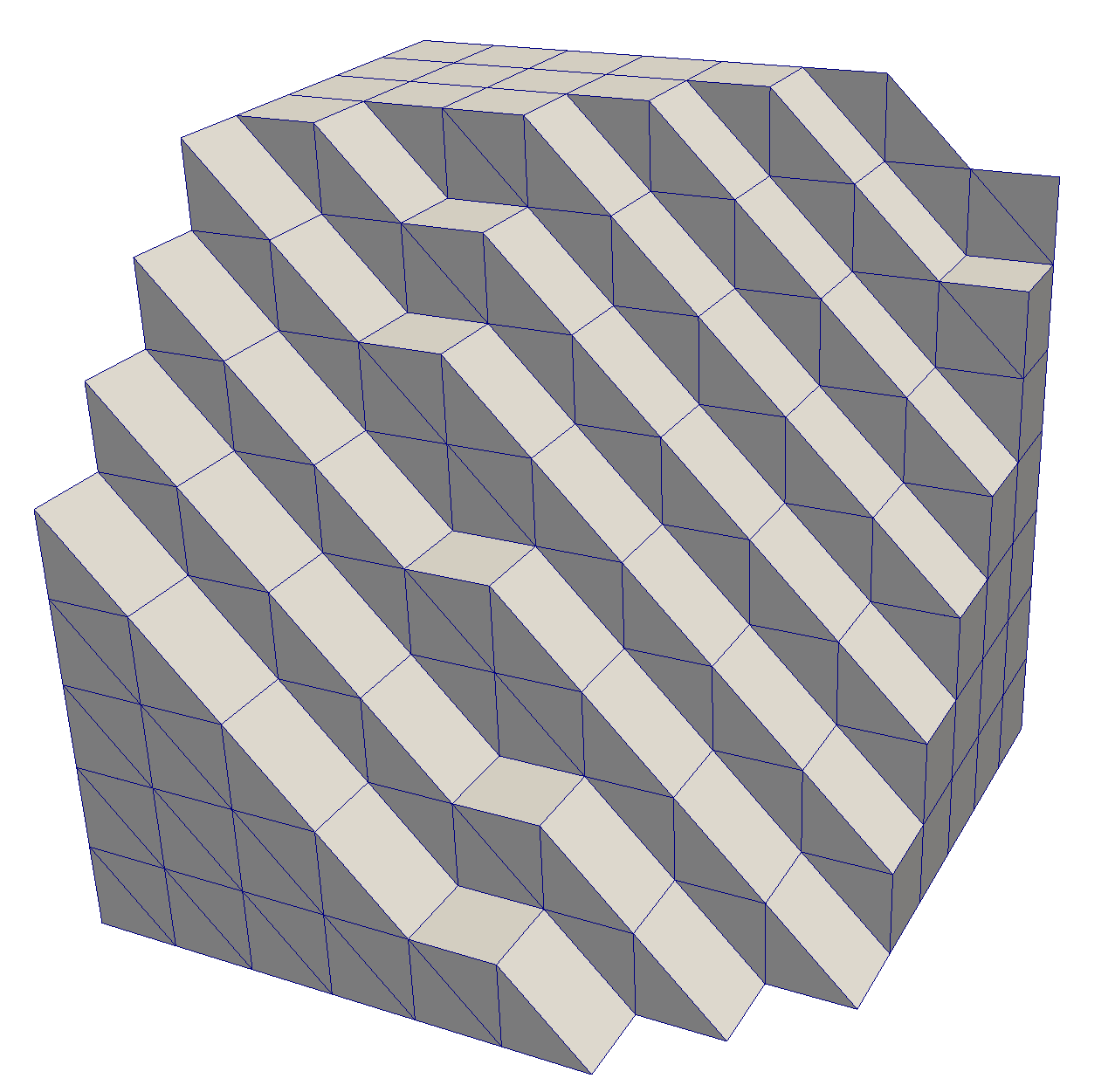} & 
\includegraphics[width=0.26\textwidth]{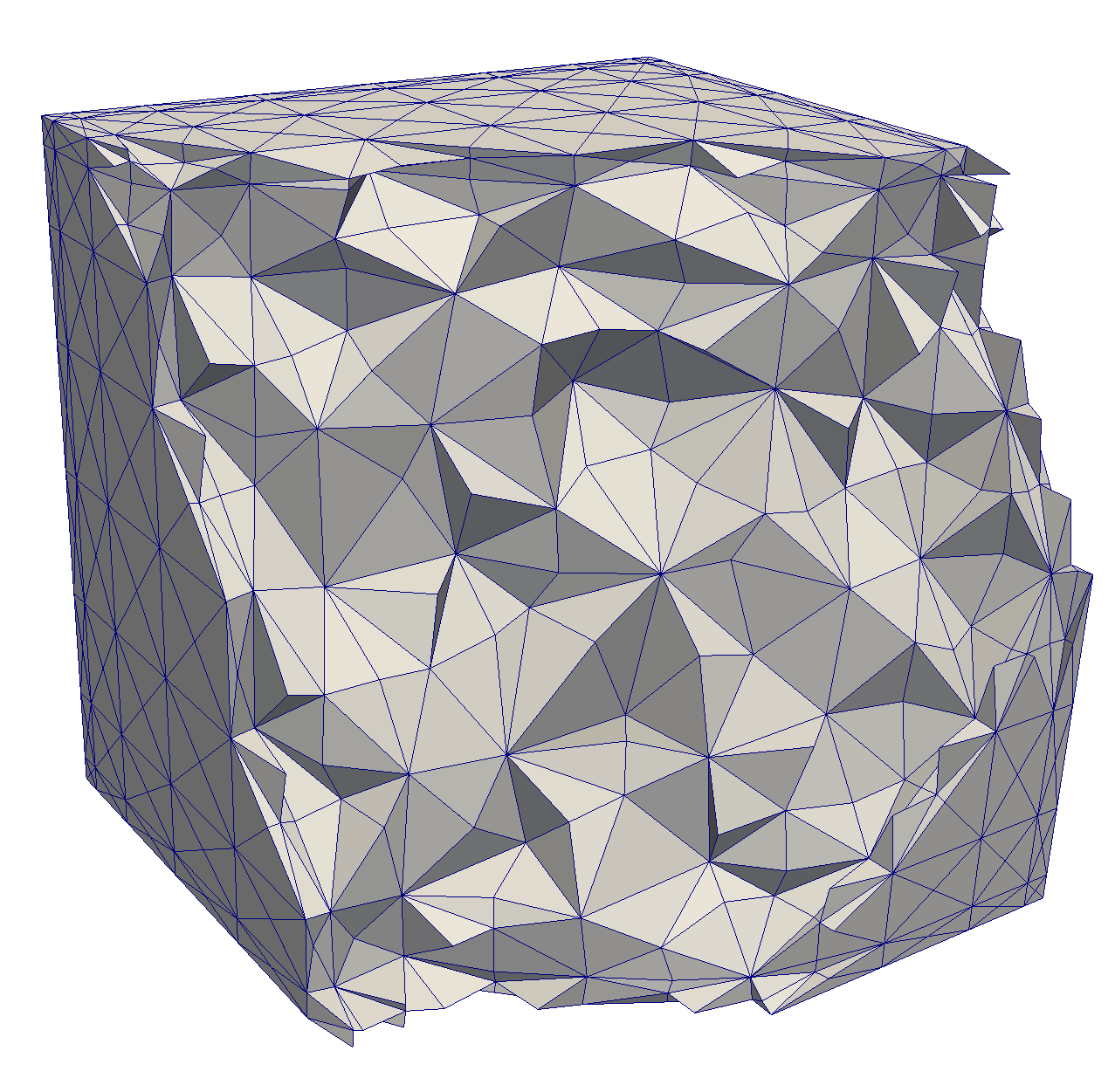} \\ 
\end{tabular}
\caption{Three-dimensional meshes (one mesh for each $h$-refined mesh sequence here considered). 
         From left to right: pyramidal mesh, prismatic mesh, graded tetrahedral mesh. \label{meshes3D}}
\end{figure}

In order to assess and compare the performance of $p$-multilevel preconditioners, we consider four $h$-refined mesh sequences of the two-dimensional domain $(-1,1)^2$, see Figure \ref{meshes2D}, and three $h$-refined mesh sequences of the three-dimensional domain $(0,1)^3$, see Figure \ref{meshes3D}.
In two space dimensions, we consider both standard and graded meshes composed of triangular and trapezoidal elements.
In three space dimensions, we consider standard meshes composed of prismatic and pyramidal elements and graded meshes composed of tetrahedral elements. 
While standard meshes have homogeneous meshsize, graded meshes feature mesh elements that become narrower and narrower while approaching 
the domain boundaries, mimicking computational grids commonly employed in CFD to capture boundary layers.  
In order to build $h$-refined graded mesh sequences, the mesh nodes are first positioned according to Gauss--Lobatto quandrature rules of increasing order and then randomly displaced by a small fraction of their distance.
Accordingly, the reduction of the meshsize is non-linear in case of graded $h$-refined mesh sequences.  

\subsection{Setting}

\subsubsection{Manufactured analytical solution}

We consider the following smooth analytical behaviours of the velocity and pressure fields:
If $d=2$, we let $\Omega\coloneq(-1,1)^2$ and set
\begin{equation*}
  \begin{alignedat}{2}
    \boldsymbol{u}(x,y) &= - e^x \, \left[
      y \, \cos(y) + \sin(y)
      \right] \, \boldsymbol{i} + e^x \, (y \, \sin(y)) \, \boldsymbol{j} &\qquad& \forall (x,y)\in\Omega,
    \\ 
    p(x,y) &= 2 \, e^x \, \sin(y), &\qquad& \forall(x,y)\in\Omega,
  \end{alignedat}
\end{equation*}
where $\{\boldsymbol{i},\boldsymbol{j}\}$ is the canonical basis of $\Real^2$ while, for $d=3$, we set $\Omega\coloneq(0,1)^3$ and
\begin{equation*}
  \begin{alignedat}{2}
    \boldsymbol{u}(x,y,z) &=  2 \, \sin(\pi \, x)\boldsymbol{i} - \pi \, y  \, \cos(\pi \, x) \boldsymbol{j}  -\pi \, z \, \cos(\pi \, x) \boldsymbol{k} &\qquad& \forall(x,y,z)\in\Omega,
    \\
    p(x,y,z) &= \sin(\pi \, x) \; \cos(\pi \, y) \; \sin(\pi \, z) &\qquad&\forall(x,y,z)\in \Omega,
  \end{alignedat}
\end{equation*}
where $\{\boldsymbol{i},\boldsymbol{j},\boldsymbol{k}\}$ is the canonical basis of $\Real^3$.
Dirichlet boundary conditions are enforced on all but one faces of $\Omega$, where Neumann boundary conditions are enforced instead.
The boundary data and forcing term are inferred from the exact solution.

\subsubsection{Multilevel solver options}
\label{sec:numResOptions}

We consider high-order and higher-order versions of the \hho{} and \dg{} schemes corresponding to the polynomial degrees $k {=} 3$ and $k {=} 6$, respectively. 
The theoretical $h$-convergence rates for \dg{} are $k{+}1$ for the velocity error in the $L^2$-norm and $k$ for the velocity gradient and the pressure error in the $L^2$-norm.
The theoretical $h$-convergence rates for \hho{} are $k{+}2$ for the velocity reconstruction error in the $L^2$-norm 
and $k{+}1$ for the gradient of the velocity reconstruction and the pressure error in the $L^2$-norm.
For the \hhohp{} scheme, both the element velocity and the reconstructed velocity display the same convergence rates, but the former is additionally divergence free on standard meshes.
For this reason, the element velocity field is used in the error computations.
For all the numerical test cases, we report in the tables the $L^2$-errors on the velocity (``$\uVec_h$'' column), velocity gradients (``$G\uVec_h$'' column), pressure (``$p_h$'' column), and divergence (``$D\uVec_h$'' column).

The solution of the linear systems is based on a FGMRES iterative solver preconditioned with a $p$-multilevel $V$-cycle iteration of three levels ($L=2$):
for $k_0{=}k{=}3$ (fine level), we set $k_1{=}2$ on the intermediate level and $k_2{=}k_L{=}1$ on the coarse level;
for $k_0{=}k{=}6$ (fine level), we set $k_1{=}3$ on the intermediate level and $k_2{=}k_L{=}1$ on the coarse level.

On the fine and intermediate levels, the pre- and post-smoothing strategy consist in two iterations of ILU preconditioned GMRES. 
On the coarse level, we employ an LU solver when working in two space dimensions and ILU preconditioned GMRES solver when working in three space dimensions. 
Since enforcing looser tolerances on the coarse level does not alter the number of outer FGMRES iterations,
we impose a three orders of magnitude decrease of the true (unpreconditioned) relative residual in three space dimensions.
The relative residual decrease for the outer FGMRES solver is set to $10^{-13}$ when $k{=}3$ and to $10^{-14}$ when $k{=}6$.

\subsubsection{Performance evaluation}

For all the numerical test cases we compare the performance and efficiency of solver strategies 
based on:
\begin{compactitem}
\item Number of FGMRES outer iterations (``ITs'' column);
\item Number of coarse solver iterations (``ITs$_L$'' column). Note that one iterations means that a direct solver is employed;
\item Wall clock time required for linear system solution (``CPU time Sol.'' column);
\item Wall clock time required for matrix assembly (``CPU time Ass.'' column);
      We remark that the computational cost of building the Schur complement is included since static condensation is performed element-by-element during matrix assembly. 
\item Wall clock time required for matrix assembly plus linear system solution (``CPU time Tot.'' column);
\item Efficiency with respect to linear scaling of the computational expense with respect to the number of DOFs (``Eff.'' column). 100\% efficiency means that for a fourfold increase of the number of DOFs we get a fourfold increase of the total (matrix assembly plus linear system solution) wall clock time.
\end{compactitem}

\subsection{Comparison based on matrix dimension and matrix non-zero entries}
\label{sec:DOFsMNZsComp}
The cost of a Krylov iteration scales linearly with the number of Matrix Non-Zero entries (MNZs) plus
the number of Krylov spaces times the matrix dimension (equal to the number of Degrees Of Freedom, DOFs), see, \eg, \cite{Quarteroni.Sacco.ea:2000}.
Multilevel Krylov solvers utilize only a few smoother iterations on the fine and intermediate levels  
and iteratively solve on the coarse level, where the number of MNZs and DOFs is favourable, see Section \ref{sec:MGViter}.
Accordingly, with respect to solver efficiency, the most relevant discretization-dependent parameters are the MNZs of the fine and coarse matrices and the number of  and DOFs of the coarse level:
fine level MNZs influence the cost of the most expensive matrix-vectors products, performed once per smoother iteration; 
coarse level MNZs influence the cost of the least expensive matrix-vector products, performed once per iteration of the coarse solver (that is, many times per multilevel iteration); 
the number of DOFs of the coarse level influences the cost of the Gram--Schmidt orthogonalization carried out within the GMREs algorithm on the coarse level.

Static condensation of the element based unknowns is an effective means of improving solver efficiency in the context of hybridized methods.
For \hhodp{}, we compare the uncondensed (\hhodp{} \uncond) implementation to the static condensation strategies described in Section \ref{sec:staticCondHHO}.
We recall that both static condensation procedures involve the local elimination of velocity unknowns attached to mesh elements, and the difference lays in the treatment of pressure degrees of freedom.
According to \eqref{eq:ST.2}, all pressure modes except the constant value are statically condensed in the \hhodp{} \vpcond{} strategy, while, according to \eqref{eq:ST.1}, pressure modes are not statically condensed in the \hhodp{} \vcond{} strategy.
For \hhohp{}, we consider static condensation of the element unknowns for both the velocity and the pressure (\hhohp{} \vpcond), so that the only skeletal unknowns appear in the global systems.

Roughly speaking, DOFs and MNZs of HHO discretizations are associated with 
\emph{element variables} and \emph{face variables}.
DG discretizations rely only on element variables.
The formulas for computing DOFs and MNZs reported in Table \ref{tab:DOFsMNZsFormula} show that:
\begin{compactitem}
\item The number of DOFs associated with element variables is proportional to the dimension of the polynomial space $\Polyd{d}{k}$ and to the number of mesh elements;
\item The number of DOFs associated with face variables is proportional to the dimension of $\Polyd{d-1}{k}$ and to the number of mesh faces;
\item The number of MNZs associated with element variables is proportional to the square of the dimension of $\Polyd{d}{k}$ and to the number of mesh elements;
\item The number of MNZs associated with face variables is proportional to the square of the dimension of $\Polyd{d-1}{k}$ and to the number of mesh faces.
\end{compactitem}
MNZs are also influenced by the stencil of the discretization and the fill-in of the Schur complement, as explained in Remark \ref{rem:StatCond}.
Since the ratio between the dimensions of $\Polyd{d}{k}$ and $\Polyd{d-1}{k}$ is $\frac{k+d}{d}$, we have the following rules of thumb:
\begin{equation}
\label{eq:dofmnzCond}
\begin{aligned}
&\text{face variables have fewer DOFs than element variables if $ \frac{\card{\Fh}}{\card{\Th}} {<} \frac{k+d}{d}$},  \\
&\text{face variables have fewer MNZs than element variables if $ \frac{\card{\Fh}}{\card{\Th}} \frac{2\card{\TF}{-}1}{\card{\TF}{+}1} {<} \left(\frac{k+d}{d}\right)^2$},
\end{aligned}
\end{equation}
where $\frac{2\card{\TF}{-}1}{\card{\TF}{+}1}$ is the ratio between the stencil of face variables and element variables, respectively. 
This simple observation allows to interpret the results of Tables \ref{tab:dofs2D}--\ref{tab:dofs3D} and \ref{tab:mnzs2D}--\ref{tab:mnzs3D}, where the DOFs and MNZs counts for the methods and implementations considered in this work are reported.
Placeholders correspond to combinations of meshes, polynomial degrees, schemes, and static condensation options that are either not possible or haven't been considered in numerical tests.
The data are reported only for the finest grids of each mesh sequence for $k\in\{1,3,6\}$ (the case $k{=}1$ is also included as it is relevant for estimating the efficiency of the coarse solver).

Some comments regarding the DOFs counts reported in Tables \ref{tab:dofs2D}--\ref{tab:dofs3D} are as follows.
As expected, the \hhodp{} \uncond{} DOFs count is the largest.
In 2D and 3D, \hhodp{} \vpcond{} and \dg{}, respectively, have the fewest DOFs count on the coarse level ($k{=}1$).
This can easily be interpreted based on \eqref{eq:dofmnzCond}, as the condition is harder to meet in 3D than in 2D. 
In 2D, the number of coarse level DOFs for \hhodp{} \vcond, \hhohp{}, and \dg{} are very similar.
In 2D and 3D, higher-order statically condensed \hho{} shows some advantage over \dg{} in terms of DOFs. 

Some comments regarding the MNZs counts reported in Tables \ref{tab:dofs2D}--\ref{tab:dofs3D} are as follows.
In 2D, \hhodp{} \vpcond{} and \vcond{} have fewer MNZs than \dg{}, at all polynomial degrees. 
In 3D, \hhodp{} \vpcond{} and \vcond{} have fewer MNZs than \dg{} for both $k{=}3$ and $k{=}6$, with \hhodp{} \vcond{} being the most efficient. 
\hhodp{} \vcond{} is very close to \dg{} for $k{=}1$.
The fact that \hhodp{} \vcond{} outperforms \hhodp{} \vpcond{} is due to increased fill-in of the Schur complement matrix arising from \eqref{eq:ST.2}, see Remark \ref{rem:StatCond}. 
\hhohp{} \vpcond{} improves \dg{} only for $k{=}6$, while \dg{} is significantly better for both $k{=}1$ and $k{=}3$.
Similar to strategy \eqref{eq:ST.2} for \hhodp, the aforementioned static condensation procedure increases the fill-in of the blocks pertaining to skeletal velocity unknowns with respect to the uncondensed operator.

\begin{table}\centering
  \begin{footnotesize}
    \renewcommand{\arraystretch}{1.4}
    \begin{tabular}{c|c|c}
      \toprule
      Scheme  &   Matrix dim. (DOFs) & Matrix non-zero entries (MNZs) \\
      \midrule
      \begin{tabular}{c}\hhodp{} \\ \uncond{} \end{tabular}    & 
      $\begin{array}{c}\card{\Th} (d{+}1) \dim(\Polyd{d}{k}) +
        \\ \card{\Fh} d \dim(\Polyd{d{-}1}{k}) \end{array}$ &  $\begin{array}{c}
        \card{\Th} (d{+}1) \dim(\Polyd{d}{k})^2 {+}
        \sum_{T \in \Th} \card{\FT} d^2 \dim(\Polyd{d}{k}) \dim(\Polyd{d{-}1}{k}){+}\\
        \sum_{F \in \Fh} \card{\TF} d^2 \dim(\Polyd{d}{k}) \dim(\Polyd{d{-}1}{k})+
        \sum_{F \in \Fh} \left(2\card{\TF}{-}1\right) d \dim(\Polyd{d{-}1}{k})^2 \\
      \end{array}$\\                                           
      \midrule
      \begin{tabular}{c}\hhodp{}  \\ \vpcond\end{tabular} & 
        $\begin{array}{c}\card{\Th} + \\\card{\Fh} d \dim(\Polyd{d{-}1}{k}) \end{array}$& $\begin{array}{c}
          \card{\Th} +
          \sum_{T \in \Th} \card{\FT} d \dim(\Polyd{d}{k}) {+}\\
          \sum_{F \in \Fh} \left(2\card{\TF}{-}1\right) d^2 \dim(\Polyd{d{-}1}{k})^2 +
          \sum_{F \in \Fh} \card{\TF} d \dim(\Polyd{d{-}1}{k})\\
        \end{array}$\\
        \midrule
        \begin{tabular}{c}\hhodp{} \\ \vcond{}  \end{tabular}  & 
        $\begin{array}{c}\card{\Th} \dim(\Polyd{d}{k}) +\\ 
          \card{\Fh} (d{+}1) \dim(\Polyd{d{-}1}{k})\end{array} $ &   $\begin{array}{c}
          \card{\Th} \dim(\Polyd{d}{k})^2{+}
          \sum_{T \in \Th} \card{\FT} d \dim(\Polyd{d}{k}) \dim(\Polyd{d{-}1}{k}){+}\\
          \sum_{F \in \Fh} \card{\TF} d \dim(\Polyd{d}{k}) \dim(\Polyd{d{-}1}{k}){+}
          \sum_{F \in \Fh} \left(2\card{\TF}{-}1\right) d \dim(\Polyd{d{-}1}{k})^2\\
        \end{array}$\\                           
        \midrule
        \begin{tabular}{c}\hhohp{} \\ \vpcond\end{tabular} & 
          $\card{\Fh} (d{+1}) \dim(\Polyd{d{-}1}{k})$  & $\sum_{F \in \Fh} \left(2\card{\FT}{-}1\right) \left[(d{+}1)^2\right] \dim(\Polyd{d{-}1}{k})^2$ \\
          \midrule
          \dg{}               &   $\card{\Th} (d{+}1) \dim(\Polyd{d}{k})$  & $\sum_{T \in \Th} \left(\card{\FT}{+}1\right) (3d{+}1) \dim(\Polyd{d}{k})^2$ \\
          \bottomrule
    \end{tabular}
  \end{footnotesize}
  \caption{Formulas for computing the matrix dimension (equal to the number of DOFs) and number of matrix non-zero entries (MNZs) of nonconforming discretization. Several static condensation options are considered for \hhodp, see text for details. \label{tab:DOFsMNZsFormula}}
\end{table}

\begin{table}
  \centering
  \begin{tabular}{ c c c | c c c | c | c }
    \multicolumn{8}{c}{Number of DOFs, 2D meshes} \\
    \toprule
    mesh seq.                  & $\card{\Th}$                    &   $k$     & \multicolumn{3}{c|}{\hhodp{}}                & \hhohp{}        & \dg{}        \\
    &                                 &           & \uncond{}         & \vpcond{}    & \vcond{}     & \vpcond{}     &              \\
    \midrule
    \multirow{2}{*}{trapz-quad}& \multirow{2}{*}{16384}          & $k{=}3$   & 7.56e+05       & 2.81e+05     & 4.28e+05   & 3.96e+05      & 4.92e+05     \\

    &                                 & $k{=}1$   & 2.8e+05        & 1.48e+05     & 2.23e+05   & 1.98e+05      & 1.47e+05     \\
    &  \multirow{2}{*}{1024}          & $k{=}6$   &    -           & 3.06e+04     & 5.82e+04   & 4.44e+04      & 8.6e+04      \\
    &                                 & $k{=}1$   &    -           & 9.47e+03     & 1.15e+04   & 1.27e+04      & 9.22e+03     \\
    \multirow{2}{*}{Del. tri}  & \multirow{2}{*}{50744}          & $k{=}3$   & 2.13e+06       & 6.62e+05     & 1.12e+06   & 9.16e+05      & 1.52e+06     \\
    &                                 & $k{=}1$   & 7.62e+05       & 3.56e+05     & 4.58e+05   & 4.58e+05      & 4.57e+05     \\
    &  \multirow{2}{*}{3120}          & $k{=}6$   &    -           & 6.95e+04     & 1.54e+05   & 9.96e+04      & 2.62e+05     \\
    &                                 & $k{=}1$   &    -           & 2.21e+04     & 2.83e+04   & 2.85e+04      & 2.81e+04     \\
    \multirow{2}{*}{dist. tri} & \multirow{2}{*}{32768}          & $k{=}3$   & 1.38e+06       & 4.28e+05     & 7.23e+05   & 5.93e+05      & 9.83e+05     \\
    &                                 & $k{=}1$   & 4.93e+05       & 2.3e+05      & 2.96e+05   & 2.96e+05      & 2.95e+05     \\
    &  \multirow{2}{*}{2048}          & $k{=}6$   &    -           & 4.6e+04      & 1.01e+05   & 6.59e+04      & 1.72e+05     \\
    &                                 & $k{=}1$   &    -           & 1.46e+04     & 1.87e+04   & 1.88e+04      & 1.84e+04     \\
    \bottomrule
  \end{tabular}
  \caption{Matrix dimension (equal to the number of DOFs) on 2D meshes. Note that graded quadrilateral meshes are not included because they coincide with trapezoidal meshes in terms of DOFs.\label{tab:dofs2D}}
\end{table}

\begin{table}
\centering
\begin{tabular}{ c c c | c c c | c | c }
\multicolumn{8}{c}{Number of MNZs, 2D meshes} \\
\toprule
mesh seq.                  & $\card{\Th}$                    &   $k$   & \multicolumn{3}{c|}{\hhodp{}}                   & \hhohp{}        & \dg{}              \\
&                                 &           & \uncond{}         & \vpcond{}    & \vcond{}      & \vpcond{}     &                  \\
\midrule
\multirow{4}{*}{trapz-quad}& \multirow{2}{*}{16384}          & $k{=}3$   & 3.98e+07       & 1.58e+07     & 2.68e+07    & 3.31e+07      & 5.7e+07          \\
                           &                                 & $k{=}1$   & 6.01e+06       & 4.21e+06     & 3.56e+06    & 8.27e+06      & 5.13e+06         \\
                           &  \multirow{2}{*}{1024}          & $k{=}6$   &   -            & 2.94e+06     & 1.21e+07    & 6.35e+06      & 2.74e+07        \\
                           &                                 & $k{=}1$   &   -            & 2.64e+05     & 2.23e+05    & 5.18e+05      & 3.14e+05        \\
\multirow{4}{*}{Del. tri}  & \multirow{2}{*}{50744}          & $k{=}3$   & 9.64e+07       & 2.69e+07     & 5.38e+07    & 5.48e+07      & 1.42e+08        \\
                           &                                 & $k{=}1$   & 1.36e+07       & 7.36e+06     & 7.16e+06    & 1.37e+07      & 1.28e+07        \\
                           &  \multirow{2}{*}{3120}          & $k{=}6$   &   -            & 4.86e+06     & 5.43e+06    & 1.03e+07      & 6.78e+07        \\
                           &                                 & $k{=}1$   &   -            & 4.53e+05     & 4.4e+05     & 8.45e+05      & 7.78e+05        \\
\multirow{4}{*}{dist. tri} & \multirow{2}{*}{32768}          & $k{=}3$   & 6.23e+07       & 1.74e+07     & 3.48e+07    & 3.54e+07      & 9.14e+07         \\
                           &                                 & $k{=}1$   & 8.75e+06       & 4.76e+06     & 4.62e+06    & 8.86e+06      & 8.23e+06         \\
                           &  \multirow{2}{*}{2048}          & $k{=}6$   &   -            & 3.2e+06      & 7.93e+06    & 6.8e+06       & 4.43e+07        \\
&                                 & $k{=}1$   &   -            & 2.98e+05     & 2.89e+05    & 5.55e+05      & 5.08e+05         \\
\bottomrule
\end{tabular}
\caption{Number of matrix non-zero entries (MNz) on 2D meshes. Note that graded quadrilateral meshes are not included because they coincide with trapezoidal meshes in terms of MNZs.\label{tab:mnzs2D}}
\end{table}

\begin{table}
\centering
\begin{tabular}{ c c c | c c c | c | c }
\multicolumn{8}{c}{Number of DOFs, 3D meshes} \\
\toprule
mesh seq.                  & $\card{\Th}$                    &   $k$   & \multicolumn{3}{c|}{\hhodp{}}          & \hhohp{}      & \dg{}         \\
&                                 &           & \uncond{}     & \vpcond{} & \vcond{}    & \vpcond{}   &             \\
\midrule
\multirow{4}{*}{dist. tet} & \multirow{2}{*}{12288}          & $k{=}3$   & 1.74e+06   & 7.73e+05  & 1.01e+06  & 1.01e+06    & 9.83e+05    \\
                           &                                 & $k{=}1$   & 4.25e+05   & 2.4e+05   & 2.77e+05  & 3.04e+05    & 1.97e+05    \\
                           &  \multirow{2}{*}{1536}          & $k{=}6$   &  -         &  -        & 4.03e+05  & 3.66e+05    & 5.16e+05    \\
                           &                                 & $k{=}1$   &  -         &  -        & 3.55e+04  & 3.92e+04    & 2.46e+04    \\
\multirow{4}{*}{prism.}    & \multirow{2}{*}{8192}           & $k{=}3$   & 1.3e+06    & 6.53e+05  & 8.09e+05  & 8.6e+05     & 6.55e+05    \\
                           &                                 & $k{=}1$   & 3.25e+05   & 2.02e+05  & 2.26e+05  & 2.58e+05    & 1.31e+05    \\
                           &  \multirow{2}{*}{1024}          & $k{=}6$   &  -         &  -        & 3.23e+05  & 3.15e+05    & 3.44e+05    \\
                           &                                 & $k{=}1$   &  -         &  -        & 2.94e+04  & 3.38e+04    & 1.64e+04    \\
\multirow{4}{*}{pyram}     & \multirow{2}{*}{24576}          & $k{=}3$   & 3.83e+06   & 1.89e+06  & 2.36e+06  & 2.49e+06    & 1.97e+06    \\
                           &                                 & $k{=}1$   & 9.53e+05   & 5.84e+05  & 6.58e+05  & 7.46e+05    & 3.93e+05    \\
                           &  \multirow{2}{*}{3072}          & $k{=}6$   &  -         &  -        & 9.19e+05  & 8.82e+05    & 1.03e+06    \\
&                                 & $k{=}1$   &  -         &  -        & 8.31e+04  & 9.45e+04    & 4.92e+04    \\
\bottomrule
\end{tabular}
\caption{Matrix dimension (equal to the number of DOFs) on 3D meshes. \label{tab:dofs3D}}
\end{table}

\begin{table}
\centering
\begin{tabular}{ c c c | c c c | c | c }
\multicolumn{8}{c}{Number of MNZs, 3D meshes} \\
\toprule
mesh seq.                  & $\card{\Th}$                    &   $k$   & \multicolumn{3}{c|}{\hhodp{}}          & \hhohp{}     & \dg{}       \\
&                                 &           & \uncond{}     & \vpcond{} & \vcond{}    & \vpcond{}  &            \\
\midrule
\multirow{4}{*}{dist. tet} & \multirow{2}{*}{12288}          & $k{=}3$   & 2.19e+08   & 1.58e+08  & 1.16e+08  &  2.76e+08  & 2.4e+08   \\
                           &                                 & $k{=}1$   & 1.37e+07   & 1.49e+07  & 8.4e+06   &  2.49e+07  & 9.58e+06  \\
                           &  \multirow{2}{*}{1536}          & $k{=}6$   & -          &   -       & 1.49e+08  &  2.72e+08  & 5.15e+08  \\
                           &                                 & $k{=}1$   &  -         &  -        & 1.05e+06  &  3.12e+06  & 1.17e+06  \\
\multirow{4}{*}{prism.}    & \multirow{2}{*}{8192}           & $k{=}3$   & 1.87e+08   & 1.69e+08  & 1.08e+08  &  2.97e+08  & 1.88e+08  \\
                           &                                 & $k{=}1$   & 1.22e+07   & 1.58e+07  & 8.08e+06  &  2.67e+07  & 7.54e+06  \\
                           &  \multirow{2}{*}{1024}          & $k{=}6$   & -          &   -       & 1.34e+08  &  2.92e+08  & 3.97e+08  \\
                           &                                 & $k{=}1$   &  -         &  -        & 1.01e+06  &  3.35e+06  & 9.01e+05  \\
\multirow{4}{*}{pyram}     & \multirow{2}{*}{24576}          & $k{=}3$   & 5.59e+08   & 5.06e+08  & 3.23e+08  &  8.86e+08  & 5.84e+08  \\
                           &                                 & $k{=}1$   & 3.66e+07   & 4.71e+07  & 2.42e+07  &  7.97e+07  & 2.33e+07  \\
                           &  \multirow{2}{*}{3072}          & $k{=}6$   & -          &   -       & 4.01e+08  &  8.69e+08  & 1.27e+09  \\
&                                 & $k{=}1$   &  -         &  -        & 3.03e+06  &  9.98e+06  & 2.89e+06  \\
\bottomrule
\end{tabular}
\caption{Number of matrix non-zero (MNZs) entries on 3D meshes. \label{tab:mnzs3D}}
\end{table}

\subsection{Comparison of static condensation strategies}
\label{sec:nrStatCond}
In this section we evaluate the performance of the multilevel solution strategy for Scheme \ref{scheme:hho} (\hhodp) comparing the two approaches for static condensation described in Section \ref{sec:staticCondHHO};
see in particular \eqref{eq:ST.2} (\hhodp{} \vpcond) and \eqref{eq:ST.1} (\hhodp{} \vcond).
We also consider the uncondensed formulation (\hhodp{} \uncond) as a reference to evaluate the performance gains.

In case of regular 2D mesh sequences, the results reported in Table \ref{tab:numResStatCondReg2D}
confirm that static condensation leads to significant gains (on average, the computation time halves) when compared with the uncondensed implementation. 
The results reported in Table \ref{tab:numResStatCondGrad2D}, where graded 2D mesh sequences are considered, 
show that the \hhodp{} \vpcond{} strategy (static condensation of both velocity element unknowns and high-order pressure modes) leads to a suboptimal performance 
of the multigrid preconditioner in case of stretched elements: notice the increasingly high number of FGMRES iterations when the mesh is refined.
A similar behavior, even if less pronounced, is observed for the uncondensed implementation.
The results reported in Table \ref{tab:numResStatCond3D}, where 3D mesh sequences are considered, confirm the strategy \hhodp{} \vcond{} (static condensation of element velocity unknowns only) leads to the best performance in terms of execution times, both in the case of standard and graded meshes. 
We remark that the gains are to be ascribed to fewer FGMRES iterations and a smaller number of
matrix non-zero entries, see Table \ref{tab:mnzs3D}. 

It is interesting to remark that accuracy and convergence rates are not influenced by the static condensation procedure 
as soon as the relative residual drop satisfies the prescribed criterion.
Solver fails to converge for \hho{} \vpcond{} over fine graded triangular meshes, see Table \ref{tab:numResStatCondGrad2D}.
Note that the prescribed maximum number of iteration (1k) of the FMGRES solver is reached and the convergence rates are spoiled.

\begin{sidewaystable}[ht]
\small
\centering
\begin{tabular}{ c | c c  c c | c c c | c c | c c c | c }
\multicolumn{4}{c}{{\cellcolor{BallBlue} \hhodp{} \uncond{} }} 
                                              & \multicolumn{5}{c}{{\cellcolor{LightBlue} \hhodp{} \vpcond{} }} 
                                              & \multicolumn{5}{c}{{\cellcolor{LightCyan} \hhodp{} \vcond{} }}\\
$\card{\Th}$ & \multicolumn{4}{c|}{error in $L^2$ norm} & \multicolumn{3}{c|}{conv. rate} & \multicolumn{2}{c|}{ITs} & \multicolumn{3}{c|}{CPU time} & Eff \\ \hline
 &  $\uVec_h$    & $G\uVec_h$   & $p_h$ & $D\uVec_h $& $\uVec_h$   & $G\uVec_h$ & $p_h$  & ITs & ITs$_L$ & Sol. & Ass. & Tot. & \\ \hline
&	\multicolumn{13}{l}{trapezoidal elements grid}\\
\hline
\rowcolor{BallBlue} 4      & 0.00201 & 0.0243 & 0.0117 & 0.0123 & - & - & -  & 10 & 1 & 0.00196 & 0.00108 & 0.00304 & -                \\ 
\rowcolor{BallBlue} 16     & 7.42e-05 & 0.00181 & 0.000951 & 0.000825 & 4.76 & 3.74 & 3.62 & 13 & 1 & 0.00722 & 0.00421 & 0.0114 & 106 \\ 
\rowcolor{BallBlue} 64     & 2.83e-06 & 0.00013 & 6.45e-05 & 6.3e-05 & 4.71 & 3.8 & 3.88 & 13 & 1 & 0.0285 & 0.0165 & 0.0451 & 101     \\ 
\rowcolor{BallBlue} 256    & 9.43e-08 & 8.43e-06 & 4.07e-06 & 4.25e-06 & 4.91 & 3.95 & 3.98 & 13 & 1 & 0.137 & 0.0663 & 0.203 & 88.8   \\ 
\rowcolor{BallBlue} 1024   & 2.95e-09 & 5.28e-07 & 2.65e-07 & 2.77e-07 & 5 & 4 & 3.94 & 13 & 1 & 0.672 & 0.265 & 0.937 & 86.7          \\ 
\rowcolor{BallBlue} 4096   & 9.4e-11 & 3.37e-08 & 1.71e-08 & 1.77e-08 & 4.97 & 3.97 & 3.96 & 13 & 1 & 3.01 & 1.07 & 4.07 & 92          \\ 
\rowcolor{BallBlue} 16384  & 2.99e-12 & 2.14e-09 & 1.09e-09 & 1.13e-09 & 4.97 & 3.98 & 3.97 & 13 & 1 & 13.6 & 4.32 & 17.9 & 91         \\ 
\rowcolor{LightBlue}      4      & 0.00201 & 0.0243 & 0.0117 & 0.0123 & - & - & -  & 4 & 1 & 0.000997 & 0.000756 & 0.00175 & -               \\
\rowcolor{LightBlue}      16     & 7.42e-05 & 0.00181 & 0.000951 & 0.000825 & 4.76 & 3.74 & 3.62 & 5 & 1 & 0.00225 & 0.00293 & 0.00518 & 135 \\
\rowcolor{LightBlue}      64     & 2.83e-06 & 0.00013 & 6.45e-05 & 6.3e-05 & 4.71 & 3.8 & 3.88 & 6 & 1 & 0.008 & 0.0115 & 0.0195 & 106       \\
\rowcolor{LightBlue}      256    & 9.43e-08 & 8.43e-06 & 4.07e-06 & 4.25e-06 & 4.91 & 3.95 & 3.98 & 6 & 1 & 0.0348 & 0.0459 & 0.0808 & 96.5  \\
\rowcolor{LightBlue}      1024   & 2.95e-09 & 5.28e-07 & 2.65e-07 & 2.77e-07 & 5 & 4 & 3.94 & 6 & 1 & 0.18 & 0.185 & 0.365 & 88.5            \\
\rowcolor{LightBlue}      4096   & 9.4e-11 & 3.37e-08 & 1.71e-08 & 1.77e-08 & 4.97 & 3.97 & 3.96 & 6 & 1 & 0.868 & 0.74 & 1.61 & 90.8        \\
\rowcolor{LightBlue}      16384  & 3.01e-12 & 2.14e-09 & 1.09e-09 & 1.13e-09 & 4.96 & 3.98 & 3.96 & 6 & 1 & 4.28 & 3 & 7.27 & 88.4           \\
\rowcolor{LightCyan}      4      & 0.00201 & 0.0243 & 0.0117 & 0.0123 & - & - & -  & 5  & 1 & 0.00107 & 0.000694 & 0.00177 & -                \\  
\rowcolor{LightCyan}      16     & 7.42e-05 & 0.00181 & 0.000951 & 0.000825 & 4.76 & 3.74 & 3.62 & 5  & 1 & 0.00245 & 0.00263 & 0.00507 & 139 \\  
\rowcolor{LightCyan}      64     & 2.83e-06 & 0.00013 & 6.45e-05 & 6.3e-05 & 4.71 & 3.8 & 3.88 & 6  & 1 & 0.00937 & 0.0103 & 0.0196 & 103     \\  
\rowcolor{LightCyan}      256    & 9.43e-08 & 8.43e-06 & 4.07e-06 & 4.25e-06 & 4.91 & 3.95 & 3.98 & 6  & 1 & 0.0404 & 0.041 & 0.0814 & 96.5   \\  
\rowcolor{LightCyan}      1024   & 2.95e-09 & 5.28e-07 & 2.65e-07 & 2.77e-07 & 5 & 4 & 3.94 & 6  & 1 & 0.229 & 0.165 & 0.394 & 82.7           \\  
\rowcolor{LightCyan}      4096   & 9.4e-11 & 3.37e-08 & 1.71e-08 & 1.77e-08 & 4.97 & 3.97 & 3.96 & 6  & 1 & 1.13 & 0.66 & 1.79 & 88.2         \\  
\rowcolor{LightCyan}      16384  & 2.99e-12 & 2.14e-09 & 1.09e-09 & 1.13e-09 & 4.97 & 3.98 & 3.97 & 6  & 1 & 5.74 & 2.68 & 8.42 & 84.9        \\ 
\hline
&	\multicolumn{13}{l}{delaunay triangular grid} \\
\hline
\rowcolor{BallBlue}  8       & 0.000697 & 0.0103 & 0.00823 & 0.00887 & - & - & -  & 12 & 1 & 0.00396 & 0.00198 & 0.00594 & -          \\ 
\rowcolor{BallBlue}  50      & 9.52e-06 & 0.000334 & 0.000281 & 0.0003 & 4.69 & 3.74 & 3.69 & 13 & 1 & 0.0186 & 0.0115 & 0.0301 & 118 \\ 
\rowcolor{BallBlue}  192     & 2.91e-07 & 2.14e-05 & 1.77e-05 & 2e-05 & 5.19 & 4.08 & 4.11 & 14 & 1 & 0.0841 & 0.0393 & 0.123 & 73.1  \\ 
\rowcolor{BallBlue}  810     & 1.03e-08 & 1.41e-06 & 1.12e-06 & 1.29e-06 & 4.64 & 3.78 & 3.83 & 14 & 1 & 0.43 & 0.167 & 0.597 & 82.8  \\ 
\rowcolor{BallBlue}  3120    & 3.43e-10 & 9.1e-08 & 7.19e-08 & 8.3e-08 & 5.04 & 4.06 & 4.08 & 14 & 1 & 1.8 & 0.644 & 2.44 & 73.2      \\ 
\rowcolor{BallBlue}  12780   & 1.01e-11 & 5.41e-09 & 4.32e-09 & 4.97e-09 & 5 & 4 & 3.99 & 14 & 1 & 8.06 & 2.66 & 10.7 & 91.2          \\ 
\rowcolor{BallBlue}  50744   & 5.64e-13 & 3.51e-10 & 2.81e-10 & 3.2e-10 & 4.19 & 3.97 & 3.96 & 14 & 1 & 35.1 & 10.6 & 45.7 & 70.4     \\ 
\rowcolor{LightBlue}      8       & 0.000697 & 0.0103 & 0.00823 & 0.00887 & - & - & -  & 9 & 1 & 0.00146 & 0.00115 & 0.0026 & -              \\
\rowcolor{LightBlue}      50      & 9.52e-06 & 0.000334 & 0.000281 & 0.0003 & 4.69 & 3.74 & 3.69 & 10 & 1 & 0.00568 & 0.00688 & 0.0126 & 124 \\
\rowcolor{LightBlue}      192     & 2.91e-07 & 2.14e-05 & 1.77e-05 & 2e-05 & 5.19 & 4.08 & 4.11 & 10 & 1 & 0.0213 & 0.0261 & 0.0474 & 79.4   \\
\rowcolor{LightBlue}      810     & 1.03e-08 & 1.41e-06 & 1.12e-06 & 1.29e-06 & 4.64 & 3.78 & 3.83 & 12 & 1 & 0.124 & 0.111 & 0.235 & 80.8   \\
\rowcolor{LightBlue}      3120    & 3.43e-10 & 9.1e-08 & 7.19e-08 & 8.3e-08 & 5.04 & 4.06 & 4.08 & 12 & 1 & 0.558 & 0.427 & 0.984 & 71.6     \\
\rowcolor{LightBlue}      12780   & 1.01e-11 & 5.41e-09 & 4.32e-09 & 4.97e-09 & 5 & 4 & 3.99 & 11 & 1 & 2.44 & 1.75 & 4.19 & 93.9            \\
\rowcolor{LightBlue}      50744   & 4.41e-13 & 3.46e-10 & 2.77e-10 & 3.18e-10 & 4.54 & 3.99 & 3.98 & 12 & 1 & 12 & 7.01 & 19 & 66.1          \\
\rowcolor{LightCyan}      8       & 0.000697 & 0.0103 & 0.00823 & 0.00887 & - & - & -  & 7  & 1 & 0.0018 & 0.00127 & 0.00307 & -             \\ 
\rowcolor{LightCyan}      50      & 9.52e-06 & 0.000334 & 0.000281 & 0.0003 & 4.69 & 3.74 & 3.69 & 9  & 1 & 0.00716 & 0.00709 & 0.0143 & 129 \\ 
\rowcolor{LightCyan}      192     & 2.91e-07 & 2.14e-05 & 1.77e-05 & 2e-05 & 5.19 & 4.08 & 4.11 & 9  & 1 & 0.0271 & 0.0243 & 0.0514 & 83.1   \\ 
\rowcolor{LightCyan}      810     & 1.03e-08 & 1.41e-06 & 1.12e-06 & 1.29e-06 & 4.64 & 3.78 & 3.83 & 9  & 1 & 0.146 & 0.102 & 0.248 & 83     \\ 
\rowcolor{LightCyan}      3120    & 3.43e-10 & 9.1e-08 & 7.19e-08 & 8.3e-08 & 5.04 & 4.06 & 4.08 & 9  & 1 & 0.643 & 0.396 & 1.04 & 71.6      \\ 
\rowcolor{LightCyan}      12780   & 1.01e-11 & 5.41e-09 & 4.32e-09 & 4.97e-09 & 5 & 4 & 3.99 & 9  & 1 & 3.03 & 1.64 & 4.67 & 89              \\ 
\rowcolor{LightCyan}      50744   & 3.51e-13 & 3.45e-10 & 2.76e-10 & 3.18e-10 & 4.88 & 3.99 & 3.99 & 9  & 1 & 14.5 & 6.48 & 21 & 66.7        \\ 
\hline
\end{tabular}
\caption{Evaluation of $p$-multilevel solution strategies for solving high-order $k{=}3$ \hhodp{} over 2D regular mesh sequences. 
         Solvers are applied to uncondensed and statically condensed matrices (identified by different colors) considering 
         two alternative Schur complement implementations, see text for details. See Section \ref{sec:numResOptions} for solver options. \label{tab:numResStatCondReg2D}}
\end{sidewaystable}

\begin{sidewaystable}
\small
\centering
\begin{tabular}{ c | c c  c c | c c c | c c | c c c | c }
\multicolumn{4}{c}{{\cellcolor{BallBlue} \hhodp{} \uncond{}}} 
                                              & \multicolumn{5}{c}{{\cellcolor{LightBlue} \hhodp{} \vpcond{} }} 
                                              & \multicolumn{5}{c}{{\cellcolor{LightCyan} \hhodp{} \vcond{} }}\\
$\card{\Th}$ & \multicolumn{4}{c|}{error in $L^2$ norm} & \multicolumn{3}{c|}{conv. rate} & \multicolumn{2}{c|}{ITs} & \multicolumn{3}{c|}{CPU time} & Eff \\ \hline
 &  $\uVec_h$    & $G\uVec_h$   & $p_h$ & $D\uVec_h $& $\uVec_h$   & $G\uVec_h$ & $p_h$  & ITs & ITs$_L$ & Sol. & Ass. & Tot. & \\ \hline
 & \multicolumn{13}{l}{graded quadrilateral elements grid}\\
\hline
\rowcolor{BallBlue} 4       & 0.00201 & 0.0243 & 0.0117 & 0.0123 & - & - & -  & 10 & 1 & 0.00202 & 0.00107 & 0.00309 & -              \\
\rowcolor{BallBlue} 16      & 0.000296 & 0.00552 & 0.00292 & 0.00202 & 2.76 & 2.14 & 2 & 12 & 1 & 0.00695 & 0.00434 & 0.0113 & 110    \\
\rowcolor{BallBlue} 64      & 1.23e-05 & 0.000398 & 0.000225 & 0.000208 & 4.58 & 3.8 & 3.7 & 13 & 1 & 0.0288 & 0.0166 & 0.0454 & 99.4 \\
\rowcolor{BallBlue} 256     & 3.52e-07 & 2.31e-05 & 1.16e-05 & 1.18e-05 & 5.13 & 4.11 & 4.28 & 15 & 1 & 0.153 & 0.0663 & 0.219 & 82.9 \\
\rowcolor{BallBlue} 1024    & 1.12e-08 & 1.46e-06 & 7.49e-07 & 7.59e-07 & 4.97 & 3.98 & 3.95 & 17 & 1 & 0.82 & 0.266 & 1.09 & 80.8    \\
\rowcolor{BallBlue} 4096    & 3.57e-10 & 9.3e-08 & 4.86e-08 & 4.93e-08 & 4.97 & 3.97 & 3.95 & 23 & 1 & 4.54 & 1.06 & 5.61 & 77.4      \\
\rowcolor{BallBlue} 16384   & 1.71e-11 & 7e-09 & 4.11e-09 & 3.28e-09 & 4.39 & 3.73 & 3.56 & 40 & 1 & 30.7 & 4.37 & 35 & 64            \\
\rowcolor{LightBlue} 4       & 0.00201 & 0.0243 & 0.0117 & 0.0123 & - & - & -  & 4 & 1 & 0.000941 & 0.000753 & 0.00169 & -             \\
\rowcolor{LightBlue} 16      & 0.000296 & 0.00552 & 0.00292 & 0.00202 & 2.76 & 2.14 & 2 & 5 & 1 & 0.004 & 0.00612 & 0.0101 & 67        \\
\rowcolor{LightBlue} 64      & 1.23e-05 & 0.000398 & 0.000225 & 0.000208 & 4.58 & 3.8 & 3.7 & 6 & 1 & 0.00801 & 0.0133 & 0.0213 & 190  \\
\rowcolor{LightBlue} 256     & 3.52e-07 & 2.31e-05 & 1.16e-05 & 1.18e-05 & 5.13 & 4.11 & 4.28 & 7 & 1 & 0.0376 & 0.0459 & 0.0835 & 102 \\
\rowcolor{LightBlue} 1024    & 1.12e-08 & 1.46e-06 & 7.49e-07 & 7.59e-07 & 4.97 & 3.98 & 3.95 & 8 & 1 & 0.213 & 0.184 & 0.397 & 84.1   \\
\rowcolor{LightBlue} 4096    & 3.57e-10 & 9.3e-08 & 4.86e-08 & 4.93e-08 & 4.97 & 3.97 & 3.95 & 11 & 1 & 1.23 & 0.736 & 1.97 & 80.7     \\
\rowcolor{LightBlue} 16384   & 1.18e-11 & 5.97e-09 & 4.48e-09 & 3.14e-09 & 4.92 & 3.96 & 3.44 & 24 & 1 & 9.53 & 2.99 & 12.5 & 62.9     \\
\rowcolor{LightCyan} 4      & 0.00201 & 0.0243 & 0.0117 & 0.0123 & - & - & -  & 5 & 1 & 0.00102 & 0.000814 & 0.00184 & -               \\
\rowcolor{LightCyan} 16     & 0.000296 & 0.00552 & 0.00292 & 0.00202 & 2.76 & 2.14 & 2 & 5 & 1 & 0.00244 & 0.00265 & 0.00508 & 145     \\
\rowcolor{LightCyan} 64     & 1.23e-05 & 0.000398 & 0.000225 & 0.000208 & 4.58 & 3.8 & 3.7 & 6 & 1 & 0.00923 & 0.0103 & 0.0196 & 104   \\
\rowcolor{LightCyan} 256    & 3.52e-07 & 2.31e-05 & 1.16e-05 & 1.18e-05 & 5.13 & 4.11 & 4.28 & 6 & 1 & 0.0406 & 0.0411 & 0.0817 & 95.8 \\
\rowcolor{LightCyan} 1024   & 1.12e-08 & 1.46e-06 & 7.49e-07 & 7.59e-07 & 4.97 & 3.98 & 3.95 & 6 & 1 & 0.232 & 0.165 & 0.396 & 82.5    \\
\rowcolor{LightCyan} 4096   & 3.57e-10 & 9.3e-08 & 4.85e-08 & 4.93e-08 & 4.97 & 3.97 & 3.95 & 6 & 1 & 1.13 & 0.658 & 1.79 & 88.5       \\
\rowcolor{LightCyan} 16384  & 1.16e-11 & 5.95e-09 & 3.2e-09 & 3.2e-09 & 4.94 & 3.97 & 3.92 & 7 & 1 & 6.07 & 2.69 & 8.75 & 81.9         \\
\hline
& \multicolumn{13}{l}{graded triangular elements grid} \\
\hline
\rowcolor{BallBlue}  8      &0.000743 & 0.00999 & 0.00762 & 0.0083 & - & - & -  & 12 & 1 & 0.00313 & 0.00185 & 0.00498 & -           \\
\rowcolor{BallBlue}  32     &0.000182 & 0.00321 & 0.00238 & 0.00286 & 2.03 & 1.64 & 1.68 & 16 & 1 & 0.013 & 0.00652 & 0.0195 & 102   \\
\rowcolor{BallBlue}  128    &7.9e-06 & 0.000253 & 0.000202 & 0.000228 & 4.53 & 3.67 & 3.56 & 20 & 1 & 0.065 & 0.0259 & 0.0908 & 85.8 \\
\rowcolor{BallBlue}  512    &2.36e-07 & 1.54e-05 & 1.13e-05 & 1.32e-05 & 5.06 & 4.04 & 4.17 & 27 & 1 & 0.419 & 0.103 & 0.522 & 69.6  \\
\rowcolor{BallBlue}  2048   &6.76e-09 & 9.11e-07 & 6.89e-07 & 7.85e-07 & 5.13 & 4.08 & 4.03 & 46 & 1 & 2.99 & 0.414 & 3.41 & 61.3    \\
\rowcolor{BallBlue}  8192   &2.21e-10 & 5.85e-08 & 4.39e-08 & 5e-08 & 4.94 & 3.96 & 3.97 & 75 & 1 & 20 & 1.66 & 21.6 & 63            \\
\rowcolor{BallBlue}  32768  &1.17e-11 & 4.9e-09 & 4.03e-09 & 3.3e-09 & 4.24 & 3.58 & 3.45 & 160 & 1 & 169 & 6.75 & 176 & 49.3        \\
\rowcolor{LightBlue}  8      &0.000743 & 0.00999 & 0.00762 & 0.0083 & - & - & -  & 7 & 1 & 0.00122 & 0.00125 & 0.00247 & -             \\ 
\rowcolor{LightBlue}  32     &0.000182 & 0.00321 & 0.00238 & 0.00286 & 2.03 & 1.64 & 1.68 & 10 & 1 & 0.00357 & 0.00436 & 0.00793 & 124 \\ 
\rowcolor{LightBlue}  128    &7.9e-06 & 0.000253 & 0.000202 & 0.000228 & 4.53 & 3.67 & 3.56 & 11 & 1 & 0.0133 & 0.0173 & 0.0306 & 104  \\ 
\rowcolor{LightBlue}  512    &2.36e-07 & 1.54e-05 & 1.13e-05 & 1.32e-05 & 5.06 & 4.04 & 4.17 & 16 & 1 & 0.0772 & 0.0689 & 0.146 & 83.7 \\ 
\rowcolor{LightBlue}  2048   &6.77e-09 & 9.11e-07 & 6.91e-07 & 7.85e-07 & 5.13 & 4.08 & 4.03 & 649 & 1 & 12.5 & 0.275 & 12.7 & 4.59    \\ 
\rowcolor{LightBlue}  8192   &0.0313 & 9.51 & 19.1 & 1.22 & -22.1 & -23.3 & -24.7 & 1000$^*$ & 1 & 82.8 & 1.1 & 83.9 & 60.7            \\ 
\rowcolor{LightBlue}  32768  &0.566 & 261 & 1.31e+03 & 72.3 & -4.18 & -4.78 & -6.1 & 1000$^*$ & 1 & 345 & 4.44 & 349 & 96.2            \\ 
\rowcolor{LightCyan}  8      & 0.000743 & 0.00999 & 0.00762 & 0.0083 & - & - & -  & 8 & 1 & 0.00154 & 0.00125 & 0.00279 & -            \\
\rowcolor{LightCyan}  32     & 0.000182 & 0.00321 & 0.00238 & 0.00286 & 2.03 & 1.64 & 1.68 & 8 & 1 & 0.00422 & 0.00415 & 0.00837 & 134 \\
\rowcolor{LightCyan}  128    & 7.9e-06 & 0.000253 & 0.000202 & 0.000228 & 4.53 & 3.67 & 3.56 & 8 & 1 & 0.0153 & 0.0163 & 0.0316 & 106  \\
\rowcolor{LightCyan}  512    & 2.36e-07 & 1.54e-05 & 1.13e-05 & 1.32e-05 & 5.06 & 4.04 & 4.17 & 8 & 1 & 0.0682 & 0.065 & 0.133 & 95    \\
\rowcolor{LightCyan}  2048   & 6.76e-09 & 9.11e-07 & 6.89e-07 & 7.85e-07 & 5.13 & 4.08 & 4.03 & 9 & 1 & 0.38 & 0.255 & 0.636 & 83.8    \\
\rowcolor{LightCyan}  8192   & 2.21e-10 & 5.85e-08 & 4.39e-08 & 5e-08 & 4.94 & 3.96 & 3.97 & 11 & 1 & 2.02 & 1.03 & 3.05 & 83.2        \\
\rowcolor{LightCyan}  32768  & 7.45e-12 & 3.81e-09 & 2.92e-09 & 3.27e-09 & 4.89 & 3.94 & 3.91 & 18 & 1 & 12.8 & 4.15 & 17 & 72.1       \\
\hline
\end{tabular}
\caption{Evaluation of $p$-multilevel solution strategies for solving high-order $k{=}3$ \hhodp{} over 2D graded mesh sequences. 
         Solvers are applied to uncondensed and statically condensed matrices (identified by different colors) considering 
         two alternative Schur complement implementations, see text for details. See Section \ref{sec:numResOptions} for solver options. \label{tab:numResStatCondGrad2D}}
\end{sidewaystable}

\begin{sidewaystable}
\small
\centering
\begin{tabular}{ c | c c  c c | c c c | c c | c c c | c }
\multicolumn{4}{c}{{\cellcolor{BallBlue} \hhodp{} \uncond{} }} 
                                              & \multicolumn{5}{c}{{\cellcolor{LightBlue} \hhodp{} \vpcond{} }} 
                                              & \multicolumn{5}{c}{{\cellcolor{LightCyan} \hhodp{} \vcond{} }}\\
$\card{\Th}$ & \multicolumn{4}{c|}{error in $L^2$ norm} & \multicolumn{3}{c|}{conv. rate} & \multicolumn{2}{c|}{ITs} & \multicolumn{3}{c|}{CPU time} & Eff \\ \hline
 &  $\uVec_h$    & $G\uVec_h$   & $p_h$ & $D\uVec_h $& $\uVec_h$   & $G\uVec_h$ & $p_h$  & ITs & ITs$_L$ & Sol. & Ass. & Tot. & \\ \hline
& \multicolumn{12}{l}{prismatic elements grid} \\
\hline
\rowcolor{BallBlue} 16   &   0.000721 & 0.0162 & 0.0128 & 0.00713 & - & - & -  & 13 & 21 & 0.0689 & 0.0503 & 0.119 & -  \\
\rowcolor{BallBlue} 128  &   2.41e-05 & 0.00105 & 0.000732 & 0.00047 & 4.91 & 3.95 & 4.13 & 14 & 62 & 0.884 & 0.416 & 1.3 & 73.4 \\
\rowcolor{BallBlue} 1024 &   7.77e-07 & 6.57e-05 & 4.21e-05 & 2.93e-05 & 4.95 & 3.99 & 4.12 & 14 & 146 & 14.2 & 3.42 & 17.6 & 59.2 \\
\rowcolor{BallBlue} 8192 &   2.47e-08 & 4.11e-06 & 2.5e-06 & 1.82e-06 & 4.97 & 4 & 4.08 & 15 & 333 & 392 & 27.9 & 420 & 33.5 \\
\rowcolor{LightBlue} 16   & 0.000721 & 0.0162 & 0.0128 & 0.00713 & - & - & -  & 9 & 12 & 0.0516 & 0.0418 & 0.0934 & -             \\
\rowcolor{LightBlue} 128  & 2.41e-05 & 0.00105 & 0.000732 & 0.00047 & 4.91 & 3.95 & 4.13 & 10 & 27 & 0.597 & 0.357 & 0.954 & 78.3 \\
\rowcolor{LightBlue} 1024 & 7.77e-07 & 6.57e-05 & 4.21e-05 & 2.93e-05 & 4.95 & 3.99 & 4.12 & 11 & 55 & 6.86 & 3 & 9.86 & 77.4     \\
\rowcolor{LightBlue} 8192 & 2.47e-08 & 4.11e-06 & 2.5e-06 & 1.82e-06 & 4.97 & 4 & 4.08 & 10 & 150 & 86.7 & 24.5 & 111 & 70.9      \\
\rowcolor{LightCyan} 16   & 0.000721 & 0.0162 & 0.0128 & 0.00713 & - & - & -  & 8 & 8 & 0.0253 & 0.0342 & 0.0595 & -             \\ 
\rowcolor{LightCyan} 128  & 2.41e-05 & 0.00105 & 0.000732 & 0.00047 & 4.91 & 3.95 & 4.13 & 9 & 18 & 0.291 & 0.284 & 0.575 & 82.7 \\ 
\rowcolor{LightCyan} 1024 & 7.77e-07 & 6.57e-05 & 4.21e-05 & 2.93e-05 & 4.95 & 3.99 & 4.12 & 9 & 42 & 2.95 & 2.35 & 5.3 & 86.8   \\ 
\rowcolor{LightCyan} 8192 & 2.47e-08 & 4.11e-06 & 2.5e-06 & 1.82e-06 & 4.97 & 4 & 4.08 & 10 & 105 & 42.2 & 19.1 & 61.3 & 69.1    \\ 
\hline
& \multicolumn{13}{l}{pyramidal elements grid}\\
\hline	
\rowcolor{BallBlue} 48    &  0.000278 & 0.00811 & 0.00306 & 0.0029 & - & - & -  & 16 & 36 & 0.3 & 0.164 & 0.464 & -  \\
\rowcolor{BallBlue} 384   &  8.93e-06 & 0.000519 & 0.000191 & 0.000193 & 4.96 & 3.97 & 4 & 19 & 86 & 4.31 & 1.33 & 5.64 & 65.7 \\
\rowcolor{BallBlue} 3072  &  2.82e-07 & 3.26e-05 & 1.2e-05 & 1.23e-05 & 4.99 & 3.99 & 4 & 19 & 186 & 82.3 & 10.6 & 92.9 & 48.6 \\
\rowcolor{BallBlue} 24576 &  8.82e-09 & 2.04e-06 & 7.48e-07 & 7.79e-07 & 5 & 4 & 4 & 19 & 467 & 2.46e+03 & 84.9 & 2.54e+03 & 29.3 \\
\rowcolor{LightBlue} 48    & 0.000278 & 0.00811 & 0.00306 & 0.0029 & - & - & -  & 18 & 25 & 0.31 & 0.141 & 0.451 & -           \\
\rowcolor{LightBlue} 384   & 8.93e-06 & 0.000519 & 0.000191 & 0.000193 & 4.96 & 3.97 & 4 & 19 & 47 & 3.39 & 1.15 & 4.55 & 79.4 \\
\rowcolor{LightBlue} 3072  & 2.82e-07 & 3.26e-05 & 1.2e-05 & 1.23e-05 & 4.99 & 3.99 & 4 & 20 & 95 & 42.6 & 9.43 & 52 & 69.9    \\
\rowcolor{LightBlue} 24576 & 8.82e-09 & 2.04e-06 & 7.48e-07 & 7.79e-07 & 5 & 4 & 4 & 20 & 168 & 694 & 76 & 770 & 54            \\
\rowcolor{LightCyan} 48    & 0.000278 & 0.00811 & 0.00306 & 0.0029 & - & - & -  & 10 & 13 & 0.1 & 0.111 & 0.211 & -             \\  
\rowcolor{LightCyan} 384   & 8.93e-06 & 0.000519 & 0.000191 & 0.000193 & 4.96 & 3.97 & 4 & 10 & 31 & 1.05 & 0.898 & 1.95 & 86.6 \\  
\rowcolor{LightCyan} 3072  & 2.82e-07 & 3.26e-05 & 1.2e-05 & 1.23e-05 & 4.99 & 3.99 & 4 & 10 & 64 & 11.6 & 7.28 & 18.9 & 82.6   \\  
\rowcolor{LightCyan} 24576 & 8.82e-09 & 2.04e-06 & 7.48e-07 & 7.79e-07 & 5 & 4 & 4 & 10 & 146 & 169 & 58.5 & 228 & 66.4         \\  
\hline
& \multicolumn{12}{l}{graded tetrahedral elements grid} \\
\hline
\rowcolor{BallBlue} 24    &  0.00156 & 0.0342 & 0.0153 & 0.0166 & - & - & -  & 15 & 28 & 0.098 & 0.0664 & 0.164 & -  \\
\rowcolor{BallBlue} 192   &  4.82e-05 & 0.00213 & 0.000893 & 0.00105 & 5.01 & 4.01 & 4.1 & 18 & 73 & 1.39 & 0.539 & 1.93 & 68.2 \\
\rowcolor{BallBlue} 1536  &  7.23e-06 & 0.000458 & 0.00016 & 0.000228 & 2.74 & 2.22 & 2.48 & 25 & 154 & 30.5 & 4.34 & 34.8 & 44.3 \\
\rowcolor{BallBlue} 12288 &  2.79e-07 & 3.31e-05 & 1.17e-05 & 1.64e-05 & 4.69 & 3.79 & 3.77 & 37 & 343 & 1.11e+03 & 34.8 & 1.15e+03 & 24.3 \\
\rowcolor{LightBlue} 24    & 0.00156 & 0.0342 & 0.0153 & 0.0166 & - & - & -  & 13 & 15 & 0.0597 & 0.0495 & 0.109 & -                        \\
\rowcolor{LightBlue} 192   & 4.82e-05 & 0.00213 & 0.000893 & 0.00105 & 5.01 & 4.01 & 4.1 & 17 & 39 & 0.863 & 0.418 & 1.28 & 68.1            \\
\rowcolor{LightBlue} 1536  & 7.23e-06 & 0.000458 & 0.00016 & 0.000228 & 2.74 & 2.22 & 2.48 & 47 & 121 & 36.3 & 3.41 & 39.7 & 25.8           \\
\rowcolor{LightBlue} 12288 & 2.79e-07 & 3.31e-05 & 1.17e-05 & 1.64e-05 & 4.69 & 3.79 & 3.77 & 256 & 560 & 5.98e+03 & 27.7 & 6.01e+03 & 5.29 \\
\rowcolor{LightCyan} 24    & 0.00156 & 0.0342 & 0.0153 & 0.0166 & - & - & -  & 11 & 12 & 0.0335 & 0.0436 & 0.077 & -              \\
\rowcolor{LightCyan} 192   & 4.82e-05 & 0.00213 & 0.000893 & 0.00105 & 5.01 & 4.01 & 4.1 & 11 & 20 & 0.371 & 0.362 & 0.734 & 84   \\
\rowcolor{LightCyan} 1536  & 7.23e-06 & 0.000458 & 0.00016 & 0.000228 & 2.74 & 2.22 & 2.48 & 13 & 45 & 4.45 & 2.9 & 7.34 & 79.9   \\
\rowcolor{LightCyan} 12288 & 2.79e-07 & 3.31e-05 & 1.17e-05 & 1.64e-05 & 4.69 & 3.79 & 3.77 & 18 & 108 & 86.6 & 23.5 & 110 & 53.4 \\
\hline
\end{tabular}
\caption{Evaluation of $p$-multilevel solution strategies for solving high-order $k{=}3$ \hhodp{} over 3D mesh sequences. 
         Solvers are applied to uncondensed and statically condensed matrices (identified by different colors) considering 
         two alternative Schur complement implementations, see text for details. See Section \ref{sec:numResOptions} for solver options. \label{tab:numResStatCond3D}}
\end{sidewaystable}

\subsection{Comparison based on accuracy and efficiency of the solver strategy}
In this section we compare the three nonconforming discretizations of the Stokes problem presented in Section \ref{sec:threeNCD} based on accuracy and performance of the multilevel solver strategy.
For the HHO scheme \hhodp{}, in accordance with the results of Section \ref{sec:nrStatCond}, the static condensation strategy \vcond{} is used for all meshes in both two and three space dimensions.
For the HHO scheme \hhohp{}, we consider static condensation of the element unknowns for both the velocity and the pressure (\hhohp{} \vpcond), so that only skeletal unknowns appear in the global systems.
The results for 2D regular and graded sequences are reported in Tables \ref{tab:numResPerfReg2DhighO}--\ref{tab:numResPerfReg2DhigherO} 
and \ref{tab:numResPerfGrad2DhighO}--\ref{tab:numResPerfGrad2DhigherO}, respectively. 
The results for 3D mesh sequences are reported in Tables \ref{tab:numResPerf3DhighO}--\ref{tab:numResPerf3DhigherO}.

As a first point, we remark that the theoretical convergence rates are confirmed for all the test cases performed on regular 2D and 3D mesh sequences.
When higher-order ($k{=}6$) discretizations are considered and machine precision is reached, the converge rates deteriorates, as expected.
Turning to graded mesh sequences, we observe a slightly suboptimal convergence of \hhohp{} with respect to \hhodp{} over graded triangular meshes at higher-order. 
Note that velocity gradients and pressure fields reach an asymptotic sixth order convergence rate 
and the divergence error is a bit higher than expected, compare for example with the result on standard meshes.
Interestingly, all the nonconforming discretizations suffer from a convergence degradation for $\card{\Th}$ between $192$ and $1546$ over the graded tetrahedral mesh sequence. 
This is probably due to mesh elements of extremely bad quality generated as a result of grading plus random node displacement, see Section \ref{sec:meshSeq}.
Overall, both \hhodp{} and \hhohp{} outperform \dg{} in terms of accuracy with order of magnitudes gains observed moving towards finer meshes. 
This is due to better asymptotic convergence rates (one order higher) as well as better accuracy on coarse meshes. 

$p$-Multilevel solvers guarantee uniform convergence with respect to the mesh density when standard 2D and 3D mesh sequences are considered:
note that the number of FGMRES iterations is almost uniform all along the mesh sequence. 
Interestingly, \hhodp{} discretizations show uniform convergence with respect to the mesh density on graded quadrilateral meshes, 
while \dg{} is the most affected by mesh grading, especially for $k{=}6$.
For \hhohp{}, the number of iterations increases with mesh density on graded quadrilateral meshes. Nevertheless,
the number of iterations over coarse meshes is remarkably small and grows up to match the iterations count of \hhodp{} over fine meshes.
The solver convergence deteriorates with the mesh density in case of graded triangular and tetrahedral mesh sequences:
the iterations increase is clearly visible but not pathological in case \hho{} discretizations. 

Interestingly, $p$-multilevel solvers deliver almost uniform convergence with respect to the polynomial degree when applied to \hho{} discretizations:
moving from high-order ($k{=}3) $ to higher-order ($k{=}6$) entails a mild iterations increase for \hho{}, while the iteration count doubles for \dg{}.  
In 2D this behaviour has a strong impact on computation times: \hho{} is up to three and eight times faster than \dg{} at high-order and higher-order, respectively.
\hhodp{} outperforms \dg{} because of the reduced number of matrix non-zero entries and the reduced matrix dimension, see Tables \ref{tab:dofs2D} and \ref{tab:mnzs2D}:
the former influences the cost of smoothing iterations while the latter strongly influences the cost of the LU factorization on the coarse level.

Let us consider the performance of the multilevel solver in 3D.
\hhodp{} is two times and four-to-five times faster than \dg{} in terms of solution times for $k{=}3$ and $k{=}6$, respectively.
\hhohp{} is slower than \hhodp{} in terms of solution times and faster than \dg{} by a small amount, with the exception of the pyramidal elements mesh sequence for $k{=}3$.
The difference in computational cost between \hhodp{} and \hhohp{} is essentially due to the number of MNZs, see Table \ref{tab:mnzs3D}, while the
number of FGMRES iterations is comparable.
Since in 3D the coarse level solver is generally more efficient for \dg, the \hho{} advantage results from the efficiency of the smoothers and the
reduced number of FGMRES iterations. In particular, we remark that \dg{} has fewer DOFs than \hho{} for $k{=}1$, see Table \ref{tab:dofs3D}.
Moreover, \dg{} and \hhodp{} \vcond{} have a comparable MNZs count for $k{=}1$, significantly smaller than the MNZs count of \hhohp{} \vpcond, see Table \ref{tab:mnzs3D}.

Overall, the gain in terms of total execution times is less significant than in 2D.
When working with \hho{} in three space dimensions, assembly times are a considerable fraction of the total computation time: 
matrix assembly is twice as expensive as linear system solution for \hhodp{} for $k{=}6$. 
As opposite, for \dg{}, solution times dominate.
Increased assembly costs are essentially due to the increased expense of solving local problems involved in static condensation. 
An important observation is that, since the assembly procedure is perfectly scalable while ILU preconditioned smoothers are not, 
\hho{} discretizations might show better scalability results as compared to \dg{} in massively parallel computations. 

We conclude this section commenting about solver efficiency (last column in Tables \ref{tab:numResPerfReg2DhighO}--\ref{tab:numResPerf3DhigherO}). 
It is clear that higher-order discretizations ($k{=}6$) achieve better efficiency than high-order discretizations ($k{=}3$), in both 2D and 3D.
This outlines the intrinsic limitation of $p$-multilevel solution strategies: when considering fine meshes, the performance of the coarse
solver might limit the efficiency because the number of DOFs and MNZs on the coarse level can not be chosen arbitrarily low. 
Accordingly, $p$-multilevel solver are best suited for those situations where arbitrarily coarse meshes with higher-order polynomials can be employed. 
\begin{sidewaystable}
\small
\centering
\begin{tabular}{ c | c c  c c | c c c | c c | c c c | c }
\multicolumn{4}{c}{{\cellcolor{LightGoldenrod} \dg }} 
                                              & \multicolumn{5}{c}{{\cellcolor{LightCyan} \hhodp{} \vcond{} }}
                                              & \multicolumn{5}{c}{{\cellcolor{DardSeeGreen} \hhohp{} \vpcond{} }}\\
$\card{\Th}$ & \multicolumn{4}{c|}{error in $L^2$ norm} & \multicolumn{3}{c|}{conv. rate} & \multicolumn{2}{c|}{ITs} & \multicolumn{3}{c|}{CPU time} & Eff \\ \hline
 &  $\uVec_h$    & $G\uVec_h$   & $p_h$ & $D\uVec_h $& $\uVec_h$   & $G\uVec_h$ & $p_h$  & ITs & ITs$_L$ & Sol. & Ass. & Tot. & \\ \hline
&	\multicolumn{13}{l}{trapezoidal elements grid}\\
\hline
\rowcolor{LightGoldenrod} 4      & 0.0128 & 0.154 & 0.0456 & 0.0411 & - & - & -  & 6 & 1 & 0.0016 & 0.000337 & 0.00194 & -                 \\
\rowcolor{LightGoldenrod} 16     & 0.000943 & 0.0221 & 0.00442 & 0.00635 & 3.76 & 2.8 & 3.37 & 8 & 1 & 0.00769 & 0.00128 & 0.00896 & 86.6  \\
\rowcolor{LightGoldenrod} 64     & 6.73e-05 & 0.00305 & 0.000628 & 0.000883 & 3.81 & 2.86 & 2.82 & 8 & 1 & 0.0329 & 0.00513 & 0.038 & 94.3 \\
\rowcolor{LightGoldenrod} 256    & 4.28e-06 & 0.000388 & 7.49e-05 & 0.000112 & 3.97 & 2.98 & 3.07 & 10 & 1 & 0.193 & 0.0208 & 0.214 & 71.2 \\
\rowcolor{LightGoldenrod} 1024   & 2.69e-07 & 4.87e-05 & 9.31e-06 & 1.43e-05 & 3.99 & 2.99 & 3.01 & 11 & 1 & 0.98 & 0.0835 & 1.06 & 80.4   \\
\rowcolor{LightGoldenrod} 4096   & 1.72e-08 & 6.17e-06 & 1.19e-06 & 1.81e-06 & 3.97 & 2.98 & 2.97 & 12 & 1 & 5.02 & 0.337 & 5.35 & 79.5    \\
\rowcolor{LightGoldenrod} 16384  & 1.09e-09 & 7.79e-07 & 1.52e-07 & 2.31e-07 & 3.98 & 2.98 & 2.97 & 12 & 1 & 26 & 1.35 & 27.4 & 78.2       \\
\rowcolor{LightCyan}      4      & 0.00201 & 0.0243 & 0.0117 & 0.0123 & - & - & -  & 5  & 1 & 0.00107 & 0.000694 & 0.00177 & -                \\  
\rowcolor{LightCyan}      16     & 7.42e-05 & 0.00181 & 0.000951 & 0.000825 & 4.76 & 3.74 & 3.62 & 5  & 1 & 0.00245 & 0.00263 & 0.00507 & 139 \\  
\rowcolor{LightCyan}      64     & 2.83e-06 & 0.00013 & 6.45e-05 & 6.3e-05 & 4.71 & 3.8 & 3.88 & 6  & 1 & 0.00937 & 0.0103 & 0.0196 & 103     \\  
\rowcolor{LightCyan}      256    & 9.43e-08 & 8.43e-06 & 4.07e-06 & 4.25e-06 & 4.91 & 3.95 & 3.98 & 6  & 1 & 0.0404 & 0.041 & 0.0814 & 96.5   \\  
\rowcolor{LightCyan}      1024   & 2.95e-09 & 5.28e-07 & 2.65e-07 & 2.77e-07 & 5 & 4 & 3.94 & 6  & 1 & 0.229 & 0.165 & 0.394 & 82.7           \\  
\rowcolor{LightCyan}      4096   & 9.4e-11 & 3.37e-08 & 1.71e-08 & 1.77e-08 & 4.97 & 3.97 & 3.96 & 6  & 1 & 1.13 & 0.66 & 1.79 & 88.2         \\  
\rowcolor{LightCyan}      16384  & 2.99e-12 & 2.14e-09 & 1.09e-09 & 1.13e-09 & 4.97 & 3.98 & 3.97 & 6  & 1 & 5.74 & 2.68 & 8.42 & 84.9        \\ 
\rowcolor{DardSeeGreen}      4      & 0.0023 & 0.0315 & 0.00957 & 1.75e-14 & - & - & -  & 3 & 1 & 0.000842 & 0.000708 & 0.00155 & -            \\
\rowcolor{DardSeeGreen}      16     & 9.32e-05 & 0.00238 & 0.000656 & 4.6e-14 & 4.63 & 3.73 & 3.87 & 4 & 1 & 0.00234 & 0.00249 & 0.00483 & 128 \\
\rowcolor{DardSeeGreen}      64     & 3.42e-06 & 0.000171 & 4.35e-05 & 1.22e-13 & 4.77 & 3.8 & 3.91 & 4 & 1 & 0.00856 & 0.00946 & 0.018 & 107  \\
\rowcolor{DardSeeGreen}      256    & 1.1e-07 & 1.08e-05 & 2.64e-06 & 4e-13 & 4.96 & 3.98 & 4.04 & 5 & 1 & 0.0494 & 0.0369 & 0.0862 & 83.6     \\
\rowcolor{DardSeeGreen}      1024   & 3.43e-09 & 6.88e-07 & 1.68e-07 & 1.36e-12 & 5 & 3.98 & 3.97 & 5 & 1 & 0.302 & 0.147 & 0.449 & 76.8       \\
\rowcolor{DardSeeGreen}      4096   & 1.11e-10 & 4.38e-08 & 1.08e-08 & 5.39e-12 & 4.95 & 3.97 & 3.97 & 5 & 1 & 1.72 & 0.584 & 2.3 & 78.2       \\
\rowcolor{DardSeeGreen}      16384  & 3.54e-12 & 2.77e-09 & 6.92e-10 & 2.14e-11 & 4.97 & 3.98 & 3.96 & 5 & 1 & 10.2 & 2.37 & 12.6 & 72.9       \\
\hline
&	\multicolumn{13}{l}{delaunay triangular grid} \\
\hline
\rowcolor{LightGoldenrod} 8       & 0.00498 & 0.0837 & 0.0559 & 0.0356 & - & - & -  & 9 & 1 & 0.00382 & 0.000568 & 0.00439 & -             \\
\rowcolor{LightGoldenrod} 50      & 0.000158 & 0.00613 & 0.00297 & 0.00287 & 3.77 & 2.85 & 3.2 & 12 & 1 & 0.0303 & 0.00337 & 0.0337 & 81.4 \\
\rowcolor{LightGoldenrod} 192     & 9.58e-06 & 0.000761 & 0.000287 & 0.000364 & 4.17 & 3.1 & 3.47 & 13 & 1 & 0.14 & 0.0127 & 0.152 & 84.9  \\
\rowcolor{LightGoldenrod} 810     & 6.73e-07 & 9.8e-05 & 3.65e-05 & 4.64e-05 & 3.69 & 2.85 & 2.87 & 13 & 1 & 0.68 & 0.0538 & 0.734 & 87.5  \\
\rowcolor{LightGoldenrod} 3120    & 4.46e-08 & 1.25e-05 & 4.83e-06 & 5.89e-06 & 4.03 & 3.05 & 3 & 13 & 1 & 2.89 & 0.211 & 3.1 & 91.2       \\
\rowcolor{LightGoldenrod} 12780   & 2.64e-09 & 1.5e-06 & 5.71e-07 & 6.98e-07 & 4.01 & 3.01 & 3.03 & 13 & 1 & 13.7 & 0.873 & 14.6 & 86.9    \\
\rowcolor{LightGoldenrod} 50744   & 1.68e-10 & 1.9e-07 & 7.16e-08 & 8.84e-08 & 4 & 2.99 & 3.01 & 13 & 1 & 68.1 & 3.66 & 71.8 & 80.8        \\
\rowcolor{LightCyan}      8       & 0.000697 & 0.0103 & 0.00823 & 0.00887 & - & - & -  & 7  & 1 & 0.0018 & 0.00127 & 0.00307 & -             \\ 
\rowcolor{LightCyan}      50      & 9.52e-06 & 0.000334 & 0.000281 & 0.0003 & 4.69 & 3.74 & 3.69 & 9  & 1 & 0.00716 & 0.00709 & 0.0143 & 129 \\ 
\rowcolor{LightCyan}      192     & 2.91e-07 & 2.14e-05 & 1.77e-05 & 2e-05 & 5.19 & 4.08 & 4.11 & 9  & 1 & 0.0271 & 0.0243 & 0.0514 & 83.1   \\ 
\rowcolor{LightCyan}      810     & 1.03e-08 & 1.41e-06 & 1.12e-06 & 1.29e-06 & 4.64 & 3.78 & 3.83 & 9  & 1 & 0.146 & 0.102 & 0.248 & 83     \\ 
\rowcolor{LightCyan}      3120    & 3.43e-10 & 9.1e-08 & 7.19e-08 & 8.3e-08 & 5.04 & 4.06 & 4.08 & 9  & 1 & 0.643 & 0.396 & 1.04 & 71.6      \\ 
\rowcolor{LightCyan}      12780   & 1.01e-11 & 5.41e-09 & 4.32e-09 & 4.97e-09 & 5 & 4 & 3.99 & 9  & 1 & 3.03 & 1.64 & 4.67 & 89              \\ 
\rowcolor{LightCyan}      50744   & 3.51e-13 & 3.45e-10 & 2.76e-10 & 3.18e-10 & 4.88 & 3.99 & 3.99 & 9  & 1 & 14.5 & 6.48 & 21 & 66.7        \\ 
\rowcolor{DardSeeGreen}      8       & 0.000871 & 0.0139 & 0.00786 & 2.48e-14 & - & - & -  & 4 & 1 & 0.00129 & 0.00106 & 0.00235 & -            \\ 
\rowcolor{DardSeeGreen}      50      & 1.57e-05 & 0.00053 & 0.000244 & 9.18e-14 & 4.38 & 3.56 & 3.79 & 8 & 1 & 0.00601 & 0.00594 & 0.0119 & 118 \\ 
\rowcolor{DardSeeGreen}      192     & 5.28e-07 & 3.56e-05 & 1.56e-05 & 2.49e-13 & 5.05 & 4.02 & 4.09 & 7 & 1 & 0.0213 & 0.0217 & 0.043 & 83.3  \\ 
\rowcolor{DardSeeGreen}      810     & 1.72e-08 & 2.31e-06 & 9.61e-07 & 1.14e-12 & 4.76 & 3.8 & 3.87 & 8 & 1 & 0.139 & 0.0908 & 0.23 & 74.8     \\ 
\rowcolor{DardSeeGreen}      3120    & 5.66e-10 & 1.53e-07 & 6.18e-08 & 4.66e-12 & 5.06 & 4.03 & 4.07 & 8 & 1 & 0.647 & 0.349 & 0.995 & 69.4    \\ 
\rowcolor{DardSeeGreen}      12780   & 1.65e-11 & 9e-09 & 3.67e-09 & 2.02e-11 & 5.01 & 4.02 & 4.01 & 8 & 1 & 3.41 & 1.42 & 4.83 & 82.4          \\ 
\rowcolor{DardSeeGreen}      50744   & 9.85e-13 & 5.83e-10 & 2.38e-10 & 7.65e-11 & 4.09 & 3.97 & 3.97 & 8 & 1 & 19.7 & 5.7 & 25.4 & 57.1        \\ 
\hline
\end{tabular}
\caption{Evaluation of $p$-multilevel solution strategies ($d{=}2$ standard meshes, $k{=}3$).
         See Section \ref{sec:numResOptions} for solver options. \label{tab:numResPerfReg2DhighO}}
\end{sidewaystable}

\begin{sidewaystable}
\small
\centering
\begin{tabular}{ c | c c  c c | c c c | c c | c c c | c }
\multicolumn{4}{c}{{\cellcolor{LightGoldenrod} \dg }} 
                                              & \multicolumn{5}{c}{{\cellcolor{LightCyan} \hhodp{} \vcond{} }}
                                              & \multicolumn{5}{c}{{\cellcolor{DardSeeGreen} \hhohp{} \vpcond{} }}\\
$\card{\Th}$ & \multicolumn{4}{c|}{error in $L^2$ norm} & \multicolumn{3}{c|}{conv. rate} & \multicolumn{2}{c|}{ITs} & \multicolumn{3}{c|}{CPU time} & Eff \\ \hline
 &  $\uVec_h$    & $G\uVec_h$   & $p_h$ & $D\uVec_h $& $\uVec_h$   & $G\uVec_h$ & $p_h$  & ITs & ITs$_L$ & Sol. & Ass. & Tot. & \\ \hline
 & \multicolumn{13}{l}{graded quadrilateral elements grid}\\
\hline
\rowcolor{LightGoldenrod} 4       &0.0128 & 0.154 & 0.0456 & 0.0411 & - & - & -  & 6 & 1 & 0.00162 & 0.000333 & 0.00195 & -               \\ 
\rowcolor{LightGoldenrod} 16      &0.00291 & 0.0487 & 0.0153 & 0.0145 & 2.13 & 1.66 & 1.57 & 8 & 1 & 0.00759 & 0.00126 & 0.00885 & 88.3   \\ 
\rowcolor{LightGoldenrod} 64      &0.000207 & 0.00656 & 0.00186 & 0.00204 & 3.81 & 2.89 & 3.04 & 9 & 1 & 0.0359 & 0.0051 & 0.041 & 86.4   \\ 
\rowcolor{LightGoldenrod} 256     &1.17e-05 & 0.000784 & 0.000178 & 0.000242 & 4.14 & 3.07 & 3.39 & 10 & 1 & 0.193 & 0.021 & 0.214 & 76.7 \\ 
\rowcolor{LightGoldenrod} 1024    &7.47e-07 & 9.94e-05 & 2.22e-05 & 3.02e-05 & 3.97 & 2.98 & 3 & 11 & 1 & 0.968 & 0.0849 & 1.05 & 81.2    \\ 
\rowcolor{LightGoldenrod} 4096    &4.75e-08 & 1.26e-05 & 2.83e-06 & 3.89e-06 & 3.97 & 2.98 & 2.97 & 12 & 1 & 5 & 0.343 & 5.34 & 78.8      \\ 
\rowcolor{LightGoldenrod} 16384   &3.01e-09 & 1.59e-06 & 3.62e-07 & 4.94e-07 & 3.98 & 2.99 & 2.97 & 13 & 1 & 27.2 & 1.38 & 28.6 & 74.9    \\ 
\rowcolor{LightCyan}      4      & 0.00201 & 0.0243 & 0.0117 & 0.0123 & - & - & -  & 5 & 1 & 0.00102 & 0.000814 & 0.00184 & -               \\
\rowcolor{LightCyan}      16     & 0.000296 & 0.00552 & 0.00292 & 0.00202 & 2.76 & 2.14 & 2 & 5 & 1 & 0.00244 & 0.00265 & 0.00508 & 145     \\
\rowcolor{LightCyan}      64     & 1.23e-05 & 0.000398 & 0.000225 & 0.000208 & 4.58 & 3.8 & 3.7 & 6 & 1 & 0.00923 & 0.0103 & 0.0196 & 104   \\
\rowcolor{LightCyan}      256    & 3.52e-07 & 2.31e-05 & 1.16e-05 & 1.18e-05 & 5.13 & 4.11 & 4.28 & 6 & 1 & 0.0406 & 0.0411 & 0.0817 & 95.8 \\
\rowcolor{LightCyan}      1024   & 1.12e-08 & 1.46e-06 & 7.49e-07 & 7.59e-07 & 4.97 & 3.98 & 3.95 & 6 & 1 & 0.232 & 0.165 & 0.396 & 82.5    \\
\rowcolor{LightCyan}      4096   & 3.57e-10 & 9.3e-08 & 4.85e-08 & 4.93e-08 & 4.97 & 3.97 & 3.95 & 6 & 1 & 1.13 & 0.658 & 1.79 & 88.5       \\
\rowcolor{LightCyan}      16384  & 1.16e-11 & 5.95e-09 & 3.2e-09 & 3.2e-09 & 4.94 & 3.97 & 3.92 & 7 & 1 & 6.07 & 2.69 & 8.75 & 81.9         \\
\rowcolor{DardSeeGreen}      4      & 0.0023 & 0.0315 & 0.00957 & 1.75e-14 & - & - & -  & 3 & 1 & 0.000833 & 0.000712 & 0.00154 & -          \\
\rowcolor{DardSeeGreen}      16     & 0.000387 & 0.0073 & 0.00187 & 1.01e-13 & 2.57 & 2.11 & 2.36 & 4 & 1 & 0.00422 & 0.00515 & 0.00937 & 66 \\
\rowcolor{DardSeeGreen}      64     & 1.49e-05 & 0.000514 & 0.000122 & 4.8e-13 & 4.7 & 3.83 & 3.93 & 4 & 1 & 0.00866 & 0.0112 & 0.0199 & 189 \\
\rowcolor{DardSeeGreen}      256    & 4.2e-07 & 2.98e-05 & 7.2e-06 & 5.05e-12 & 5.15 & 4.11 & 4.09 & 5 & 1 & 0.0489 & 0.0369 & 0.0858 & 92.6 \\
\rowcolor{DardSeeGreen}      1024   & 1.33e-08 & 1.88e-06 & 4.44e-07 & 3.18e-11 & 4.99 & 3.98 & 4.02 & 5 & 1 & 0.303 & 0.147 & 0.45 & 76.3   \\
\rowcolor{DardSeeGreen}      4096   & 4.29e-10 & 1.21e-07 & 2.9e-08 & 2.98e-10 & 4.95 & 3.96 & 3.94 & 5 & 1 & 1.74 & 0.589 & 2.33 & 77.1     \\
\rowcolor{DardSeeGreen}      16384  & 1.37e-11 & 8.04e-09 & 2.82e-09 & 2.26e-09 & 4.97 & 3.91 & 3.36 & 6 & 1 & 10.6 & 2.34 & 13 & 71.9       \\
\hline
& \multicolumn{13}{l}{graded triangular elements grid} \\
\hline
\rowcolor{LightGoldenrod}  8      & 0.00512 & 0.0804 & 0.0575 & 0.034 & - & - & -  & 8 & 1 & 0.0034 & 0.000539 & 0.00394 & -                \\
\rowcolor{LightGoldenrod}  32     & 0.00163 & 0.032 & 0.0213 & 0.0127 & 1.65 & 1.33 & 1.43 & 10 & 1 & 0.0152 & 0.00205 & 0.0173 & 91.1      \\
\rowcolor{LightGoldenrod}  128    & 0.000131 & 0.00454 & 0.00247 & 0.00206 & 3.64 & 2.82 & 3.11 & 11 & 1 & 0.0711 & 0.00833 & 0.0794 & 87.2 \\
\rowcolor{LightGoldenrod}  512    & 7.96e-06 & 0.000557 & 0.000262 & 0.000257 & 4.04 & 3.03 & 3.24 & 11 & 1 & 0.343 & 0.0339 & 0.377 & 84.2 \\
\rowcolor{LightGoldenrod}  2048   & 4.91e-07 & 6.9e-05 & 3.29e-05 & 3.19e-05 & 4.02 & 3.02 & 2.99 & 17 & 1 & 2.26 & 0.135 & 2.4 & 63        \\
\rowcolor{LightGoldenrod}  8192   & 3.15e-08 & 8.78e-06 & 4.11e-06 & 4.07e-06 & 3.96 & 2.97 & 3 & 31 & 1 & 17 & 0.548 & 17.6 & 54.5         \\
\rowcolor{LightGoldenrod}  32768  & 2.04e-09 & 1.12e-06 & 5.36e-07 & 5.18e-07 & 3.94 & 2.97 & 2.94 & 50 & 1 & 114 & 2.22 & 116 & 60.8       \\
\rowcolor{LightCyan}       8      & 0.000743 & 0.00999 & 0.00762 & 0.0083 & - & - & -  & 8 & 1 & 0.00154 & 0.00125 & 0.00279 & -            \\
\rowcolor{LightCyan}       32     & 0.000182 & 0.00321 & 0.00238 & 0.00286 & 2.03 & 1.64 & 1.68 & 8 & 1 & 0.00422 & 0.00415 & 0.00837 & 134 \\
\rowcolor{LightCyan}       128    & 7.9e-06 & 0.000253 & 0.000202 & 0.000228 & 4.53 & 3.67 & 3.56 & 8 & 1 & 0.0153 & 0.0163 & 0.0316 & 106  \\
\rowcolor{LightCyan}       512    & 2.36e-07 & 1.54e-05 & 1.13e-05 & 1.32e-05 & 5.06 & 4.04 & 4.17 & 8 & 1 & 0.0682 & 0.065 & 0.133 & 95    \\
\rowcolor{LightCyan}       2048   & 6.76e-09 & 9.11e-07 & 6.89e-07 & 7.85e-07 & 5.13 & 4.08 & 4.03 & 9 & 1 & 0.38 & 0.255 & 0.636 & 83.8    \\
\rowcolor{LightCyan}       8192   & 2.21e-10 & 5.85e-08 & 4.39e-08 & 5e-08 & 4.94 & 3.96 & 3.97 & 11 & 1 & 2.02 & 1.03 & 3.05 & 83.2        \\
\rowcolor{LightCyan}       32768  & 7.45e-12 & 3.81e-09 & 2.92e-09 & 3.27e-09 & 4.89 & 3.94 & 3.91 & 18 & 1 & 12.8 & 4.15 & 17 & 72.1       \\
\rowcolor{DardSeeGreen}       8      & 0.000919 & 0.0167 & 0.00719 & 3.62e-14 & - & - & -  & 5 & 1 & 0.00111 & 0.00104 & 0.00215 & -            \\
\rowcolor{DardSeeGreen}       32     & 0.000248 & 0.00587 & 0.00215 & 1.7e-13 & 1.89 & 1.51 & 1.74 & 6 & 1 & 0.00312 & 0.00374 & 0.00686 & 125  \\
\rowcolor{DardSeeGreen}       128    & 1.07e-05 & 0.000402 & 0.000157 & 1.05e-12 & 4.54 & 3.87 & 3.77 & 8 & 1 & 0.0141 & 0.0147 & 0.0288 & 95.3 \\
\rowcolor{DardSeeGreen}       512    & 3.01e-07 & 2.45e-05 & 9.62e-06 & 8.43e-12 & 5.15 & 4.03 & 4.03 & 9 & 1 & 0.078 & 0.058 & 0.136 & 84.6    \\
\rowcolor{DardSeeGreen}       2048   & 1.01e-08 & 1.62e-06 & 5.72e-07 & 5.98e-11 & 4.9 & 3.92 & 4.07 & 10 & 1 & 0.452 & 0.228 & 0.68 & 80       \\
\rowcolor{DardSeeGreen}       8192   & 3.24e-10 & 1.03e-07 & 3.68e-08 & 5.44e-10 & 4.96 & 3.97 & 3.96 & 14 & 1 & 2.79 & 0.91 & 3.7 & 73.6       \\
\rowcolor{DardSeeGreen}       32768  & 1.07e-11 & 9.91e-09 & 3.82e-09 & 3.94e-09 & 4.92 & 3.38 & 3.27 & 21 & 1 & 18.4 & 3.6 & 22 & 67.2         \\
\hline
\end{tabular}
\caption{Evaluation of $p$-multilevel solution strategies ($d{=}2$ graded meshes, $k{=}3$). 
         See Section \ref{sec:numResOptions} for solver options. \label{tab:numResPerfGrad2DhighO}}
\end{sidewaystable}

\begin{sidewaystable}
\small
\centering
\begin{tabular}{ c | c c  c c | c c c | c c | c c c | c }
\multicolumn{4}{c}{{\cellcolor{LightGoldenrod} \dg }} 
                                              & \multicolumn{5}{c}{{\cellcolor{LightCyan} \hhodp{} \vcond{} }}
                                              & \multicolumn{5}{c}{{\cellcolor{DardSeeGreen} \hhohp{} \vpcond{} }}\\
$\card{\Th}$ & \multicolumn{4}{c|}{error in $L^2$ norm} & \multicolumn{3}{c|}{conv. rate} & \multicolumn{2}{c|}{ITs} & \multicolumn{3}{c|}{CPU time} & Eff \\ \hline
 &  $\uVec_h$    & $G\uVec_h$   & $p_h$ & $D\uVec_h $& $\uVec_h$   & $G\uVec_h$ & $p_h$  & ITs & ITs$_L$ & Sol. & Ass. & Tot. & \\ \hline
&	\multicolumn{13}{l}{trapezoidal elements grid}\\
\hline
\rowcolor{LightGoldenrod} 4      & 1.16e-05 & 0.000255 & 0.000132 & 0.000108 & - & - & -  & 9 & 1 & 0.0105 & 0.0016 & 0.0121 & -           \\
\rowcolor{LightGoldenrod} 16     & 1.35e-07 & 5.43e-06 & 1.79e-06 & 2.41e-06 & 6.43 & 5.55 & 6.2 & 13 & 1 & 0.0714 & 0.00682 & 0.0782 & 62 \\
\rowcolor{LightGoldenrod} 64     & 1.23e-09 & 9.88e-08 & 3.33e-08 & 4.32e-08 & 6.77 & 5.78 & 5.75 & 18 & 1 & 0.467 & 0.0287 & 0.496 & 63   \\
\rowcolor{LightGoldenrod} 256    & 1.02e-11 & 1.61e-09 & 5.25e-10 & 7.03e-10 & 6.92 & 5.94 & 5.99 & 21 & 1 & 2.32 & 0.119 & 2.44 & 81.3    \\
\rowcolor{LightGoldenrod} 1024   & 3.19e-13 & 2.96e-11 & 4.06e-11 & 1.57e-11 & 5 & 5.77 & 3.69 & 21 & 1 & 9.69 & 0.472 & 10.2 & 96.1       \\
\rowcolor{LightCyan}      4      & 1.12e-06 & 2.34e-05 & 1.77e-05 & 1.83e-05 & - & - & -  & 7 & 1 & 0.00246 & 0.00321 & 0.00566 & -        \\  
\rowcolor{LightCyan}      16     & 7.95e-09 & 2.93e-07 & 2.17e-07 & 2.4e-07 & 7.14 & 6.32 & 6.36 & 7 & 1 & 0.00828 & 0.0131 & 0.0214 & 106 \\  
\rowcolor{LightCyan}      64     & 3.53e-11 & 2.64e-09 & 1.93e-09 & 2.11e-09 & 7.82 & 6.8 & 6.81 & 8 & 1 & 0.0375 & 0.0535 & 0.091 & 93.8  \\  
\rowcolor{LightCyan}      256    & 1.51e-13 & 2.18e-11 & 1.59e-11 & 1.75e-11 & 7.86 & 6.92 & 6.93 & 8 & 1 & 0.179 & 0.217 & 0.396 & 91.9   \\  
\rowcolor{LightCyan}      1024   & 2.2e-13 & 2.02e-12 & 2.11e-12 & 1.47e-12 & -0.539 & 3.43 & 2.91 & 8 & 1 & 0.818 & 0.882 & 1.7 & 93.2    \\  
\rowcolor{DardSeeGreen}      4      & 1.58e-06 & 3.45e-05 & 1.47e-05 & 1.11e-13 & - & - & -  & 3 & 1 & 0.00142 & 0.00317 & 0.00459 & -        \\
\rowcolor{DardSeeGreen}      16     & 1.03e-08 & 4.28e-07 & 1.7e-07 & 2.62e-13 & 7.26 & 6.33 & 6.43 & 4 & 1 & 0.00491 & 0.0118 & 0.0167 & 110 \\
\rowcolor{DardSeeGreen}      64     & 4.91e-11 & 3.87e-09 & 1.5e-09 & 5.62e-13 & 7.71 & 6.79 & 6.83 & 5 & 1 & 0.0225 & 0.0464 & 0.0689 & 97.2 \\
\rowcolor{DardSeeGreen}      256    & 2.39e-13 & 3.25e-11 & 1.18e-11 & 2.07e-12 & 7.68 & 6.9 & 6.99 & 6 & 1 & 0.133 & 0.184 & 0.317 & 86.8    \\
\rowcolor{DardSeeGreen}      1024   & 1.03e-13 & 8.05e-12 & 5.91e-12 & 7.14e-12 & 1.22 & 2.01 & 0.997 & 8 & 1 & 0.783 & 0.734 & 1.52 & 83.6   \\
\hline
&	\multicolumn{13}{l}{delaunay triangular grid} \\
\hline
\rowcolor{LightGoldenrod} 8       &1.87e-06 & 4.97e-05 & 3.39e-05 & 2.16e-05 & - & - & -  & 20 & 1 & 0.0378 & 0.00288 & 0.0407 & -         \\ 
\rowcolor{LightGoldenrod} 50      &4.01e-09 & 2.73e-07 & 1.71e-07 & 1.36e-07 & 6.7 & 5.68 & 5.78 & 27 & 1 & 0.416 & 0.0184 & 0.435 & 58.5  \\ 
\rowcolor{LightGoldenrod} 192     &2.95e-11 & 4.11e-09 & 2.02e-09 & 2.08e-09 & 7.3 & 6.24 & 6.6 & 30 & 1 & 1.89 & 0.0736 & 1.96 & 85.2     \\ 
\rowcolor{LightGoldenrod} 810     &5.31e-13 & 7.62e-11 & 5.04e-11 & 3.75e-11 & 5.58 & 5.54 & 5.13 & 30 & 1 & 8.26 & 0.313 & 8.57 & 96.4    \\ 
\rowcolor{LightGoldenrod} 3120    &7.33e-13 & 2.44e-11 & 5.41e-11 & 1.58e-11 & -0.478 & 1.69 & -0.106 & 31 & 1 & 33.5 & 1.21 & 34.7 & 95.1 \\ 
\rowcolor{LightCyan}      8       &1.3e-07 & 3.21e-06 & 2.84e-06 & 2.98e-06 & - & - & -  & 11 & 1 & 0.0047 & 0.00569 & 0.0104 & -         \\ 
\rowcolor{LightCyan}      50      &1.5e-10 & 8.28e-09 & 7.71e-09 & 8.33e-09 & 7.38 & 6.5 & 6.45 & 11 & 1 & 0.026 & 0.0355 & 0.0616 & 101  \\ 
\rowcolor{LightCyan}      192     &5.39e-13 & 6.25e-11 & 5.6e-11 & 6.2e-11 & 8.36 & 7.26 & 7.32 & 12 & 1 & 0.123 & 0.137 & 0.26 & 71      \\ 
\rowcolor{LightCyan}      810     &1.24e-13 & 1.9e-12 & 2.25e-12 & 1.61e-12 & 2.04 & 4.86 & 4.47 & 12 & 1 & 0.603 & 0.58 & 1.18 & 88      \\ 
\rowcolor{LightCyan}      3120    &2.47e-13 & 3.52e-12 & 4.51e-12 & 2.7e-12 & -1.02 & -0.917 & -1.03 & 12 & 1 & 2.44 & 2.24 & 4.68 & 75.9 \\ 
\rowcolor{DardSeeGreen}      8       &2.3e-07 & 6.83e-06 & 3.44e-06 & 3.16e-13 & - & - & -  & 4 & 1 & 0.00232 & 0.00625 & 0.00857 & -         \\
\rowcolor{DardSeeGreen}      50      &2.61e-10 & 1.93e-08 & 8.2e-09 & 8.56e-13 & 7.4 & 6.41 & 6.59 & 12 & 1 & 0.0172 & 0.0312 & 0.0484 & 106  \\
\rowcolor{DardSeeGreen}      192     &9.83e-13 & 1.44e-10 & 6.19e-11 & 1.51e-12 & 8.3 & 7.28 & 7.26 & 12 & 1 & 0.08 & 0.119 & 0.199 & 73      \\
\rowcolor{DardSeeGreen}      810     &1.05e-13 & 9.71e-12 & 1.02e-11 & 6.33e-12 & 3.11 & 3.74 & 2.51 & 13 & 1 & 0.439 & 0.502 & 0.941 & 84.7  \\
\rowcolor{DardSeeGreen}      3120    &3.64e-13 & 3.18e-11 & 1.75e-11 & 2.44e-11 & -1.85 & -1.76 & -0.806 & 14 & 1 & 1.92 & 1.92 & 3.84 & 73.5 \\
\hline
\end{tabular}
\caption{Evaluation of $p$-multilevel solution strategies for solving higher-order ($d=2$ standard meshes, $k{=}6$).
         See Section \ref{sec:numResOptions} for solver options. \label{tab:numResPerfReg2DhigherO}}
\end{sidewaystable}

\begin{sidewaystable}
\small
\centering
\begin{tabular}{ c | c c  c c | c c c | c c | c c c | c }
\multicolumn{4}{c}{{\cellcolor{LightGoldenrod} \dg }} 
                                              & \multicolumn{5}{c}{{\cellcolor{LightCyan} \hhodp{} \vcond{} }}
                                              & \multicolumn{5}{c}{{\cellcolor{DardSeeGreen} \hhohp{} \vpcond{} }}\\
$\card{\Th}$ & \multicolumn{4}{c|}{error in $L^2$ norm} & \multicolumn{3}{c|}{conv. rate} & \multicolumn{2}{c|}{ITs} & \multicolumn{3}{c|}{CPU time} & Eff \\ \hline
 &  $\uVec_h$    & $G\uVec_h$   & $p_h$ & $D\uVec_h $& $\uVec_h$   & $G\uVec_h$ & $p_h$  & ITs & ITs$_L$ & Sol. & Ass. & Tot. & \\ \hline
 & \multicolumn{13}{l}{graded quadrilateral elements grid}\\
\hline
\rowcolor{LightGoldenrod} 4       &1.16e-05 & 0.000255 & 0.000132 & 0.000108 & - & - & -  & 9 & 1 & 0.0105 & 0.0016 & 0.0121 & -             \\ 
\rowcolor{LightGoldenrod} 16      &1.18e-06 & 3.32e-05 & 1.17e-05 & 1.39e-05 & 3.31 & 2.94 & 3.5 & 13 & 1 & 0.0714 & 0.00682 & 0.0782 & 61.9 \\ 
\rowcolor{LightGoldenrod} 64      &1.27e-08 & 6.61e-07 & 2.04e-07 & 2.82e-07 & 6.53 & 5.65 & 5.83 & 17 & 1 & 0.448 & 0.0287 & 0.477 & 65.6   \\ 
\rowcolor{LightGoldenrod} 256     &8.41e-11 & 9.01e-09 & 2.42e-09 & 3.86e-09 & 7.24 & 6.2 & 6.4 & 19 & 1 & 2.13 & 0.117 & 2.24 & 84.9        \\ 
\rowcolor{LightGoldenrod} 1024    &2.64e-12 & 1.84e-10 & 3.27e-10 & 1.04e-10 & 4.99 & 5.61 & 2.89 & 22 & 1 & 10 & 0.472 & 10.5 & 85.3        \\ 
\rowcolor{LightCyan}      4      & 1.12e-06 & 2.34e-05 & 1.77e-05 & 1.83e-05 & - & - & -  & 7 & 1 & 0.00247 & 0.00319 & 0.00565 & -         \\ 
\rowcolor{LightCyan}      16     & 8.24e-08 & 2.39e-06 & 1.65e-06 & 1.81e-06 & 3.77 & 3.29 & 3.42 & 7 & 1 & 0.0083 & 0.0132 & 0.0215 & 105  \\ 
\rowcolor{LightCyan}      64     & 5.34e-10 & 2.68e-08 & 1.84e-08 & 2.12e-08 & 7.27 & 6.48 & 6.49 & 7 & 1 & 0.0347 & 0.0536 & 0.0883 & 97.2 \\ 
\rowcolor{LightCyan}      256    & 1.8e-12 & 1.78e-10 & 1.33e-10 & 1.47e-10 & 8.21 & 7.23 & 7.11 & 7 & 1 & 0.164 & 0.218 & 0.382 & 92.5     \\ 
\rowcolor{LightCyan}      1024   & 5.17e-13 & 3.06e-11 & 3.39e-11 & 2.74e-11 & 1.8 & 2.54 & 1.97 & 7 & 1 & 0.754 & 0.882 & 1.64 & 93.4      \\ 
\rowcolor{DardSeeGreen}      4      & 1.58e-06 & 3.45e-05 & 1.47e-05 & 1.11e-13 & - & - & -  & 3 & 1 & 0.00138 & 0.00319 & 0.00458 & -         \\
\rowcolor{DardSeeGreen}      16     & 1.2e-07 & 3.6e-06 & 1.23e-06 & 5.14e-13 & 3.71 & 3.26 & 3.58 & 4 & 1 & 0.00489 & 0.0119 & 0.0168 & 109   \\
\rowcolor{DardSeeGreen}      64     & 6.28e-10 & 3.67e-08 & 1.39e-08 & 3.44e-12 & 7.58 & 6.61 & 6.47 & 5 & 1 & 0.0222 & 0.0465 & 0.0687 & 97.6 \\
\rowcolor{DardSeeGreen}      256    & 2.96e-12 & 3.02e-10 & 9.68e-11 & 2.5e-11 & 7.73 & 6.93 & 7.17 & 6 & 1 & 0.13 & 0.184 & 0.314 & 87.6      \\
\rowcolor{DardSeeGreen}      1024   & 4.74e-13 & 1.97e-10 & 1.67e-10 & 1.76e-10 & 2.64 & 0.616 & -0.784 & 8 & 1 & 0.787 & 0.734 & 1.52 & 82.5  \\
\hline
& \multicolumn{13}{l}{graded triangular elements grid} \\
\hline
\rowcolor{LightGoldenrod}  8      &2.73e-06 & 5.64e-05 & 4.58e-05 & 2.39e-05 & - & - & -  & 14 & 1 & 0.0281 & 0.00283 & 0.0309 & -         \\ 
\rowcolor{LightGoldenrod}  32     &3.83e-07 & 1.01e-05 & 7.37e-06 & 4.17e-06 & 2.83 & 2.48 & 2.64 & 18 & 1 & 0.165 & 0.0116 & 0.176 & 70.1 \\ 
\rowcolor{LightGoldenrod}  128    &3.69e-09 & 2.24e-07 & 1.33e-07 & 1.02e-07 & 6.7 & 5.5 & 5.79 & 23 & 1 & 0.943 & 0.048 & 0.991 & 71.2    \\ 
\rowcolor{LightGoldenrod}  512    &3.7e-11 & 3.8e-09 & 2.42e-09 & 1.7e-09 & 6.64 & 5.88 & 5.78 & 23 & 1 & 4.01 & 0.195 & 4.21 & 94.2       \\ 
\rowcolor{LightGoldenrod}  2048   &6.96e-13 & 8.9e-11 & 1.71e-10 & 5.37e-11 & 5.73 & 5.42 & 3.83 & 36 & 1 & 24.5 & 0.784 & 25.2 & 66.7     \\ 
\rowcolor{LightCyan}       8      & 1.73e-07 & 3.84e-06 & 3.4e-06 & 3.57e-06 & - & - & -  & 10 & 1 & 0.00395 & 0.00548 & 0.00943 & -         \\
\rowcolor{LightCyan}       32     & 2.06e-08 & 5.65e-07 & 4.75e-07 & 5.12e-07 & 3.07 & 2.77 & 2.84 & 11 & 1 & 0.0151 & 0.0222 & 0.0373 & 101 \\
\rowcolor{LightCyan}       128    & 1.26e-10 & 6.84e-09 & 6.13e-09 & 6.74e-09 & 7.35 & 6.37 & 6.28 & 11 & 1 & 0.0688 & 0.0898 & 0.159 & 94.2 \\
\rowcolor{LightCyan}       512    & 6.57e-13 & 6.17e-11 & 5.23e-11 & 5.71e-11 & 7.59 & 6.79 & 6.87 & 12 & 1 & 0.353 & 0.362 & 0.715 & 88.7   \\
\rowcolor{LightCyan}       2048   & 8.53e-13 & 2.42e-11 & 2.6e-11 & 2.1e-11 & -0.376 & 1.35 & 1.01 & 13 & 1 & 1.63 & 1.46 & 3.09 & 92.5      \\
\rowcolor{DardSeeGreen}       8      & 3.39e-07 & 1.03e-05 & 4.11e-06 & 3.53e-13 & - & - & -  & 6 & 1 & 0.002 & 0.00524 & 0.00724 & -         \\
\rowcolor{DardSeeGreen}       32     & 2.63e-08 & 1.17e-06 & 5.14e-07 & 4.38e-12 & 3.69 & 3.13 & 3 & 8 & 1 & 0.00783 & 0.0201 & 0.0279 & 104  \\
\rowcolor{DardSeeGreen}       128    & 2.91e-10 & 2.01e-08 & 7.05e-09 & 4.89e-11 & 6.5 & 5.87 & 6.19 & 9 & 1 & 0.0354 & 0.0795 & 0.115 & 97.2 \\
\rowcolor{DardSeeGreen}       512    & 1.06e-12 & 3.21e-10 & 1.14e-10 & 1.93e-10 & 8.1 & 5.97 & 5.95 & 12 & 1 & 0.234 & 0.316 & 0.549 & 83.6  \\
\rowcolor{DardSeeGreen}       2048   & 4.71e-12 & 1.21e-08 & 1.33e-09 & 2.48e-09 & -2.15 & -5.23 & -3.54 & 15 & 1 & 1.31 & 1.26 & 2.57 & 85.6 \\
\hline
\end{tabular}
\caption{Evaluation of $p$-multilevel solution strategies for solving higher-order ($d{=}2$ graded meshes, $k{=}6$).
         See Section \ref{sec:numResOptions} for solver options. \label{tab:numResPerfGrad2DhigherO}}
\end{sidewaystable}

\begin{sidewaystable}
\small
\centering
\begin{tabular}{ c | c c  c c | c c c | c c | c c c | c }
\multicolumn{4}{c}{{\cellcolor{LightGoldenrod} \dg }} 
                                              & \multicolumn{5}{c}{{\cellcolor{LightCyan} \hhodp{} \vcond{} }}
                                              & \multicolumn{5}{c}{{\cellcolor{DardSeeGreen} \hhohp{} \vpcond{} }}\\
$\card{\Th}$ & \multicolumn{4}{c|}{error in $L^2$ norm} & \multicolumn{3}{c|}{conv. rate} & \multicolumn{2}{c|}{ITs} & \multicolumn{3}{c|}{CPU time} & Eff \\ \hline
 &  $\uVec_h$    & $G\uVec_h$   & $p_h$ & $D\uVec_h $& $\uVec_h$   & $G\uVec_h$ & $p_h$  & ITs & ITs$_L$ & Sol. & Ass. & Tot. & \\ \hline
& \multicolumn{12}{l}{prismatic elements grid} \\
\hline
\rowcolor{LightGoldenrod} 16   & 0.00575 & 0.129 & 0.0988 & 0.0365 & - & - & -  & 9 & 10 & 0.046 & 0.0111 & 0.0571 & -               \\
\rowcolor{LightGoldenrod} 128  & 0.000347 & 0.0162 & 0.00846 & 0.00538 & 4.05 & 2.99 & 3.55 & 12 & 22 & 0.66 & 0.0873 & 0.747 & 61.2 \\
\rowcolor{LightGoldenrod} 1024 & 2.15e-05 & 0.00202 & 0.000887 & 0.000707 & 4.01 & 3 & 3.25 & 13 & 61 & 7.31 & 0.698 & 8.01 & 74.6   \\
\rowcolor{LightGoldenrod} 8192 & 1.34e-06 & 0.000252 & 0.000102 & 9e-05 & 4 & 3 & 3.12 & 15 & 160 & 105 & 5.65 & 111 & 57.9          \\
\rowcolor{LightCyan}      16   & 0.000721 & 0.0162 & 0.0128 & 0.00713 & - & - & -  & 8 & 8 & 0.0253 & 0.0342 & 0.0595 & -             \\ 
\rowcolor{LightCyan}      128  & 2.41e-05 & 0.00105 & 0.000732 & 0.00047 & 4.91 & 3.95 & 4.13 & 9 & 18 & 0.291 & 0.284 & 0.575 & 82.7 \\ 
\rowcolor{LightCyan}      1024 & 7.77e-07 & 6.57e-05 & 4.21e-05 & 2.93e-05 & 4.95 & 3.99 & 4.12 & 9 & 42 & 2.95 & 2.35 & 5.3 & 86.8   \\ 
\rowcolor{LightCyan}      8192 & 2.47e-08 & 4.11e-06 & 2.5e-06 & 1.82e-06 & 4.97 & 4 & 4.08 & 10 & 105 & 42.2 & 19.1 & 61.3 & 69.1    \\ 
\rowcolor{DardSeeGreen}      16   & 0.000989 & 0.0307 & 0.0122 & 3.46e-14 & - & - & -  & 5 & 5 & 0.0766 & 0.0359 & 0.112 & -             \\
\rowcolor{DardSeeGreen}      128  & 3.11e-05 & 0.00199 & 0.000697 & 1.37e-13 & 4.99 & 3.95 & 4.13 & 7 & 12 & 0.851 & 0.279 & 1.13 & 79.6 \\
\rowcolor{DardSeeGreen}      1024 & 9.65e-07 & 0.000127 & 4.03e-05 & 5.08e-13 & 5.01 & 3.97 & 4.11 & 8 & 32 & 8.45 & 2.2 & 10.7 & 84.9   \\
\rowcolor{DardSeeGreen}      8192 & 3e-08 & 7.97e-06 & 2.4e-06 & 1.84e-12 & 5.01 & 3.99 & 4.07 & 8 & 93 & 87.7 & 17.6 & 105 & 80.9       \\
\hline
& \multicolumn{13}{l}{pyramidal elements grid}\\
\hline	
\rowcolor{LightGoldenrod} 48    & 0.00257 & 0.0765 & 0.0396 & 0.0221 & - & - & -  & 11 & 11 & 0.261 & 0.0329 & 0.294 & -            \\
\rowcolor{LightGoldenrod} 384   & 0.000164 & 0.00979 & 0.0044 & 0.00291 & 3.98 & 2.97 & 3.17 & 13 & 28 & 2.65 & 0.266 & 2.91 & 80.8 \\
\rowcolor{LightGoldenrod} 3072  & 1.02e-05 & 0.00123 & 0.000522 & 0.000371 & 4 & 3 & 3.08 & 13 & 62 & 25 & 2.12 & 27.1 & 85.9       \\
\rowcolor{LightGoldenrod} 24576 & 6.38e-07 & 0.000153 & 6.45e-05 & 4.68e-05 & 4 & 3 & 3.02 & 14 & 155 & 317 & 16.8 & 333 & 65.1     \\
\rowcolor{LightCyan}      48    & 0.000278 & 0.00811 & 0.00306 & 0.0029 & - & - & -  & 10 & 13 & 0.1 & 0.111 & 0.211 & -             \\  
\rowcolor{LightCyan}      384   & 8.93e-06 & 0.000519 & 0.000191 & 0.000193 & 4.96 & 3.97 & 4 & 10 & 31 & 1.05 & 0.898 & 1.95 & 86.6 \\  
\rowcolor{LightCyan}      3072  & 2.82e-07 & 3.26e-05 & 1.2e-05 & 1.23e-05 & 4.99 & 3.99 & 4 & 10 & 64 & 11.6 & 7.28 & 18.9 & 82.6   \\  
\rowcolor{LightCyan}      24576 & 8.82e-09 & 2.04e-06 & 7.48e-07 & 7.79e-07 & 5 & 4 & 4 & 10 & 146 & 169 & 58.5 & 228 & 66.4         \\  
\rowcolor{DardSeeGreen}      48    & 0.000302 & 0.0149 & 0.00335 & 1.19e-13 & - & - & -  & 11 & 16 & 0.447 & 0.103 & 0.55 & -            \\
\rowcolor{DardSeeGreen}      384   & 1.04e-05 & 0.000991 & 0.000203 & 4.2e-13 & 4.85 & 3.92 & 4.05 & 12 & 32 & 4.31 & 0.813 & 5.12 & 86  \\
\rowcolor{DardSeeGreen}      3072  & 3.4e-07 & 6.34e-05 & 1.25e-05 & 2.31e-12 & 4.94 & 3.97 & 4.02 & 12 & 76 & 44.7 & 6.47 & 51.2 & 80   \\
\rowcolor{DardSeeGreen}      24576 & 1.08e-08 & 4.01e-06 & 7.74e-07 & 8.84e-12 & 4.97 & 3.98 & 4.01 & 12 & 201 & 733 & 51.6 & 785 & 52.2 \\
\hline
& \multicolumn{12}{l}{graded tetrahedral elements grid} \\
\hline
\rowcolor{LightGoldenrod} 24    & 0.00992 & 0.215 & 0.13 & 0.0553 & - & - & -  & 14 & 17 & 0.1 & 0.0141 & 0.114 & -                  \\
\rowcolor{LightGoldenrod} 192   & 0.000671 & 0.0273 & 0.0125 & 0.00731 & 3.89 & 2.98 & 3.38 & 16 & 33 & 1.16 & 0.112 & 1.27 & 72.3   \\
\rowcolor{LightGoldenrod} 1536  & 0.000164 & 0.00915 & 0.00423 & 0.00239 & 2.03 & 1.57 & 1.56 & 19 & 63 & 13.2 & 0.879 & 14.1 & 71.9 \\
\rowcolor{LightGoldenrod} 12288 & 1.24e-05 & 0.00129 & 0.000646 & 0.000326 & 3.73 & 2.82 & 2.71 & 23 & 98 & 181 & 7.03 & 188 & 59.9  \\
\rowcolor{LightCyan}      24    & 0.00156 & 0.0342 & 0.0153 & 0.0166 & - & - & -  & 11 & 12 & 0.0335 & 0.0436 & 0.077 & -              \\
\rowcolor{LightCyan}      192   & 4.82e-05 & 0.00213 & 0.000893 & 0.00105 & 5.01 & 4.01 & 4.1 & 11 & 20 & 0.371 & 0.362 & 0.734 & 84   \\
\rowcolor{LightCyan}      1536  & 7.23e-06 & 0.000458 & 0.00016 & 0.000228 & 2.74 & 2.22 & 2.48 & 13 & 45 & 4.45 & 2.9 & 7.34 & 79.9   \\
\rowcolor{LightCyan}      12288 & 2.79e-07 & 3.31e-05 & 1.17e-05 & 1.64e-05 & 4.69 & 3.79 & 3.77 & 18 & 108 & 86.6 & 23.5 & 110 & 53.4 \\
\rowcolor{DardSeeGreen}      24    & 0.00165 & 0.0624 & 0.0137 & 5.37e-14 & - & - & -  & 6 & 5 & 0.0721 & 0.0443 & 0.116 & -              \\
\rowcolor{DardSeeGreen}      192   & 5.65e-05 & 0.00412 & 0.000728 & 1.95e-13 & 4.87 & 3.92 & 4.24 & 8 & 10 & 0.777 & 0.352 & 1.13 & 82.5 \\
\rowcolor{DardSeeGreen}      1536  & 8.31e-06 & 0.000905 & 0.000142 & 1.31e-12 & 2.77 & 2.19 & 2.36 & 11 & 27 & 8.86 & 2.74 & 11.6 & 77.9 \\
\rowcolor{DardSeeGreen}      12288 & 3.25e-07 & 6.77e-05 & 1.04e-05 & 1.06e-11 & 4.68 & 3.74 & 3.77 & 14 & 54 & 105 & 21.8 & 127 & 73.1   \\
\hline
\end{tabular}
\caption{Evaluation of $p$-multilevel solution strategies ($d{=}3$, $k{=}3$). 
         See Section \ref{sec:numResOptions} for solver options. \label{tab:numResPerf3DhighO}}
\end{sidewaystable}

\begin{sidewaystable}[ht]
\small
\centering
\begin{tabular}{ c | c c  c c | c c c | c c | c c c | c }
\multicolumn{4}{c}{{\cellcolor{LightGoldenrod} \dg}} 
                                              & \multicolumn{5}{c}{{\cellcolor{LightCyan} \hhodp{} \vcond{} }}
                                              & \multicolumn{5}{c}{{\cellcolor{DardSeeGreen} \hhohp{} \vpcond{} }}\\
$\card{\Th}$ & \multicolumn{4}{c|}{error in $L^2$ norm} & \multicolumn{3}{c|}{conv. rate} & \multicolumn{2}{c|}{ITs} & \multicolumn{3}{c|}{CPU time} & Eff \\ \hline
 &  $\uVec_h$    & $G\uVec_h$   & $p_h$ & $D\uVec_h $& $\uVec_h$   & $G\uVec_h$ & $p_h$  & ITs & ITs$_L$ & Sol. & Ass. & Tot. & \\ \hline
& \multicolumn{12}{l}{prismatic elements grid} \\
\hline
\rowcolor{LightGoldenrod} 16   & 2.38e-06 & 0.000105 & 6.49e-05 & 3.21e-05 & - & - & -  & 13 & 10 & 1.01 & 0.347 & 1.36 & -          \\
\rowcolor{LightGoldenrod} 128  & 1.87e-08 & 1.62e-06 & 8.46e-07 & 4.58e-07 & 6.99 & 6.02 & 6.26 & 17 & 25 & 12.7 & 2.8 & 15.5 & 69.8 \\
\rowcolor{LightGoldenrod} 1024 & 1.43e-10 & 2.48e-08 & 1.16e-08 & 6.25e-09 & 7.02 & 6.03 & 6.19 & 21 & 76 & 133 & 22.5 & 156 & 79.6  \\
\rowcolor{LightCyan}      16   & 1.47e-06 & 4.06e-05 & 4.74e-05 & 2.61e-05 & - & - & -  & 11 & 10 & 0.417 & 0.915 & 1.33 & -          \\
\rowcolor{LightCyan}      128  & 6.19e-09 & 3.43e-07 & 3.97e-07 & 2.11e-07 & 7.89 & 6.89 & 6.9 & 12 & 20 & 3.76 & 7.85 & 11.6 & 91.8  \\
\rowcolor{LightCyan}      1024 & 2.48e-11 & 2.72e-09 & 3.15e-09 & 1.69e-09 & 7.97 & 6.98 & 6.98 & 12 & 48 & 30.3 & 65.2 & 95.5 & 97.3 \\
\rowcolor{DardSeeGreen}      16   & 2.73e-06 & 0.000122 & 4.26e-05 & 3.39e-13 & - & - & -  & 7 & 4 & 1.26 & 0.96 & 2.22 & -            \\
\rowcolor{DardSeeGreen}      128  & 1.2e-08 & 1.06e-06 & 3.65e-07 & 9.84e-13 & 7.83 & 6.84 & 6.87 & 9 & 14 & 12.1 & 7.55 & 19.6 & 90.4 \\
\rowcolor{DardSeeGreen}      1024 & 4.89e-11 & 8.59e-09 & 2.93e-09 & 2.87e-12 & 7.94 & 6.95 & 6.96 & 10 & 38 & 105 & 60 & 165 & 95     \\
\hline
& \multicolumn{13}{l}{pyramidal elements grid}\\
\hline	
\rowcolor{LightGoldenrod} 48   & 9.56e-07 & 5.26e-05 & 2.26e-05 & 1.45e-05 & - & - & -  & 15 & 12 & 5.5 & 1.06 & 6.56 & -           \\ 
\rowcolor{LightGoldenrod} 384  & 7.53e-09 & 8.31e-07 & 3.69e-07 & 2.25e-07 & 6.99 & 5.98 & 5.94 & 21 & 29 & 56 & 8.47 & 64.5 & 81.4 \\ 
\rowcolor{LightGoldenrod} 3072 & 5.86e-11 & 1.3e-08 & 5.67e-09 & 3.47e-09 & 7.01 & 6 & 6.02 & 22 & 70 & 479 & 67.9 & 547 & 94.3     \\ 
\rowcolor{LightCyan} 48        & 1.18e-07 & 4.89e-06 & 4.58e-06 & 2.11e-06 & - & - & -  & 14 & 13 & 1.44 & 3.07 & 4.52 & -            \\
\rowcolor{LightCyan} 384       & 5.01e-10 & 4.06e-08 & 3.66e-08 & 1.86e-08 & 7.88 & 6.91 & 6.96 & 13 & 33 & 11.3 & 24.9 & 36.2 & 99.8 \\
\rowcolor{LightCyan} 3072      & 2.02e-12 & 3.23e-10 & 2.87e-10 & 1.54e-10 & 7.96 & 6.98 & 7 & 14 & 55 & 98.9 & 203 & 302 & 96        \\
\rowcolor{DardSeeGreen} 48        & 2.14e-07 & 1.39e-05 & 4.28e-06 & 1.08e-12 & - & - & -  & 17 & 13 & 6.36 & 2.81 & 9.17 & -            \\
\rowcolor{DardSeeGreen} 384       & 9.19e-10 & 1.18e-07 & 3.44e-08 & 3.22e-12 & 7.86 & 6.89 & 6.96 & 20 & 37 & 55.9 & 22.3 & 78.3 & 93.7 \\
\rowcolor{DardSeeGreen} 3072      & 3.75e-12 & 9.53e-10 & 2.72e-10 & 8.67e-12 & 7.94 & 6.95 & 6.98 & 20 & 76 & 467 & 178 & 645 & 97      \\
\hline
& \multicolumn{12}{l}{graded tetrahedral elements grid} \\
\hline
\rowcolor{LightGoldenrod} 24    & 1.38e-05 & 0.000562 & 0.000273 & 0.00021 & - & - & -  & 20 & 19 & 2.04 & 0.458 & 2.5 & -           \\
\rowcolor{LightGoldenrod} 192   & 1.16e-07 & 8.34e-06 & 3.94e-06 & 2.86e-06 & 6.89 & 6.07 & 6.11 & 30 & 36 & 25.3 & 3.66 & 29 & 69.1 \\
\rowcolor{LightGoldenrod} 1536  & 1.11e-08 & 1.03e-06 & 4.97e-07 & 3.46e-07 & 3.39 & 3.02 & 2.99 & 33 & 61 & 236 & 29.2 & 265 & 87.6 \\
\rowcolor{LightCyan} 24         & 1.87e-06 & 6.14e-05 & 4.26e-05 & 3.44e-05 & - & - & -  & 14 & 11 & 0.517 & 1.25 & 1.77 & -          \\
\rowcolor{LightCyan} 192        & 9.04e-09 & 5.61e-07 & 4.41e-07 & 3.52e-07 & 7.7 & 6.78 & 6.59 & 15 & 22 & 4.39 & 10.5 & 14.9 & 95.2 \\
\rowcolor{LightCyan} 1536       & 5.31e-10 & 4.58e-08 & 3.13e-08 & 3.18e-08 & 4.09 & 3.61 & 3.81 & 20 & 56 & 44.3 & 85.5 & 130 & 91.7 \\
\rowcolor{DardSeeGreen} 24         & 3.63e-06 & 0.000222 & 3.96e-05 & 1.15e-12 & - & - & -  & 9 & 5 & 1.07 & 1.25 & 2.32 & -             \\
\rowcolor{DardSeeGreen} 192        & 1.96e-08 & 2.2e-06 & 4.21e-07 & 8.59e-11 & 7.53 & 6.66 & 6.56 & 11 & 12 & 9.88 & 9.93 & 19.8 & 93.8 \\
\rowcolor{DardSeeGreen} 1536       & 9.6e-10 & 1.6e-07 & 2.91e-08 & 1.39e-10 & 4.35 & 3.78 & 3.85 & 14 & 26 & 91.8 & 78.9 & 171 & 92.9   \\
\hline
\end{tabular}
\caption{Evaluation of $p$-multilevel solution strategies ($d{=}3, k{=}3$).
         See Section \ref{sec:numResOptions} for solver options. \label{tab:numResPerf3DhigherO}}
\end{sidewaystable}

\subsection{Scalability}

In this section we include basic scalability results for p-multilevel solvers applied to \hhodp{} discretizations.
Even if a complete analysis and comparison of the parallel performance of nonconforming discretizations is outside the scope of the paper,
we ought to show that Additive Schwarz Method (ASM) preconditioners are an effective means of achieving satisfactory parallel efficiency.
We consider the finest grid of the pyramidal elements mesh sequence (counting of 24k elements) and a \hhodp{} scheme with $k=5$.
Static condensation acts on the sole velocity unknowns (\hhodp{} \vcond), as described in \eqref{eq:ST.1}.
The multilevel solver strategy is the same employed in serial computations for $k=6$,
but smoother preconditioners are suitably designed, as outlined in what follows.

The parallel implementation is based on the distributed memory paradigm and requires to partition the computational mesh in several subdomains. 
In case of HHO methods, not only the mesh but also the mesh skeleton needs to be partitioned: as a result, each mesh entity (element or face) belongs to one and only subdomain.
Each subdomain is assigned to a different computing unit that performs matrix assembly for the \emph{local} mesh elements pertaining to the subdomain. 
Mesh partitioning directly reflects into matrix partitioning in the sense that all entries of the matrix rows 
(PETSc matrix implementation is row-major) pertaining to local mesh entities are allocated and stored in local memory. 
Once matrix assembly is completed, the linear system is approximately solved in each subdomain. 
Depending on the preconditioner strategy, the solver performance might degrade increasing the number of subdomains, see \eg~\cite{SmithDD96}. 

A commonly used ASM preconditioner strategy for DG discretizations consists in employing
an ILU decomposition in each subdomain matrix suitably extended to include the matrix rows of \emph{ghost} elements, 
that is, neighbors of local mesh elements that pertain to a different subdomain.
This implies that the local matrix is extended to encompass the stencil of the DG discretizations, 
see \cite{Franciolini2020} for additional details.
We consider a similar strategy for HHO discretizations: each subdomain matrix is extended to
include the matrix rows of ghost faces, that is, faces of the local mesh elements that pertain to a different subdomain. 
Interestingly, even if the resulting local matrix does not encompass the stencil of the HHO discretization,
mass conservation defect takes into account all element's faces.

As a result of the ASM described above, the amount of \emph{overlap} between subdomain matrices, 
i.e., the number of matrix entries that are repeated in more than one subdomain, is smaller for HHO than for DG. 
Consider, for example, two subdomains sharing a face: if the face is local for subdomain $A$, it is a ghost face for subdomain $B$ and vice-versa.
Accordingly, only one of the two subdomain matrices is extended for HHO discretizations.
As opposite, since each of the two mesh elements sharing the face has a ghost neighbor, both subdomain matrices are extended for DG discretizations.

Scalability is measured on an AMD EPYC cluster of four nodes and 256 cores, increasing the number of execution units from 16 to 256:
in particular we consider a total of five steps doubling the number of execution units at each step.
Notice that, when running on 256 subdomains, each subdomain counts of approximately 96 local elements.
The results reported in Table \ref{tab::HHOdpParPerf} confirm that the ASM preconditioner strategy provides satisfactory parallel performance:
the number of outer FGMRES iterations is uniform while increasing the number of execution units, and only a mild increase in the iteration count is observed for the ASM preconditioned GMRES solvers on the coarse level.
The efficiency parameter (last column in Table \ref{tab::HHOdpParPerf}) measures strong scalability: 100\% efficiency with $N$ execution units
would imply a $N/16$ fold reduction of total computation time with respect to the baseline computation performed with 16 execution units. 

\begin{table}
\centering
\begin{tabular}{ c | c c  c c | c c | c c c | c }
  \toprule
  Ex.Units &\multicolumn{4}{c|}{error in $L^2$ norm} & \multicolumn{2}{c|}{ITs} & \multicolumn{3}{c|}{CPU time} & Eff \\
  \midrule
	&  $\uVec_h$    & $G\uVec_h$   & $p_h$ & $D\uVec_h $& ITs & ITs$_L$& Sol. & Ass. & Tot. & \\
  \midrule
	\rowcolor{LightCyan} 16 &  2.1e-12 & 2.56e-10 & 1.82e-10 & 2.64e-09     & 12 & 132 & 157 & 78 & 235 & -  \\      
	\rowcolor{LightCyan}  32&   3.33e-12 & 2.58e-10 & 1.88e-10 & 2.86e-09  & 12 & 124 & 82.4 & 39.4 & 122 & 96  \\     
	\rowcolor{LightCyan}  64&   9.78e-13 & 2.54e-10 & 1.8e-10 & 2.56e-09   & 13 & 183 & 40.4 & 19.6 & 60 & 98  \\     
	\rowcolor{LightCyan}    128&  1.12e-12 & 2.55e-10 & 1.8e-10 & 2.62e-09  & 13 & 192 & 19.8 & 10.1 & 29.9 & 98  \\    
	\rowcolor{LightCyan}  256&   1.1e-12 & 2.55e-10 & 1.8e-10 & 2.58e-09    & 13 & 186 & 10.7 & 5.11 & 15.8 & 93  \\
  \bottomrule
\end{tabular}
\caption{Parallel performance of $p$-multilevel solution strategies applied to \hhodp{} discretizations with $k=5$. \label{tab::HHOdpParPerf}}
\end{table}

\section{Conclusions}\label{sec:conclusion}

The multilevel $V$-cycle iteration based on $p$-coarsened operators and ILU preconditioned Krylov smoothers 
is an effective solution strategy for high-order accurate HHO discretizations of the Stokes equations. 
The global linear system resulting from the spatial HHO discretization can be solved up to machine precision 
in a reasonable amount of $V$-cycle preconditioned FGMRES iterations (less than 20).
This is remarkable considering that severely graded mesh sequences have been tackled in both 2D and 3D.

Comparing $p$-multilevel solvers for HHO and DG discretizations based on FGMRES iteration count, we can conclude 
that the former are more robust than the latter with respect to both the meshsize and the polynomial degree.
When standard $h$-refined mesh sequences are considered,
HHO formulations show uniform convergence with respect to the meshsize, irrespectively of the considered polynomial degree.
On graded $h$-refined mesh sequences, the iteration count increases over finer meshes, more severely so for DG discretizations.
Similarly, when doubling the polynomial degree (passing from $k=3$ to $k=6$) for a fixed meshsize,
we observe that the iteration count is more stable for HHO schemes.

Since code ruse and code optimization are still possible
(note that HHO implementation is more recent and probably less optimized), we avoid drawing conclusions regarding computation times. 
Nevertheless, the following observations suggest that $p$-multilevel solution strategies are a compelling choice in case of HHO formulations:
\begin{itemize}
\item HHO has a clear advantage over DG both in terms of matrix dimension and number of non-zero entries when the polynomial degree is sufficiently high; 
\item $p$-multilevel solvers for HHO show better solver robustness with respect to the polynomial degree. 
\end{itemize}

\section*{Acknowledgements}
Daniele Di Pietro acknowledges the support of \emph{Agence Nationale de la Recherche} grant fast4hho (ANR-17-CE23-0019).


\bibliographystyle{plain}
\bibliography{pMGhho}

\end{document}